\documentclass[aap,MSNbibl,seceqn,dvips]{arximspdf}

\doi{10.1214/12-AAP897} %kopijuoti is PTS
\volume{23}
\issue{6}
\pubyear{2013}
\firstpage{2161}
\lastpage{2211}

\makeatletter

\newcommand{\rrvert}{\vert}
\newcommand{\llvert}{\vert}
\newcommand{\eqref}[1]{(\ref{#1})}

\newcommand{\eps}{\varepsilon}
\newcommand{\yy}{\mathbf{y}}
\newcommand{\vv}{\mathbf{v}}
\newcommand{\tX}{\tilde{\mathcal{X}}}
\newcommand{\I}{\mathcal{I}}

\newtheorem{theo}{Theorem}[section]
\newtheorem{coro}{Corollary}[section]

\newtheorem{lemm}{Lemma}[section]
\newproclaim{defn}{Definition}[section]
\newproclaim{assu}{Assumption}[section]
\newproclaim{rema}{Remark}[section]
\newproclaim{remarks}{Remarks}

\def\N{\mathbb{N}}
\def\E{\mathbb{E}}
\def\S{\mathcal{S}}
\def\0{\mathbf{0}}
\def\Z{\mathbb{Z}}
\def\R{\mathbb{R}}
\def\dist{\operatorname{dist}}
\def\B{{B}}
\def\T{\mathcal{T}}
\def\epsilon{\varepsilon}
\renewcommand{\E}{\mathbb{E}}
\newcommand{\tH}{\tilde{\mathcal{H}}}

\newcommand{\Po}{\mathcal{P}}
\newcommand{\X}{\mathcal{X}}
\def\la{{\lambda}}
\newcommand{\M}{\mathcal{M}}
\newcommand{\Y}{\mathcal{Y}}
\newcommand{\U}{\mathcal{U}}
\renewcommand{\H}{\mathcal{H}}
\renewcommand{\P}{\mathcal{P}}
\newcommand{\A}{\mathcal{A}}
\newcommand{\Var}{\operatorname{Var}}

\newcommand{\card}{\operatorname{card}}
\newcommand{\diam}{\operatorname{diam}}
\newcommand{\K}{\mathcal{K}}
\newcommand{\F}{\mathcal{F}}

\newcommand{\tk}{\kappa}

\newcommand{\tm}{{\tilde{m}}}
\newcommand{\NN}{\mathcal{N}}

\newcommand{\1}{\mathbf{1}}

\def\R{\mathbb{R}}

\def\de{{\delta}}
\def\A{\mathcal{A}}

\def\la{\lambda}

\def\ka{\tilde{\kappa}}

\def\tka{{\tk}}

\def\VRdelta{\beta}
\def\Cl{\mathcal{C}} % clique number
\newcommand{\NND}{N} % nearest neighbour distance
\newcommand{\tN}{\tilde{N}} % modified nearest neighbour distance

\makeatother

\begin{document}
\begin{frontmatter}

\title{Limit theory for point processes in manifolds}
\runtitle{Limit theory for point processes in manifolds}

\begin{aug}
\author[A]{\fnms{Mathew D.} \snm{Penrose}\thanksref{t1}}
\and
\author[B]{\fnms{J.~E.} \snm{Yukich}\corref{}\thanksref{t2}\ead[label=e2]{jey0@lehigh.edu}}
\thankstext{t1}{Supported in part by the
Alexander von Humboldt Foundation through a Friedrich
Wilhelm Bessel Research Award.}
\thankstext{t2}{Supported in part by NSF Grant DMS-08-05570.}
\runauthor{M.~D. Penrose and J.~E. Yukich}
\affiliation{University of Bath and Lehigh University}
\address[A]{Department of Mathematical Sciences\\
University of Bath\\
Claverton Down, Bath BA2 7AY\\
United Kingdom} %adresu isvedimo komanda gale!
\address[B]{Department of Mathematics\\
Lehigh University\\
Bethlehem, Pennsylvania 18015\\
USA}
\end{aug}

% HISTORY:
\received{\smonth{4} \syear{2011}}
\revised{\smonth{4} \syear{2012}}

% ABSTRACT
%
\begin{abstract}
Let $Y_i, i \geq1$, be i.i.d. random variables
having values in an $m$-dimensional manifold $\M\subset\R^d$ and
consider sums
$\sum_{i=1}^n \xi(n^{1/m}Y_i,\break  \{n^{1/m}Y_j\}_{j=1}^n)$,
where $\xi$ is a real
valued function defined on pairs $(y, \Y)$, with $y \in\R^d$
and $\Y\subset\R^d$ locally finite. Subject to $\xi$ satisfying a weak
spatial dependence and continuity condition, we show that such sums
satisfy weak laws of large numbers, variance asymptotics and
central limit theorems. We show that the limit behavior is
controlled by the value of $\xi$ on homogeneous Poisson point
processes on
$m$-dimensional hyperplanes tangent to $\M$.
We apply the general results to
establish the limit theory of dimension and
volume content estimators, R\'enyi and
Shannon entropy estimators and clique counts in the
Vietoris--Rips complex on $\{Y_i\}_{i=1}^n$.
\end{abstract}

% KEYWORDS
\begin{keyword}[class=AMS]
\kwd[Primary ]{60F05}
\kwd[; secondary ]{60D05}
\end{keyword}
\begin{keyword}
\kwd{Manifolds}
\kwd{dimension estimators}
\kwd{entropy estimators}
\kwd{Vietoris--Rips complex}
\kwd{clique counts}
\end{keyword}

\end{frontmatter}

%s1 #&#
\section{Introduction}\label{INTRO}

There has been recent interest in the statistical, topological and
geometric properties of high-dimensional nonlinear data sets.
Typically the data sets may be modeled as realizations of i.i.d.
random variables $\{Y_i\}_{i=1}^n$
having support on an unknown nonlinear manifold $\M$ embedded in $\R^d$.
Given a sample $\{Y_i\}_{i=1}^n$, whose pairwise distances are
given, but whose coordinate representation is not, can one determine
geometric characteristics of the manifold, including its intrinsic
dimension and volume content?
Can one
recover global properties of the distribution of $\{Y_i\}_{i=1}^n$
such as its intrinsic entropy? These properties, as well as graph
theoretic functionals such as clique counts in the Vietoris--Rips
complex
may be studied via the statistics of the form
%
%e1.1 #&#
\begin{equation}
\label{lstat} \sum_{i
=1}^n \xi
\bigl(Y_i, \{Y_j\}_{j=1}^n\bigr),
\end{equation}
where $\xi( \cdot, \cdot)$ is a
real-valued measurable function defined on pairs $(y, \Y)$, where $y
\in\Y$ and $\Y\subset\R^d$ is locally finite, with
$\xi(y,\Y)$ locally determined in some sense.

Our goal is to establish the dependency between the large-$n$
behavior of the statistics (\ref{lstat}), the underlying point
density $\tk$ of the $\{Y_i\}_{i=1}^n$ and the manifold $\M$. For
$\M= \R^d$ there is a large literature describing limit theorems
for (\ref{lstat}) \cite{BY2,PeEJP,PeBer,PY1,PY4,Yubk} whereas
when $\M\neq\R^d$, there is a relative dearth of results, although
spatial data generated by a
curved surface embedded in three dimensional space are arguably more
natural than those involving
data generated by flat surfaces.
This paper partly redresses this situation under reasonably general
conditions on $\tk$ and $\M$.

In Section \ref{secresults} we present laws of large numbers and
central limit
theorems for (i)~the Levina--Bickel dimension estimator for data
$\{Y_i\}_{i=1}^n$ supported on a manifold, (ii) R\'enyi %Tsallis,
and Shannon entropy estimators for $\{Y_i\}_{i=1}^n$,
(iii)~volume estimators for the support of
$\{Y_i\}_{i=1}^n$ and (iv)
the order $k$ clique count in the
Vietoris--Rips complex on $\{Y_i\}_{i=1}^n$. %, and (iv) laws of
The asymptotic normality results for dimension and entropy
estimators appear to be new even in the setting of linear
manifolds. The mean and variance asymptotics for the Levina--Bickel
dimension estimator %and the component count in the nearest
depend only on the dimension of $\M$ and are invariant with respect
to $\tk$, whereas for most of the other functionals considered here,
the mean and variance asymptotics explicitly depend on $\tk$ via the
integral $\int_{\M} (\tk(y))^p \,dy$ for some $p \in\R$.

We shall derive these results
%show that these results can all be deduced
from
general theorems governing the limit theory of $\sum_{i =1}^n
\xi(n^{1/m}Y_i, \{n^{1/m}Y_j\}_{j=1}^n)$, where $m$ is the intrinsic
dimension of $\M$, $n^{1/m}$ is a dilation factor and $\xi$ belongs to
a general class of translation invariant
functionals %$\xi(y,\Y)$
that are determined by
the locations of either the $k$ nearest neighbors of $y$ for some
fixed $k \in\N$, or the points of $\Y$ within some fixed distance
of $y$. For $\xi$ locally
determined in this way, then since the manifolds are themselves
local, one might expect as the number of sample points increases,
that the local contribution of each $\xi$ at $y \in\M$ converges in
distribution to its linearized version on the tangent space to $y$.
In other words, for locally determined and translation invariant
$\xi$, one might expect that the large $n$ behavior of
$\xi(n^{1/m}y, \{n^{1/m}Y_j\}_{j=1}^n)= \xi(\0, \{n^{1/m}(Y_j-
y)\}_{j=1}^n) $ is controlled by the behavior of $\xi(\0,
\H_{\kappa(y)})$, with $\H_{\kappa(y)}$ a homogeneous Poisson point
process of intensity $\kappa(y)$ on a \emph{tangent hyperplane} of
Euclidean dimension
$m$ %, that is to say the Grassmannian manifold $Gr_m(d)$ of
%$m$-dimensional linear subspaces of $\R^d$
(here and elsewhere $\0$
denotes a point at the origin of $\R^k$). Subject to moment
conditions on $\xi$, this is indeed the case, as shown by the
general results of Section \ref{gensec}.
%Our general results are presented in Section \ref{gensec}.
%The limit theory makes use of the fact that for such $\xi$,
%the large $n$ behavior with respect to the
%shifted and dilated point sets $ \{n^{1/m}(Y_j- Y_1)\}_{j=1}^n $
%is controlled by the behavior of $\xi$ on homogeneous Poisson point
%processes on \emph{tangent hyperplanes} of Euclidean dimension $m$,
%that is to say %\Comment{Grassmannian defined. JY}
%the Grassmannian manifold $Gr_m(d)$ of $m$-dimensional linear
%subspaces of $\R^d$.
The locally  defined behavior of~$\xi$, quantified in terms of
dependency graphs involving radii of stabilization of~$\xi$, yields
central limit theorems via Stein's method.

In fact, our methods should work in still greater generality. Most of
the examples considered in \cite{BY2,PeEJP,PeBer,PY1,PY4,Yubk}
have $\xi$ stabilizing, that is, they are locally determined in some
sense, and
it should be possible to adapt our methods to most of these
examples. For example, we anticipate that our methods
can be extended to
establish the limit theory for the total
edge length and other stabilizing
functionals of the Delaunay and Voronoi graphs on random
point sets in manifolds.
We also expect that our methods extend to give the limit theory of
statistics of germ-grain models, coverage processes and random
sequential adsorption models generated by data on manifolds.
Moreover,
we anticipate that the theory presented
here can be modified to establish
limit theorems for
generalized spacing statistics based on $k$ nearest
neighbor distances for random points in a manifold, including
estimators of relative entropy such as those considered in
\cite{BPY}, although these involve consideration of nontranslation invariant~$\xi$ so they
do not automatically fall within
the scope of this paper.

In many examples the
functional $\xi(y,\Y)$ is determined by
inter-point distances in the vicinity of $y$. For $\M$ an arbitrary
Riemannian manifold, it may be possible to derive similar limit
results for such $\xi$ using the geodesic distance rather than the
extrinsic distance in $\R^d$ (by the Nash embedding theorem, such
$\M$ can always be embedded into some $\R^d$). However, this lies
beyond the scope of the present paper.

%s2 #&#
\section{Stochastic functionals on manifolds}
\label{secresults}

%s2.1 #&#
\subsection{Terminology and definitions}
\label{sectermin}

For $k \in\N$, let $\| \cdot\|$ be the Euclidean norm in~$\R^k$.
Recall that $\0$ denotes a point at the origin of $\R^k$.
For $r \in(0, \infty)$ and $z \in
\R^k$, let $B_r(z):= \{y \in\R^k\dvtx \|y-z\| \leq r\}$. Given $F
\subset\R^k$, and $y \in\R^k$, $a >0$, set $y + F :=\{y + z\dvtx z
\in F\}$ and
$a F := \{az \dvtx z \in F \}$.
If $F$ is locally finite, let $\card(F)$ denote the cardinality
(number of
elements) of $F$.
If also $y \in\R^k$ and $j \in\Z^+ := \{0,1,2,\ldots\}$,
then let $\NND_j(y,F)$ be
the Euclidean distance between $y$ and its
$j$th nearest neighbor in $F \setminus\{y\}$, that is,
%
%e2.1 #&#
\begin{equation}\label{NNDdef}
\NND_j(y,F): = \inf\bigl\{r \geq0\dvtx \card\bigl(F \cap
B_r(y) \setminus\{y\} \bigr) \geq j \bigr\}
\end{equation}
with the infimum of the empty set
taken to be $+\infty$. In particular, $N_0(y,F)=0$. Let $\Phi$ be
the distribution function for the standard normal random variable
$\NN(0,1)$, and let  $\stackrel{{P}}{\longrightarrow}$  denote
convergence in probability.
For  $\sigma>0$, let  $\NN(0,\sigma^2)$  denote the random variable $\sigma\NN(0,1)$.

Let $m \in\N$
and $d \in\N$ with $m \leq d$. A nonempty subset $\M$ of $\R^d$,
endowed with the subset topology,
is called an $m$-dimensional $C^1$
submanifold of $\R^d$ if for each $y \in\M$ there exists an open
subset
$U$
of $\R^m$ and a continuously differentiable
injection
$g $
from
$U $ to $\R^d$, such that (i) $y \in g(U) \subseteq\M$,
and (ii) $g$ is an open map from $U$ to $\M$, and
(iii) the linear map $ g'(u)$ has full rank for all $ u \in U$; see,
for example, Theorem 2.1.2(v) of \cite{BG}. The pair $(U,g)$ is called a
\emph{chart}. Let ${\mathbb{M}}:= \mathbb{M}(m,d)$ denote the class of all
$m$-dimensional $C^1$ submanifolds of $\R^d$ which are also closed
subsets of $\R^d$.

Given $\M\in\mathbb{M}$, using a routine compactness argument we can
choose an index set $\I\subset\N$, and a set
$\{(y_i,\delta_i,U_i,g_i),i \in{\mathcal I} \}$ of ordered quadruples
with $y_i \in\M$, $\delta_i \in(0, \infty)$ and
$(U_i,g_i)$ a chart for each $i$, such that
(i)~$\M\cap B_{3 \delta_i}(y_i) \subset g_{y_i}(U_i)$ for each
$i$, and (ii) $\M\subset\bigcup_{i \in{\mathcal I}} B_{\delta_i}(y_i)$.

We refer to $((U_i,g_i), i \in{\mathcal I})$ as an \emph{atlas} for
$\M$. Given such an atlas, we can find a \emph{partition of unity}
$\{\psi_i\}$ subordinate to the atlas, that is, a collection of
functions $(\psi_i,i \in{\mathcal I})$ from $\M$ to $[0,1]$, such that
$\sum_{i \in{\mathcal I}} \psi_i(y) =1$ for all $y \in\M$, and such
that for each $i$, $ \psi_i(y) = 0 $ for $y \notin g_i(U_i)$, and
$\psi_i \circ g_i$ is a measurable function on~$U_i$.
The more common definition of a partition of unity has some extra
differentiability conditions on $\psi_i$ but these are not needed
here. With our more relaxed definition, the existence of a partition
of unity is completely elementary to prove.

Given $i \in{\mathcal I} $ and $x \in U_i$, let $D_{g_i}(x):= \sqrt{
\det(J_{g_i}(x))^T(J_{g_i}(x)) }$, with $J_{g_i}$ standing for the
Jacobian of ${g_i}$.
For bounded measurable $h\dvtx \M\to\R$, the integral $\int_{\M} h(y)
\,dy$ is defined by
%
%e2.2 #&#
\begin{equation}
\label{int-form}\int_{\M} h(y) \,dy = \sum
_{i \in{\mathcal I}} \int_{U_i} h\bigl(g_i(x)
\bigr) \psi_i \bigl(g_i(x)\bigr) D_{g_i}(x)\,dx,
\end{equation}
which is
well-defined in the sense that it does not depend on the choice of
atlas or the partition of unity. Equation (\ref{int-form}) is discussed
further at (\ref{0508d}) below.

Given any manifold $\M\in\mathbb{M}$ and nonempty $\K\subset\M$, the relative
interior of $\K$ consists of all those $y \in\K$ such that $y $ has
a neighborhood in $\M$ that is contained in $\K$. The \emph{boundary}
of $\K$ is the set of all other $y \in\K$ (possibly empty).
Also, we
%Given $\M\in\m$ and $\K\subset\M$,
% $\kappa\in\p_b(\M)$,
set
$
\diam( \K) := \sup\{\|x - y\|\dvtx x,y \in\K\}$
(possibly infinite).
We say
%$\K\subset\M$ is \emph{locally conic} if
$\K$ is \emph{locally conic} if
$0 < \diam(\K) < \infty$ and
%for all $\tk\in\p_c(\M)$
% there is $L:= L(m, \kappa) \in(0,
% $r \in(0, \diam(\K)]$ and all $w \in\K$ we have
%
%e2.3 #&#
\begin{equation}
\label{lconic} %\inf_{r \in(0, 1],w \in\K}
\inf\biggl\{ r^{-m} \int
_{B_r(w) \cap\K} \,dy \dvtx r \in\bigl(0, \diam(K)\bigr],w \in\K\biggr\} >0.
\end{equation}
We say
%$\K\subset\M$
$\K$ is an $m$-dimensional $C^1$
\emph{submanifold-with-boundary} of $\M$ if for all $y$ in the
boundary of $\K$, there exists a choice of chart $(U,g)$ for $\M$
such that $\0 \in U$, and $g(\0)=y$, and $g([0,\infty) \times
\R^{m-1}) = g(U) \cap\K$. This includes the possibility that $\K$
has empty boundary.
If $\K$ is a compact $m$-dimensional $C^1$
submanifold-with-boundary of $\M$ then
it is locally conic; see Remark \ref{bdyrmk}.

Given $\M\in\mathbb{M}$, a \emph{probability density function}
on $\M$ is
a nonnegative scalar field
$\tk$ on $\M$
satisfying
$\int_{\M} \tk(y) \,dy = 1$. Let $\mathbb{P}(\M)$ denote the class of
probability density functions on $\M$. %Given $\tka\in\p(\M)$, let
%$\K(\tka)$ denote the support of $\tka$, i.e. the smallest closed
%set $\K\subset\M$ such that $\int_\K\tka(y) \,dy =1$.
Given $\tka\in\mathbb{P}(\M)$, let $\K(\tka)$ denote the support of
$\tka$, that is, the smallest closed set $\K\subset\M$ such that
$\int_\K\tka(y) \,dy =1$.
Given also
$\rho\in\R$,
define the integral
%
%e2.4 #&#
\begin{equation}\label{Idef}
I_\rho(\tk) := \int_{\K(\tk)} \bigl(\tk(y)
\bigr)^{\rho}\,dy := \int_\M \bigl(\tk(y)
\bigr)^{\rho} \mathbf{1}\bigl\{\K(\tk)\bigr\}(y) \,dy.
\end{equation}

Let ${\mathbb{P}}_b(\M)$ denote the class of bounded probability
density functions
$\tk\in\mathbb{P}(\M)$, such that
$\K(\tk)$ is compact.
Let ${\mathbb{P}}_c(\M)$ denote those probability density functions
$\tk\in
\mathbb{P}_b(\M)$ whose support $\K(\tk)$ is locally conic
%a compact $C^1$ sub-manifold-with-boundary of $\M$,
and which are bounded away from
zero and infinity on their support. The motivation for considering
these classes of probability densities appears in Remark 3 following
Theorem \ref{smainCLT}.

Suppose $\M\in\mathbb{M}$ and $\tka\in\mathbb{P}(\M)$ are given. Let
$Y_1, Y_2,\ldots$ be i.i.d.
random variables
%with support $\M$, having
with probability density function $\tk$ with respect to the
Riemannian volume element $dy$. Define the binomial point process
$\Y_n := \{Y_i\}_{i=1}^n$.
Let $\P_\la$
denote the Poisson point process on $\M$ (and also the associated
counting measure)
having intensity density $\lambda\tk(\cdot)$, that
is,
%
%e2.5 #&#
\begin{equation}
\label{PPP} \E\P_\la(dy) = \la\tk(y)\,dy.
\end{equation}

Recall that $\xi(\cdot, \cdot)$ denotes a measurable function
defined on pairs $(y,\Y)$, where $y \in\Y\subset\R^d$ and $\Y$ is
locally finite.
When
$y \notin\Y$, we write $\xi(y,\Y)$ instead of $\xi(y, \Y\cup\{y\})$;
also, we sometimes write $\Y^y$ for $\Y\cup\{y\}$.

Let $\H$ denote a homogeneous Poisson process of unit intensity in
$\R^m$ with $\R^m$
embedded in $\R^d$ (since $m \leq d$) so that the random variable
$\xi(\0,\H)$ is well defined.
In keeping with
notation from \cite{BY2,PeEJP},
for all such functionals $\xi$
we set
%
%e2.6 #&#
\begin{equation}\label{Videnti}
V^{\xi}:= \E\xi(\0, \H)^2 + \int_{\R^m}
\bigl\{ \E\xi\bigl(\0, \H^{u}\bigr) \xi\bigl(u, \H^{\0}\bigr)
- \bigl(\E\xi(\0, \H)\bigr)^2 \bigr\} \,du %\label{Vdef1}
\end{equation}
and
%
%e2.7 #&#
\begin{equation}
\label{Didenti} \delta^{\xi} := \E\xi(\0, \H) + \int
_{\R^m} \E\bigl[\xi\bigl( \0, \H^{u} \bigr) - \xi(
\0, \H)\bigr] \,du,
\end{equation}
whenever these integrals are defined.

%s2.2 #&#
\subsection{Estimators of intrinsic dimension, manifold learning}
\label{subseclearn}

Given data embedded in a high-dimensional vector space,
a natural problem in manifold learning, signal processing and
statistics is to discover the low-dimensional structure of the data,
namely the intrinsic dimension of the hypersurface containing the
data. Levina and Bickel \cite{LB} propose a dimension estimator
making use of nearest neighbor statistics. Their estimator, which
uses distances between a given sample point and its $k$ nearest
neighbors, estimates the dimension of random variables lying on a
manifold $\M$ of unknown dimension $m$ embedded in $\R^d,  d \geq
m$. Specifically, for all $k = 3,4,\ldots,$ the Levina and Bickel
estimator of the dimension of a finite data cloud $\Y\subset\M\in
{\mathbb{M}}$, is given by
%
%e2.8 #&#
\begin{equation}
\label{estim} \hat{m}_k:= \hat{m}_k(\Y) := \bigl(
\operatorname{card}(\Y)\bigr)^{-1} \sum_{y \in\Y}
\zeta_k(y,\Y),
\end{equation}
where for
all $y \in\Y$ we have
%
%e2.9 #&#
\begin{equation}
\label{estim1} \zeta_k(y,\Y): = \cases{\displaystyle(k - 2) \Biggl( \sum
_{j = 1}^{k -1} \log \frac{\NND_k(y)}{ \NND_j(y)}
\Biggr)^{-1}, &\quad $\mbox{if}  \card(\Y) \geq k+1$, \vspace*{2pt}
\cr
0, & \quad $\mbox{otherwise},$}
\end{equation}
where $\NND_j(y):=\NND_j(y,\Y)$, as
given by (\ref{NNDdef}). For all $\rho> 0$,
we also define
%
%e2.10 #&#
\begin{equation}
\label{estim1r} \zeta_{k,\rho}(y,\Y): = \zeta_k(y,\Y) \mathbf{1}
\bigl\{\NND_k(y) \leq\rho\bigr\},
\end{equation}
using the convention $0 \times\infty=0$ if necessary,
and we put
%
%e2.11 #&#
\begin{equation}
\label{estimr} \hat{m}_{k,\rho}:= \hat{m}_{k,\rho}(\Y) := \bigl(
\operatorname{card}(\Y)\bigr)^{-1} \sum_{y
\in\Y}
\zeta_{k,\rho}(y,\Y).
\end{equation}

Given $\Y_n := \{Y_i\}_{i=1}^n$ as in Section \ref{sectermin},
Levina and Bickel \cite{LB} argue that $\hat{m}_k(\Y_n)$ estimates
the intrinsic dimension of $\M$. Our purpose here is to
substantiate this claim and to provide further distributional
results. The following shows (i) consistency of the dimension
estimator $\hat{m}_k$ over Poisson and binomial samples,
and (ii) a~central limit
theorem for $\hat{m}_{k,\rho}(\Y_n)$, $\rho$ fixed and small,
addressing a question raised by Peter Bickel.
Here
$\de^{\zeta_{k} }$
is
given by
taking $\xi\equiv\zeta_k$ in
(\ref{Videnti}) and (\ref{Didenti}).

%
%th2.1 #&#
\begin{theo} \label{dimLLN} Let $\M\in\mathbb{M}(m,d)$ and let $\tk
\in\mathbb{P}_c(\M)$.
For all $k \geq3$ and $\rho> 0$, we have
%
%e2.12 #&#
\begin{equation}
\label{close} \lim_{n \to\infty} \mathbf{1} \bigl\{\hat{m}_{k,\rho} (
\Y_n) \neq\hat{m}_{k} (\Y_n)\bigr\} = 0  \qquad \mbox{a.s.},
\end{equation}
and if $k \geq4$, then
%
%e2.13 #&#
\begin{equation}
\label{dimestLLN} \hat{m}_k (\Y_n) \stackrel{{P}} {
\longrightarrow}m    \qquad \mbox{as }  n \to\infty,
\end{equation}
while if $k \geq11$, then (\ref{dimestLLN}) holds a.s.
If $\tk$ is a.e. continuous, then there exists $\rho_1 > 0$ such that
if $\rho\in(0, \rho_1)$ and $k \geq7$, then
%
%e2.14 #&#
\begin{equation}
\label{devarlimit} \lim_{n \to\infty} n \Var \bigl[\hat{m}_{k,\rho}(
\Y_n)\bigr] = {\sigma}^2(\zeta_{k}):=
\frac{m^2}{k-3} - \bigl(\de^{\zeta_{k} }\bigr)^2 > 0,
\end{equation}
and also we have as $n \to\infty$ that
%
%e2.15 #&#
\begin{equation}
\label{LBCLT} n^{1/2}\bigl(\hat{m}_{k,\rho}(\Y_n) -
\E\hat{m}_{k,\rho}(\Y_n)\bigr) \stackrel{{\mathcal D}} {
\longrightarrow} \NN\bigl(0, {\sigma}^2(\zeta_{k})\bigr).
\end{equation}
\end{theo}

\begin{remarks*}
(i) (\emph{Counterexamples}.) Without further
conditions on $\M\in\mathbb{M}$, in general
$\hat{m}_{k}$
(as opposed to $\hat{m}_{k,\rho}$)
might not satisfy the variance asymptotics (\ref{devarlimit}) or the
central limit theorem (\ref{LBCLT}).
To see this, suppose $m = 1$ and $d = 3$, and suppose
$\M$ is a compact $1$-manifold which includes a segment $S$ on the
$z$-axis, say from $z=0$ to $z=1$, as well as an arc of the unit
circle in the $(x,y)$ plane. If $Y_i$ are i.i.d. with the uniform
measure on $\M$, then there is a positive probability that $Y_1 \in
S$ and that its $k$ nearest neighbors are all on the arc. In this
case the $\NND_j(Y_1, \Y_n), j = 1,\ldots,k-1,$ all coincide and none
of the moments of $\hat{m}_k (\Y_n)$ exist.

For a more general counterexample along similar lines,
fix $m$ and put \mbox{$d = 2(m + 1)$}. Let $S$ be the set of unit vectors
in $\R^d$ and
let $S_1$ be those $x \in S$ such that the first $m + 1$ coordinates
are zero and let $S_2$ be those $x \in S$ such that the last $m + 1$
coordinates are zero. Then $\M:= S_1 \cup S_2$ is an
$m$-dimensional manifold in $\mathbb{M}$. If $Y_i$ are i.i.d. with the
uniform volume measure on $\M$, then there is a positive probability
that $Y_1 \in S_1$ and that its $k$ nearest neighbors all belong to
$S_2$. In this case $\NND_j(Y_1, \Y_n) = 2$ for all $j = 1,\ldots,k-1$,
showing that $\zeta_k(Y_1,\Y_n)$ is infinite with positive
probability, and therefore none of the moments of $\hat{m}_k (\Y_n)$
exist.\vspace*{-6pt}
\begin{longlist}[(iii)]
\item[(ii) ](\emph{On the constant $\rho_1$}.) The constant $\rho_1$,
loosely speaking, reflects the maximum amount of curvature of the
manifold $\M$. The larger this is, the smaller $\rho_1$ has to be
taken. See Lemma \ref{claim0} below.

\item[(iii)] (\emph{Relation to previous work}.) Levina and Bickel
(Section 3.1 of \cite{BB}) argue on heuristic grounds that $\Var
[\hat{m}_k(\Y_n)] = O(n^{-1})$ whenever the $Y_i$ are the image
under a sufficiently smooth map $g$ of random variables having a
smooth density, but the last counterexample shows that this bound
is not true in general.

Chatterjee \cite{Ch} provides a rate of normal approximation for
$\hat{m}_k(\Y_n)$ for all $k > 9$ whenever $\M$ is a ``nice''
manifold and under minimal assumptions on the distribution of
$Y_i$. His rates are with respect to the
Kantorovich--Wasserstein distance and are of order
$n^{(9-k)/(2k -2)}$,
subject to the validity of $\Var[\hat{m}_k(\Y_n)] = \Theta(n^{-1})$.

Bickel and Yan (Theorems 1 and 3 of Section 4 of \cite{BY})
establish a central limit theorem for $\hat{m}_k(\Y_n)$ for $\M=
\R^d$. The methods of Bickel and Breiman \cite{BB} or \cite{BY2,PY1,PeEJP} could be used to establish the asymptotic normality of
$\hat{m}_k(\Y_n)$ if the $Y_i$ had a density with respect to
Lebesgue measure on $\R^d$. These methods do not appear applicable
in the present situation. Under strong assumptions on~$\M$,
Yukich~\cite{YuBlau} outlines an approach giving a rate of normal
approximation for $\hat{m}_k(\P_\la)$, but does not provide
consistency results or variance asymptotics for either
$\hat{m}_k(\P_\la)$ or $\hat{m}_k(\Y_n)$. % nor does it establish

\item[(iv)] (\emph{Limits are dimension dependent only}.) The mean and
variance asymptotics (\ref{dimestLLN}) and (\ref{devarlimit}) are
invariant with respect to $\tk$ and depend only on $\operatorname{dim}
(\M)$, and (in the case of variance) the parameter $k$. These
results appear to be new even for $\M= \R^d$.

\item[(v)] (\emph{Variance asymptotics}.)
For Poisson samples, a similar result to Theorem~\ref{dimLLN}
holds (see Theorem \ref{PoCLT} below) but
with the limiting variance $\sigma^2(\zeta_k)$
modified to simply $m^2/(k-3)$.
Thus at least for Poisson samples, the
limiting variance of the dimension estimator $\zeta_k$ decreases with
increasing $k$.
\end{longlist}
\end{remarks*}

%s2.3 #&#
\subsection{Estimators of intrinsic entropy and volume content}

\textit{R\'enyi entropies}. %Let $\{C_j\}$ be a partition of a linear
Let $Y_1$ be as in Section \ref{sectermin}.
Given $\rho
> 0$ with $ \rho\neq1$, the R\'enyi $\rho$-entropy \cite{Re} of
$Y_1$, denoted
$H_\rho^*(\tk) $,
and
the closely related Tsallis entropy (or Havrda and Charv\'at entropy
\cite{HC}), denoted $H_\rho(\tk) $,
are given, respectively, by
\[
H_\rho^*(\tk) := (1 - \rho)^{-1} \log I_{\rho}(\tk)
;  \qquad     H_\rho(\tk) := (\rho- 1)^{-1} \bigl(1-
I_{\rho}(\tk) \bigr),
\]
where $I_\rho(\tk):= \int_{\K(\kappa)} (\kappa(y))^{\rho}\,dy$ is
at (\ref{Idef}).
When $\rho$ tends to $1$, $H_\rho^*(\tk)$ and $H_\rho(\tk)$ tend
to
the Shannon differential entropy
%
%e2.16 #&#
\begin{equation}
\label{Shan}H_1(\tk) := - \int_{\M} \tk(y) \log
\bigl(\tk(y)\bigr)\,dy.
\end{equation}
R\'enyi and
Tsallis entropies are used in the study of nonlinear Fokker--Planck
equations, fractal random walks, parameter estimation in
semi-parametric modeling, and data compression; see \cite{CHH} and
the introduction of \cite{LPS} for details and references. The
gradient $\lim_{\rho\to1} (dH^*/d\rho)$ equals $(-1/2) \Var[\log
\tka(Y_1)]$, a measure of the shape of the distribution \cite{LPS,So}
which also appears in the statement of Theorem \ref{ShCLT} below.

A problem of interest is to estimate the R\'enyi and Tsallis
entropies given only the sample $\{Y_i\}_{i=1}^n$ and their
pairwise distances.
Here we show consistency results, variance asymptotics and central
limit theorems for nonparametric estimators of $I_\rho(\tk) $.
The authors of the papers \cite{CHH,JK,LPS,Wa,PY6} consider
estimators of $I_\rho(\tk)$ in terms of the $k$ nearest-neighbor
graph; here we restrict to $k = 1$, but this is for presentational
purposes only. The approach taken here also yields consistent
estimators of
volume content.

Recall that $\NND_1(y, \Y)$ is the distance between $y$
and its nearest neighbor in $\Y$. For all
$\alpha\in(-\infty, \infty)$ and
finite $\Y\subset\R^d$
put
\[
R^{\alpha}(\Y):= \sum_{y \in\Y}
\NND_1(y, \Y)^\alpha.
\]

For $r \in(0,\infty)$,
define the critical moment $r_c (\tk) \in[0, \infty]$ by
%
%e2.17 #&#
\begin{equation}
\label{rc} r_c(\tk) := \sup\bigl\{ r \geq0\dvtx \E\|Y_1
\|^r < \infty\bigr\} .
\end{equation}
The next result spells out conditions under which
$n^{-1}R^{\alpha}(n^{1/m} \Y_n)$ consistently estimates a scalar
multiple of $I_{1 - \alpha/m}(\tk)$ in the $L^q$ sense.
Let $\omega_m:= \pi^{m/2}[\Gamma(1 + m/2)]^{-1}$ be
the volume of the unit radius $m$-dimensional ball.

%
%th2.2 #&#
\begin{theo} \label{NNGCLT1}
If $\tk\in\mathbb{P}_c(\M)$ and $\alpha\in(0,\infty)$, then
%
%e2.18 #&#
\begin{equation}
\label{NNWLLN1}\quad  n^{-1} R^{\alpha}\bigl(n^{1/m}
\Y_n\bigr) \to \omega_m^{-\alpha/m} \Gamma\biggl(1 +
\frac{\alpha}{ m}\biggr) I_{1 - \alpha/m}(\tk)   \qquad  \mbox{as }   n \to\infty
\end{equation}
with both $L^2$ and a.s. convergence.
If instead $\kappa\in\mathbb{P}(\M) $ is bounded,
and $\alpha\in(-m/q,0)$
for $q=1$ or $q=2$,
then (\ref{NNWLLN1}) holds with $L^q$ convergence.
\end{theo}

Putting $m = \alpha$,
we obtain consistent estimators of $I_0(\kappa)$, that is, the
$m$-dimensional
content of the support of $\kappa$.\eject

%
%co2.1 #&#
\begin{coro} If $\tk\in\mathbb{P}_c(\M)$, then
$
\omega_m R^{m}(\Y_n) \to
I_0( \tk)
$
as $ n \to\infty$,
with both $L^2$ and a.s. convergence.
\end{coro}

We now state variance asymptotics and a central limit theorem for\break
$R^{\alpha}(n^{1/m} \Y_n)$. Define $V^{\NND_1^\alpha}$ and
$\de^{\NND_1^\alpha}$ by taking
$\xi\equiv\NND_1^\alpha$ in (\ref{Videnti}) and (\ref{Didenti}).

%
%th2.3 #&#
\begin{theo} \label{NNGCLT3} Suppose $\tk\in\mathbb{P}_c(\M)$
is a.e. continuous, and
$\alpha\in(-m/2,\break0) \cup(0, \infty)$. Then
%
%e2.19 #&#
\begin{eqnarray}
\label{binomNNvar}
\quad\lim_{n \to\infty} n^{-1} \Var
\bigl[R^\alpha\bigl(n^{1/m} \Y_n\bigr) \bigr] &= &{
\sigma}^2\bigl(N_1^\alpha,\tka\bigr)
\nonumber
\\[-8pt]
\\[-8pt]
\nonumber
&:=&
V^{\NND_1^\alpha} I_{1
- 2\alpha/m}(\tk) - \bigl( \delta^{\NND_1^\alpha}
I_{1 -
\alpha/m}(\tk) \bigr)^2 %> 0
\end{eqnarray}
and, as $n \to\infty$
%
%e2.20 #&#
\begin{equation}
\label{binomNNCLT} n^{-1/2} \bigl(R^{\alpha}\bigl(n^{1/m}
\Y_n\bigr) - \E R^{\alpha}\bigl(n^{1/m}
\Y_n\bigr)\bigr) \stackrel{{\mathcal D}} {\longrightarrow }\NN\bigl(0,
{\sigma}^2\bigl(N^\alpha_1,\tka\bigr)\bigr).
\end{equation}
Also, $\sigma^2(N_1^\alpha,\tka)>0$ and
%
%e2.21 #&#
\begin{equation}\label{NNde}
\delta^{N_1^\alpha} = (1 - \alpha/m) \omega_m^{-\alpha/m}
\Gamma\biggl(1+ \frac{\alpha}{m}\biggr) .
\end{equation}
\end{theo}

\begin{remarks*}
(i) [\emph{Comparison of \textup{(}\ref{NNWLLN1}\textup{).} with previous
work}.] Limit (\ref{NNWLLN1}) extends the law of large numbers
limit theory for entropy estimators developed by Costa and Hero on
manifolds \cite{CHH}, who restrict to $\alpha\in(0,m)$
and to compact manifolds. They extend the
consistency results of Leonenko et al. (Theorem~3.2 of~\cite{LPS}),
Theorem~2.1 of Wade~\cite{Wa}, Theorems~2.1--2.3 of~\cite{PY6} and
Theorem~2.4 of~\cite{PY4}, all of which restrict to $\M= \R^m$. In
\cite{PY6} it is shown that if $\M= \R^m$, then~(\ref{NNWLLN1})
holds whenever $\alpha\in(0,m/q)$, $I_{1 - \alpha/m}(\tk) <
\infty$ and
$r_c (\tk) > q\alpha m/(m - q \alpha)$.
For $\M= \R^m $ and
% $\tk\in\p_c(\M)$,
$\tk$ supported by an $m$-dimensional submanifold-with-boundary
of~$\M$,
under a Lipschitz assumption on $\kappa$ (along with an assumption on its
gradient), Liiti\"{a}inen et al. \cite{LLC} develop a closed form
expansion for the moments of $R^\alpha$, $\alpha>0$.

(ii) [\emph{Comparison of \textup{(}\ref{binomNNvar}\textup{)} and \textup{(}\ref{binomNNCLT}\textup{)}
with previous work}.] Existing central  limit theorems and variance
asymptotics for the entropy estimators  $R^\alpha$
and the volume content estimator
$\omega_m R^m$ (e.g., Theorem 6.1 of \cite{PY1}) assume that  $\M=
\R^m$ and $\alpha>0$.  Theorem \ref{NNGCLT3}  allows us to relax
both assumptions.
Subject to  $\Var[R^{\alpha}(n^{1/m}
\Y_n)] = \Theta(n)$, \cite{Ch} yields a rate of convergence in
(\ref{binomNNCLT}) with respect to the Kantorovich--Wasserstein
distance under minimal assumptions on the distribution of the $Y_i$.
\end{remarks*}

\textit{Shannon entropy}.
The Shannon entropy $H_1(\tk)$ defined at (\ref{Shan}) is an
information theoretic measure of how the data $\{Y_i\}_{i=1}^n$ is
``spread out;'' low entropy implies that the data is confined to a
small volume whereas high entropy indicates the data is widely
dispersed. Accurate estimation of differential entropy is widely
used in pattern recognition, source coding, quantization,
parameter estimation, and goodness-of-fit tests; cf. the survey
\cite{BDGM}.

Shannon differential entropy is commonly estimated by first
estimating the density $\tk$ and then evaluating $H_1(\tk_0)$
where $\tk_0$ is the estimated density;
such methods involve
technical complications involving bin-width selection for histogram
methods and window width for kernel methods and are usually
restricted to $\M= \R^d$. We bypass these technical issues and use
only inter-point data distances to estimate entropy; the methods are
thus applicable to general nonlinear manifolds. This extends
\cite{KL,LPS}, which restricts to $\M= \R^d$, and it lends rigor
to the arguments in~\cite{NK1,NK2}.

As in \cite{KL,LPS}, we shall consider estimators of $H_1(\tk)$ in
terms of nearest neighbor distances. Put $\psi(y, \Y) := \log
(e^{\gamma} \omega_m \NND_1^m(y,\Y))$
where $\gamma:= 0.57721\ldots$ is Euler's constant; see, for example,
\cite{Havil}.

For finite $\Y$ put $S(\Y):= \sum_{y \in\Y} \psi(y,
\Y)$. Define $V^\psi$ by
taking $\xi\equiv\psi$ in
(\ref{Videnti}).
The next result provides the limit theory for the
\emph{Shannon entropy estimators} $S(n^{1/m} \Y_n)$.

%
%th2.4 #&#
\begin{theo} \label{ShCLT}
Suppose $\tk\in\mathbb{P}(\M)$ and
that either \textup{(i)} $\tk\in\mathbb{P}_b(\M)$
or \textup{(ii)}~$\M= \R^m$ and $r_c(\tk) > 0$. Then
as $n \to\infty$,
%
%e2.22 #&#
\begin{equation}
\label{ShLLN} n^{-1} S\bigl(n^{1/m} \Y_n\bigr)
\to H_1(\tk) \qquad \mbox{in }   L^2.   \end{equation}
If $\tk\in\mathbb{P}_c(\M)$
is a.e. continuous, then
%
%e2.23 #&#
\begin{eqnarray}
\label{ShVar}\quad \lim_{n \to\infty} n^{-1} \Var\bigl[S
\bigl(n^{1/m} \Y_n\bigr) \bigr] &=& {\sigma}^2(
\psi, \kappa)
\nonumber
\\[-8pt]
\\[-8pt]
\nonumber
&:=& V^\psi- m^{-2} + \Var\bigl[\log
\tka(Y_1)\bigr] > 0
\end{eqnarray}
and
%
%e2.24 #&#
\begin{equation}
\label{binomShCLT} n^{-1/2} \bigl(S\bigl(n^{1/m}
\Y_n\bigr) - \E S\bigl(n^{1/m} \Y_n\bigr)\bigr)
\stackrel{{\mathcal D}} {\longrightarrow}\NN\bigl(0, {\sigma}^2(\psi,
\kappa)\bigr).
\end{equation}
\end{theo}

\begin{remarks*}
The papers \cite{KL,LPS} show that $n^{-1} S(n^{1/m}
\Y_n)$ consistently estimates Shannon entropy when $\M= \R^m$, but
they do not treat variance asymptotics, distributional results or
general manifolds. Theorem \ref{ShCLT} redresses this. Subject to
$\Var[S(n^{1/m}\Y_n)] = \Theta(n)$, \cite{Ch} yields
(\ref{binomShCLT})
under minimal assumptions on the distribution of the $Y_i$. %, but
\end{remarks*}

%s2.4 #&#
\subsection{Vietoris--Rips clique counts}
\label{subsecVR}

Let $\Y\subset\R^d$ be locally finite, and let $\VRdelta\in(0,
\infty)$ be a scale parameter.
The Vietoris--Rips
complex ${\mathcal R}^{\VRdelta}(\Y)$, also called the Vietoris complex
or Rips complex, is the simplicial complex whose $k$-simplices
correspond to unordered $(k + 1)$ tuples of points of $\Y$ which are
pairwise within Euclidean distance $\VRdelta$ of each other. Thus,
if there is a subset $S$ of $\Y$ of size $k + 1$ with all points of
$S$ distant at most $\VRdelta$ from each other, then $S$ is a
$k$-simplex in the
complex. %For example, if a triple of points are all within
The Vietoris--Rips complex has received attention in connection with
the statistical analysis of high-dimensional data sets \cite{CGOP},
manifold reconstruction~\cite{CO} and gaps in communication
coverage and sensor networks \cite{SG1,SG2}, because of its close
relation to the \v{C}ech complex of a set of balls (which contains a
simplex for every finite subset of balls with nonempty
intersection). It has also received attention amongst topologists
\cite{Ca}, who, given a data cloud $\Y$, allow the scale parameter~$\VRdelta$
to vary to obtain homological signatures of the
Vietoris--Rips complexes, which when taken together, yield clustering
and connectivity information about $\Y$. The limit
theory for the $k$th Betti number of Vietoris--Rips complexes
generated by random points in $\R^d$ is given in \cite{Ka,KM}. The
general methods of this paper may be useful in extending these
results to the setting of manifolds.

Given $\Y, \VRdelta$ and
$k$, let $\Cl_k^{(\VRdelta)}(\Y)$ be the number of $k$-simplices (i.e.,
cliques of order $k + 1$) in ${\mathcal R}^{\VRdelta}(\Y)$. For example,
\[
\Cl_1^{(\VRdelta)}(\Y_n) = \sum
_{1 \leq i < j \leq n} {\1} \bigl\{\|Y_i - Y_j\| <
\VRdelta\bigr\},
\]
which is the empirical version of the so-called correlation integral
\[
\int\int\1 \bigl\{(x,y)\dvtx \|x - y\| < \VRdelta \bigr\} \tka(x) \tka (y) \,dx \,dy.
\]

The
quantity $n^{-1}\Cl_1^{(\VRdelta)}(\Y_n)$, with $\Y_n :=
\{Y_i\}_{i=1}^n$ and with $Y_i$ having unknown distribution $\mu$,
is the widely used sample correlation integral of Grassberger and
Procaccia \cite{GP}. Grassberger and Procaccia use least
squares linear regression of $\log n^{-1}\Cl_1^{(\VRdelta)}(\Y_n)$
versus $\log\VRdelta$ to estimate the ``correlation dimension'' of
$\mu$, that is, the exponent when the correlation integral is assumed to
follow a power law as $\VRdelta\downarrow0$.
The
quantity $n^{-1}\Cl_1^{(\VRdelta)}(\Y_n)$ also features in
estimators of the $K$ function, as discussed in, for example, Chapter 8.2.6
of Cressie \cite{Cr}.

The next result provides a law of large numbers and central limit
theorem for the clique count $n^{-1}\Cl_k^{(\beta)}(\Y_n)$ for any
$k \in\N$.
Define $h_k\dvtx (\R^m)^{k+1} \to\R$ by
$h_k(x_1,\ldots,x_{k+1}): = \prod_{1 \leq i < j \leq k+1}
\mathbf{1} \{ \|x_i-x_j\| \leq1\}$, that is,
% on $(\R^m)^{k+1}$ to be
the indicator
of the event that $x_1,\ldots,x_{k+1}$ are
all within unit distance of each other.
Given $\VRdelta\in(0, \infty)$, put
%
%e2.25 #&#
\begin{eqnarray}\label{Jdef}
 &&J_{k,j} : = \int_{\R^m} \cdots\int
_{\R^m} \biggl( \frac{h(\0,x_1,\ldots,x_{k}) h(\0,x_1, \ldots,x_{j-1},
x_{k+1},\ldots,x_{2k+1 -j}) }{ j! (k+1-j)!^2} \biggr)\hspace*{-25pt}
\nonumber
\\[-6pt]
\\[-10pt]
\nonumber
&&\hspace*{86pt}dx_1 \cdots \,dx_{2k+1-j} ,
\end{eqnarray}
so in particular $J_{k,k+1}:= \int\cdots\int h(\0,x_1,\ldots,x_k)
\,dx_1 \cdots \,dx_k
/(k+1)!$.
Set
%
%e2.26 #&#
\begin{eqnarray}\label{1230a}
\qquad {\sigma}^2_k &:= &{\sigma}^2_k(
\beta,\tka)
\nonumber
\\[-8pt]
\\[-8pt]
\nonumber
&:=& \Biggl( \sum_{j=1}^{k+1}
J_{k,j} \beta^{m(2k+1 - j)} I_{2k+2-j}(\tk) \Biggr)
-
\bigl((k+1)\beta^{mk} J_{k,k+1}I_k(\tka)
\bigr)^2.
\end{eqnarray}

%
%th2.5 #&#
\begin{theo} \label{VRipsthm} Let $\tk$ be bounded on $\M\in\mathbb{M}$.
For all $k = 1, 2,\ldots$ and all $\VRdelta\in(0, \infty)$ we have
%
%e2.27 #&#
\begin{equation}
\label{CliqLLN} \lim_{n \to\infty} n^{-1} \Cl_k^{(\VRdelta)}
\bigl(n^{1/m} \Y_n\bigr) = \VRdelta^{mk}
J_{k,k+1} I_{k+1}(\tk) \qquad  \mbox{in }  L^2  \mbox{ and a.s.}
\end{equation}
If $\tk$ is a.e. continuous and if $\tk\in\mathbb{P}_b(\M)$,
then
%
%e2.28 #&#
\begin{equation}
\label{binomcliqvar} \lim_{n \to\infty} n^{-1} \Var\bigl[
\Cl_k^{(\VRdelta)}\bigl(n^{1/m} \Y_n\bigr)
\bigr] = {\sigma}^2_k(\beta,\tka) > 0,
\end{equation}
and as $n \to\infty$,
%
%e2.29 #&#
\begin{equation}
\label{binomcliqCLT} n^{-1/2} \bigl(\Cl_k^{(\VRdelta)}
\bigl(n^{1/m} \Y_n\bigr) - \E\Cl_k^{(\VRdelta)}
\bigl(n^{1/m} \Y_n\bigr)\bigr) \stackrel{{\mathcal D}} {
\longrightarrow}\NN\bigl(0, {\sigma }^2_k(\beta,\tka)
\bigr).
\end{equation}
\end{theo}
\begin{rema*}[(Related work)]
Bhattacharya
and Ghosh \cite{BhGh} used the limit theory of $U$-statistics
to obtain a central limit theorem similar
to (\ref{binomcliqCLT}) for Poisson input,
in the case where
$\M= \R^m$ and $\kappa$ is uniform on
the unit cube.
The limit (\ref{CliqLLN}) extends the results in Penrose~\cite{Pe} (Proposition 3.1, Theorem 3.17) to nonlinear
manifolds whereas (\ref{binomcliqvar}) and (\ref{binomcliqCLT})
extend Theorem 3.13 of \cite{Pe} to nonlinear manifolds.
\end{rema*}
%s3 #&#
\section{General limit theorems}
\label{gensec} %\allco

As in Section \ref{sectermin},
$\xi(y, \Y)$ is a real-valued
functional defined on locally finite $\Y\subset\R^d$
and $y \in\Y$,
and $Y_i, i \geq1$ are i.i.d.
with density $\tk$. % defined on pairs $(y, \Y)$.
%with $\xi$ assumed translation-invariant
%being determined by inter-point distances in $\Y$.
We shall say that $\xi$ is \emph{translation invariant} if
$\xi(y, \Y) = \xi(z + y , z + \Y)$ for all $z
\in\R^d$ and all $(y,\Y)$.
We say $\xi$ is \emph{rotation invariant} if $\xi(\0,\Y) $
is invariant under rotations of $\Y$, for all $\Y$.
That is, we say $\xi$ is rotation invariant if
$\xi(\0,\Y)= \xi(\0,A\Y)$ for all orthogonal $m \times m$
matrices $A$, where $A\Y:= \{A z\dvtx \in\Y\}$.
In this section
we provide a general limit theory for the sums $\sum_{i=1}^n
\xi(n^{1/m} Y_i, \{n^{1/m}
Y_i\}_{i=1}^n)$, namely
Theorems \ref{smainLLN}, \ref{smainCLT} and \ref{PoCLT}.
We shall
use these general results
to prove the results in Sections \ref{subseclearn}--\ref{subsecVR}.

We introduce a scaled version of $\xi$, dilating the pair $(y, \Y)$
by the
factor $\la^{1/m}$; this scaling is natural when $\Y\subset\M$ has
cardinality approximately $\la$, and $\M$ is an $m$-dimensional
manifold in $\R^d$. Thus, given $k \in\N$, for $\la, \rho\in(
0,\infty)$ set
%
%e3.1 #&#
%e3.2 #&#
\begin{eqnarray}\label{scalexi}
\label{scalexirho} \xi_{\la,k,\rho}(y, \Y)&:=& \xi\bigl(\la^{1/m}y,
\la^{1/m} \Y\bigr) \mathbf{1}\bigl\{ \NND_k(y,\Y) < \rho\bigr\}
;
\\
\xi_\la(y,\Y) &:= &\xi_{\lambda,k,\infty}(y,\Y) : = \xi\bigl(
\la^{1/m}y, \la^{1/m} \Y\bigr).
\end{eqnarray}
Here we are allowing for a finite macroscopic cutoff parameter
$\rho$ because for some manifolds $\xi_n(y,\Y_n)$ may suffer from
nonlocal effects even when $n$ becomes large; see the
counterexamples in remark (i) of Section \ref{subseclearn}.

Given $k \in\Z^+$ and $r>0$, let $ \Xi(k,r)$ be the class of
translation and rotation invariant functionals $\xi$ such that (i)
for all $y,\Y$ with $\card( \Y\setminus\{y\} ) \geq k$ we have
\[
\xi(y,\Y) = \xi \bigl(y, \Y\cap B_{\max(r,N_k(y,\Y))}(y) \bigr)
\]
and (ii) for all $n$, Lebesgue-almost every $(y_1,\ldots,y_n) \in
(\R^m)^{n}$ (with $\R^m$ embedded in $\R^d$)
is at a continuity point of the mapping from $(\R^d)^{n}
\to\R$ given by
\[
(y_1,\ldots,y_n) \mapsto\xi \bigl(\0, \{y_1,
\ldots,y_n \} \bigr).
\]

For finite $\Y\subset\M$, and for $k \in\Z^+$, $\lambda\in
(0,\infty)$ and $\rho\in(0,\infty]$
define %\Comment{$H_\lambda$. MP}
%
%e3.3 #&#
\begin{equation}
\label{rf} H^{\xi}_{\la,k,\rho}(\Y) := \sum
_{y \in\Y} \xi_{\la,k,\rho}(y, \Y) ;   \qquad    H^\xi_{\la}
(\Y) := H_{\la,k,\infty}(\Y). %\ \ \ \mbox{and} \ \ \ \ H^{\xi}_n := \sum_{i=1}^n \xi_n(Y_i, \Y_n).
\end{equation}

Recalling that $\Y_n := \{Y_i\}_{i=1}^n$ and that $\P_\la$ is a
Poisson point process on $\M$ of intensity $\la\tk(y)\,dy$ defined at
(\ref{PPP}), we now give a general law of large numbers for scaled
versions of the linear statistics (\ref{lstat}) when $\xi\in
\Xi(k,r)$. For $i \in\Z^+$, let $\S_i$ be the collection of all
subsets of $\K(\tka)$ of cardinality at most $i$ (including the
empty set). Consider the following moment conditions on $\xi$:
%
%
%e3.4 #&#
%e3.5 #&#
\begin{eqnarray}
\label{smoment0} \sup_n \E\bigl|\xi_{n,k,\rho}(Y_1,
\Y_n)\bigr|^p &<& \infty,
\\
\label{mom5}
\sup_{n \ge1, y \in\K(\tka),
{\mathcal A} \in\S_3} \sup_{(n/2) \le\ell\le(3 n/2)} \E \bigl|\xi_{n,k,\rho} (y,
\Y_\ell\cup{\mathcal A})\bigr|^{p} &<& \infty
\end{eqnarray}
[noting that (\ref{mom5}) implies (\ref{smoment0})], and
%
%
%e3.6 #&#
\begin{equation}
\label{moment1} % \sup_{\lambda\geq1, \ y, y' \in\K}
\sup_{\lambda\geq1, y \in\K(\tka), {\mathcal A} \in\S_1} \E\bigl| \xi_{\la,k,\rho}(y,
\P_\la\cup{\mathcal A} ) \bigr|^p < \infty.
\end{equation}

%
%th3.1 #&#
\begin{theo}[{[Laws of large numbers for $\xi\in\Xi(k,r)$]}]\label{smainLLN}
Let $\M\in\mathbb{M}$, $\tka\in\mathbb{P}(\M)$, $k \in\Z^+$
and $\rho\in
(0,\infty]$,
and put $q = 1$ or $q = 2$. Let $\xi\in\Xi(k,r)$, and suppose there exists
$p
> q$ such that
(\ref{smoment0}) holds.
Then as $n \to\infty$ we have $L^q$ convergence
%
%e3.7 #&#
\begin{equation}
\label{sWLLN} n^{-1} H_{n,k,\rho}^\xi(
\Y_n) \to \int_{\M} \E\bigl[ \xi(\0,
\H_{\kappa(y)}) \bigr] \kappa(y) \,dy.
\end{equation}
If also
(\ref{mom5}) holds for some $p>5$,
and also $\tka\in\mathbb{P}_b(\M)$
and moreover either
$k =0$ or
$\tka\in\mathbb{P}_c(\M)$,
then (\ref{sWLLN}) holds a.s. %with almost sure
\end{theo}

Given $a >0$, let $\H_a$ denote a homogeneous Poisson
process of intensity $a$ in $\R^m$ (embedded in
$\R^d$). Extending the earlier definitions (\ref{Videnti}) and
(\ref{Didenti}),
we set
%
%e3.8 #&#
\begin{equation}\label{Vdef}
\qquad V^{\xi}(a):= \E\xi(\0, \H_{a})^2 +  a \int
_{\R^m} \bigl\{ \E\xi\bigl(\0, \H_{a}^{u}
\bigr) \xi\bigl(u, \H_{a}^{\0}\bigr) - \bigl(\E\xi(\0,
\H_{a})\bigr)^2 \bigr\} \,du
\end{equation}
and
%
%e3.9 #&#
\begin{equation}
\label{ddef} \delta^{\xi}(a) := \E\xi(\0, \H_{a } ) + a \int
_{\R^m} \E\bigl[\xi\bigl( \0, \H_a^{u}
\bigr) - \xi(\0, \H_a )\bigr] \,du,
\end{equation}
so in particular
$
V^{\xi}(1)
=
V^\xi
$
and $\de^\xi(1) =\de^{\xi}$.
For all $\tka\in\mathbb{P}(\M)$ and
$\xi$,
we define
%
%e3.10 #&#
\begin{equation}
\label{hatsig} {\sigma}^2(\xi, \tka) := \int_{\M}
V^{\xi}\bigl(\tka(y) \bigr) \tka(y) \,dy - \biggl( \int
_{\M} \de^{\xi} \bigl(\tka(y)\bigr) \tka(y) \,dy
\biggr)^2,
\end{equation}
provided that both integrals in
(\ref{hatsig}) exist and are finite.

%
%th3.2 #&#
\begin{theo}[{[Variance asymptotics and CLT for $\xi\in\Xi(k,r)$]}] \label{smainCLT}
Let $\M\in\mathbb{M}$ and let $\tk\in\mathbb{P}_b(\M)$ be a.e.
continuous.
Let $k \in\Z^+$, $r \geq0$, and $\rho\in(0,\infty]$.
Assume $k=0$ or $\kappa\in\mathbb{P}_c(\M)$.
Let $\xi\in\Xi(k,r)$ and suppose that $\xi$ satisfies
(\ref{mom5}) and (\ref{moment1})
for some $p > 2$.
Then $\sigma^2(\xi,\kappa)< \infty$
and
%
%e3.11 #&#
\begin{equation}
\label{sbinomialVAR} \lim_{n \to\infty} n^{-1} \Var
\bigl[H_{n,k,\rho}^\xi(\Y_n)\bigr] = {\sigma
}^2(\xi, \kappa)
\end{equation}
and as $n \to\infty$,
%
%e3.12 #&#
\begin{equation}
\label{sbinomCLTlim} n^{-1/2} \bigl(H_{n,k,\rho}^{\xi}(
\Y_n) - \E H_{n,k,\rho}^{\xi}(\Y_n)\bigr)
\stackrel{{\mathcal D}} {\longrightarrow }\NN\bigl(0, {\sigma}^2(\xi,
\kappa)\bigr).
\end{equation}
\end{theo}

\begin{remarks*} (i) (\textit{Related work}.) Under a slightly different
set of assumptions, Chatterjee \cite{Ch} provides estimates for the
Kantorovich--Wasserstein distance between the distribution of
$H^\xi_{n} (\Y_n)$ and the normal, which imply a central limit
theorem subject to the validity of $\Var[H_{n}^\xi(\Y_n)] =
\Theta(n)$. Indeed, if in Theorem 3.2 we take $\rho=\infty$, $r=0$
and $p >8$, and if also $\sigma^2(\xi,\kappa) >0$, then combining
(\ref{sbinomialVAR})
with Theorem 3.4 of \cite{Ch}
shows that the convergence (\ref{sbinomCLTlim})
is at rate $O(n^{4/p - 1/2}) $ with respect to
the Kantorovich--Wasserstein distance.

(ii) (\emph{Simplification of mean and variance asymptotics for
homogeneous $\xi$}.) The limits (\ref{sWLLN}) and
(\ref{sbinomialVAR}) take a simpler form when
there is some $\beta\in\R$ such that
$\xi$ is \emph{homogeneous of order $\beta$}, meaning
that for all $a \in(0, \infty)$ and all $y \in\Y\subset\R^m$,
we have $\xi(ay, a\Y) = a^\beta\xi(y, \Y)$. In this case
$\xi(\0,\H_a) \stackrel{{\mathcal D}}{=}\xi(\0,a^{-1/m}\H) =
a^{-\beta/m} \xi(\0,\H)$,
so
(\ref{sWLLN}) becomes
%
%e3.13 #&#
\begin{equation}
\label{newedgeWLLN} n^{-1} H_n^\xi(
\Y_n) \to I_{1-\beta/m}(\tka) \E \xi(\0, \H)
\end{equation}
and similarly,
by definitions (\ref{Videnti}), (\ref{Didenti}), we can show that
(\ref{Vdef}) and (\ref{ddef}) simplify to
$V^\xi(a) = a^{-2 \beta/m} V^\xi$ and $\delta^\xi(a) =
a^{-\beta/m} \delta^\xi$. Hence, in this case
(\ref{hatsig}) becomes
%
%e3.14 #&#
\begin{equation}
\label{binomialedgevar} \sigma^2(\xi,\tka) = V^{\xi}
I_{1-2\beta/m} (\tka) - \bigl( \delta^\xi I_{1 - \beta/m}(\tk)
\bigr)^2.  \end{equation}
If $\xi$ is homogeneous of order 0, we say it is
\emph{scale invariant.}

(iii) (\emph{Assumptions on $\kappa$}.) While some of the general
results allow for $\tk\in\mathbb{P}_b(\M)$, it is often necessary to
require $\tk\in\mathbb{P}_c(\M)$ in order to guarantee that functionals
$\xi\in\Xi(k,r)$ satisfy the spatial localization property termed
exponential stabilization, described in Section \ref{sec6}.
Additionally, in Section \ref{secresults}, the requirement $\tk\in
\mathbb{P}_c(\M)$
is needed when verifying that the estimators of that section
satisfy the general results given in Theorems~\ref{smainLLN} and~\ref{smainCLT}.

(iv) (\emph{Limit theory for random measures on
manifolds}.) Consider the point measures
%
%e3.15 #&#
\begin{equation}
\label{rm} \mu^\xi_{\la,k,\rho}:= \sum
_{i=1}^n \xi_{n,k,\rho}(Y_i,
\Y_n) \delta_{Y_i},
\end{equation}
where $\delta_y$ denotes the unit point mass at $y$. As in
\cite{BY2,PeEJP,PeBer}, Theorem \ref{smainLLN} admits an extension
to the random measures (\ref{rm}) as follows. Let $\B(\M)$ be the
space of bounded, measurable, real-valued functions on $\M$ and for $f
\in\B(\M)$, and $\mu$ a measure on $\M$, let $\langle f, \mu
\rangle$ denote the integral of $f$ with respect to $\mu$.
Put $q = 1$ or
$q = 2$. Let $\xi\in\Xi(k,r)$ and suppose that there is a $p >
q$ such that (\ref{smoment0}) and
(\ref{moment1})
are satisfied. It
can be shown that for all $f \in\B(\M)$,
%
%e3.16 #&#
\begin{equation}
\label{2WLLN} \lim_{n \to\infty} n^{-1} \bigl\langle f,
\mu_{n,k,\rho}^\xi\bigr\rangle = \int_{\M}
f(y) \E\bigl[ \xi(\0, \H_{\kappa(y)}) \bigr] \kappa(y) \,dy   \qquad\mbox{in }   L^q.  \end{equation}

Similarly, if $\xi\in\Xi(k,r)$ satisfies the moment assumptions
of Theorem \ref{smainCLT},
it can be shown that
$n^{-1/2} (
\langle f, \mu_{n,k,\rho}^\xi\rangle- \E\langle f,
\mu_{n,k,\rho}^\xi\rangle)$ converges in distribution to a
mean zero normal random variable with variance
\[
\int_{\M} f(y)^2 V^\xi \bigl({
\tka}(y) \bigr) {\tka}(y) \,dy - \biggl( \int_\M
\delta^\xi \bigl(\tka(y) \bigr) f(y) \tka(y)  \,dy
\biggr)^2 .
\]
We refer to
\cite{BY2,PeEJP} for details.

(v) (\emph{Other functionals}.) The approach also works for more
general functionals than $\xi\in\Xi(k,r)$. Along with the moment
conditions already discussed, the key properties $\xi$ needs to
satisfy are
\emph{exponential stabilization}
and \emph{continuity}, which we discuss in the proof and which
represent $\xi$ being locally determined in some sense; cf. the remark
at the end of Section \ref{sec6}. Also, the approach
works for marked point processes,
where the points carry independent identically distributed
marks. We would expect other functionals to satisfy this, as has
been considered for a variety of functionals in the case $\M= \R^m$
\cite{BY2,PeEJP}.
Finally,
%we believe
%that
the assumption that $\xi$ is rotation invariant
could be relaxed; one would need to
change definition (\ref{scalexi}) to
\[
\xi_\la(y,\Y) = \xi \bigl(y,y+\la^{1/m}(-y+ \Y) \bigr),
\]
modify (\ref{scalexirho}) similarly, and
in (\ref{sWLLN})--(\ref{hatsig})
take $\H_{\kappa(y)}$ to
be a homogeneous Poisson process in the
hyperplane tangent to $-y +\M$ at $\0$,
and
%in the $V^\xi(\kappa(y))$ appearing in
take the integrals in
(\ref{Vdef}) and (\ref{ddef}) to be over this
tangent hyperplane rather than over $\R^m$.

% This comment seems to make the paper shorter!

(vi) (\emph{Noisy input}.) It is arguably more realistic to consider
input having a $d$-dimensional noise component. Consider the
situation where data $Y_i, i \geq1$, is corrupted with a noise
component $n^{-1/m}Z_i$, with $Z_i, i \geq1$, being i.i.d.
$\R^d$-valued random variables which are independent of $Y_i, i \geq
1$, and which are assumed rotation invariant, for example, with
i.i.d. mean zero
normal components.

For all $x \in\R^d$, let $Z_x$ denote a copy of $Z_1$.
%Put ${\mathcal Z}_n := \{X_i \}_{i=1}^n.$
For all $u \in\R^m, a \in(0, \infty)$,
put $\H_a^{u, {\mathcal Z}}:= \{ x + Z_x\dvtx x \in\H_a \cup\{u \} \},$
and put
$\H_a^{ {\mathcal Z}}:= \{ x + Z_x\dvtx x \in\H_a \}.$

Then it can be shown that the law of large numbers
(Theorem \ref{smainLLN}) takes the form
\begin{eqnarray*}
&&\lim_{n \to\infty} \sum_{i=1}^n
\xi_n \bigl(Y_i + n^{-1/m}Z_i,
\bigl\{Y_j + n^{-1/m} Z_j \bigr
\}_{j=1}^n % \Y_n + n^{-1/m}{\mathcal Z}_n
\bigr) \\
&&\qquad= \int_{\M}
\E \bigl[ \xi \bigl(\0 + Z_\0, \H_{\kappa(y)}^{ {\mathcal Z}}
\bigr) \bigr] \kappa(y) \,dy.
\end{eqnarray*}
Moreover, it can be shown that
Theorem \ref{smainCLT} still holds
if in the statement of that result we replace $\Y_n$ by
$\{Y_j + n^{-1/m} Z_j\}_{j=1}^n$ and replace
definitions
(\ref{Vdef})
and (\ref{ddef}), respectively, by
\begin{eqnarray*}
V^{\xi}(a)&:=& \E\xi\bigl(\0, \H^{{\mathcal Z}}_{a}
\bigr)^2 \\
&&{}+  a \int_{\R^m} \bigl\{ \E\xi\bigl(\0
+Z_\0, \H_{a}^{u,{\mathcal Z}}\bigr) \xi\bigl(u
+Z_u, \H_{a}^{\0,{\mathcal Z}}\bigr)
- \bigl(\E\xi\bigl(\0+Z_\0, \H_{a}^{{\mathcal Z}}\bigr)
\bigr)^2 \bigr\} \,du
\end{eqnarray*}
and
\begin{eqnarray*}
\delta^{\xi}(a) := \E\xi\bigl(\0+Z_\0,
\H_{a }^{{\mathcal Z}} \bigr) + a \int_{\R^m} \E
\bigl[\xi\bigl( \0, \H_a^{u,{\mathcal
Z}} \bigr) - \xi\bigl(\0,
\H_a^{{\mathcal Z}} \bigr)\bigr] \,du.
\end{eqnarray*}

(vii) (\emph{Poisson input}.)
In (\ref{sWLLN})
we have presented the law of large numbers
for binomial samples, but the same
$L^q$ limit holds for functionals of the form $\lambda^{-1}
H_{\la,k,\rho}^\xi(\Po_\lambda)$, with $\Po_\lambda$ as in
(\ref{PPP}). Likewise, there is a Poisson analog to~(\ref{2WLLN}).

We also have variance asymptotics and a central limit theorem
for\break $H_{\la,k,\rho}^\xi(\Po_\lambda)$, similar
to Theorem \ref{smainCLT} but with a different limiting variance.
Moreover, in the Poisson setting, we have a bound
on the rate of convergence to the normal, using
the Kolmogorov distance.
The result goes as follows.
\end{remarks*}

%
%th3.3 #&#
\begin{theo}
\label{PoCLT}
Let $\M\in\mathbb{M}$ and let $\tk\in\mathbb{P}_b(\M)$ be a.e.
continuous.
Let $k \in\Z^+$, $r >0$ and $\rho\in(0,\infty]$, and
suppose $k=0$ or $\tk\in\mathbb{P}_c(\M)$.
Let $\xi\in\Xi(k,r)$ and suppose that $\xi$ satisfies
(\ref{moment1}) for some $p > 2$. Then
%
%e3.17 #&#
\begin{equation}
\label{sVAR} \qquad \lim_{\la\to\infty} \la^{-1} \Var\bigl[
H^\xi_{\la,k,\rho}(\P_\la)\bigr] = \tau^2(
\xi, \kappa) := \int_{\M} V^{\xi}\bigl(y,\kappa(y)
\bigr) \kappa(y) \,dy < \infty,
\end{equation}
and as $\la\to\infty$,
%
%e3.18 #&#
\begin{equation}
\la^{-1/2} \bigl( H^\xi_{\la,k,\rho}(\P_\la) -
\E H^\xi_{\la,k,\rho}(\P_\la) \bigr) \stackrel{{
\mathcal D}} {\longrightarrow}\NN\bigl(0, \tau^2(\xi, \kappa ) \bigr).
\label{0324a}
\end{equation}
Additionally, if $\xi$ satisfies
(\ref{moment1})
for some $p > 3$ %\emph{[can we change 3 to 2? NO, JY]} \Comment{MP}
and if $\sigma^2(\xi, \kappa) > 0$, then there
exists a finite constant $C$ depending on $d, k, \xi$, $\tk$, $\rho$
and $p$
such that for all $\la\geq2$,
%
%e3.19 #&#
\begin{equation}
\label{sSteinbd}\qquad \sup_{t \in\R} \biggl\llvert P \biggl[ \frac{ H_{\la,k,\rho}^{\xi}(\P_\la
) -
\E
H_{\la,k,\rho}^{\xi}(\P_\la) }{ \sqrt{ \Var[H_{\la,k,\rho
}^{\xi}(\P_\la) ] } }
\leq t \biggr] - \Phi(t) \biggr\rrvert \leq C ( \log \la)^{3m}
\la^{-1/2}.
\end{equation}
\end{theo}
As well as being of independent interest,
Theorem \ref{PoCLT} is used in our proof of
Theorem \ref{smainCLT}.
Equation (\ref{sSteinbd}) is the counterpoint for manifolds to the rate
of normal convergence result in \cite{PY5}. Theorem \ref{PoCLT}
could be used to provide Poisson analogs to the results presented in
Sections \ref{subseclearn}--\ref{subsecVR}.

%s4 #&#
\section{Geometrical preliminaries}
\label{balls}
The lemmas in this section, concerned with
properties of manifolds, have no probabilistic content.
The first of these relates
distances in the manifold to distances in a chart.

%
%le4.1 #&#
\begin{lemm}
\label{DGlem}
Suppose $(U,g)$ is a chart for $\M\in\mathbb{M}(m,d)$.
Suppose $F \subset g(U)$ is a compact subset of $\M$.
Then
%
%e4.1 #&#
\begin{equation} \label{101221a}
 0< \inf_{y,z \in F\dvtx y \neq z} \frac{\|g^{-1}(z) - g^{-1}(y)\|}{\|z-y\|} \leq \sup_{y,z \in F\dvtx y \neq z}
\frac{\|g^{-1}(z) - g^{-1}(y)\|}{\|z-y\|} < \infty.\hspace*{-25pt}
\end{equation}
\end{lemm}
\begin{pf}
Suppose (\ref{101221a}) fails. Then by compactness,
we can find a sequence
$(y_n,z_n), n \in\N$ with $y_n \in F $,
$z_n \in F \setminus\{y_n\}$,
and $y_n \to y$ for some $y \in F$, such that
setting $u_n := g^{-1}(y_n)$ and $v_n:=g^{-1}(z_n)$,
we have either
%
%e4.2 #&#
\begin{equation}\label{101221b}
\|v_n-u_n\| / \|z_n -y_n\| \to
\infty     \qquad  \mbox{as }    n \to\infty
\end{equation}
or
%
%e4.3 #&#
\begin{equation}\label{1101a}
\|v_n-u_n\| / \|z_n -y_n\| \to0
      \qquad\mbox{as }    n \to\infty.
\end{equation}
Since $g$ is an open map, the set $g^{-1}(F)$ is
compact. Hence $\|v_n-u_n\|$ remains
bounded.

Suppose (\ref{101221b}) holds. Then
$\|z_n - y_n\| \to0$,
and hence $z_n \to y$ as $n \to\infty$.
Setting $u:= g^{-1}(y)$, by continuity of $g^{-1}$
we have $u_n \to u$ and $v_n \to u$. Hence,
arguing componentwise using the mean value theorem and the
continuity of $g'$, we have that
%
%e4.4 #&#
\begin{equation}\label{0105a}
\bigl\| z_n - y_n - g'(u) (v_n
-u_n) \bigr\| = o\bigl(\|v_n -u_n \| \bigr)
\end{equation}
and therefore since $g'(u)$ has full rank,
\[
\liminf_{n \to\infty} \frac{\| z_n - y_n\|}{\| v_n -u_n\|} = \liminf_{n \to\infty}
\frac{\| g'(u) ( v_n - u_n) \|}{\| v_n -u_n\|} >0,
\]
which contradicts (\ref{101221b}). On the other hand,
if (\ref{1101a}) holds we can show by a similar argument
that $
\limsup_{n \to\infty}
( \| z_n - y_n\|/ \| v_n -u_n\|) < \infty$,
again giving a contradiction.
\end{pf}

% Recall that $\K(\tka)$ denotes the
%support of $\tk$, and if $\tk\in\p_c(\M)$, then \Comment{ }
%$\K(\tka)$ is
%locally conic and satisfies \eqref{lconic}.
%a compact $C^1$
%submanifold-with-boundary of $\M$.
For $w \in E \subseteq\R^d$ and $r
> 0$, let $ \label{balldef} B_r^E(w):= B_r(w) \cap E $. Recall that
$\omega_m:= \pi^{m/2}[\Gamma(1 + m/2)]^{-1}$ is the volume of the
ball $B_1(\0)$ in $\R^m$.

%
%le4.2 #&#
\begin{lemm}
\label{Lem0729}
Let $\M\in\mathbb{M}$.
%Suppose $\tk\in\p_c(\M)$.
Suppose
$y_\infty\in\M$,
and suppose for $n \in\N$ we are given $y_n \in
\M,
r_n >0, a_n >0$
with $y_n \to y_\infty$, $r_n \to0$ and $a_n \to0$
as $n \to\infty$. Then %\Comment{removed all mention of
%$\K$ round here}
%
%e4.5 #&#
\begin{equation}
\limsup_{n \to\infty} \biggl( r_n^{-m} \int
_{B^\M_{r_n}(y_n) }\,dy \biggr) \leq\omega_m \label{0728a},
\end{equation}
%
%and \bea\liminf_{n \to\infty} \left( r_n^{-m}
and putting $s_n = r_n(1-a_n)$, we have
%
%e4.6 #&#
\begin{equation}\label{0729a}
\limsup_{n \to\infty} \biggl( \bigl(r_n^m
-s_n^m\bigr)^{-1} \int_{B^\M_{r_n}(y_n) \setminus B^\M_{s_n}(y_n) }\,dy
\biggr) < \infty.
\end{equation}
\end{lemm}

\begin{pf} Let $(U,g)$ be a chart such that $\0 \in U$ and $g(\0)
= y_\infty$.
%Assume also the chart is chosen so that $U \cap
%g^{-1}(\K) \supseteq U \cap( \R^+ \times\R^{m-1} ) $ (with
%equality if $y_\infty$ is on the boundary of $\K(\tka)$).
Set $U_n = g^{-1}(B^\M_{r_n}(y_n)) \subseteq U$, and set
$V_n = g^{-1}(B^\M_{s_n}(y_n)) $.
Let ${\mathcal L}$ denote Lebesgue measure, and note
that by continuity $\sup_{u \in U_n} |(D_g(u)/D_g(\0))-1| $ vanishes
as $n \to\infty$.
Thus there exists $n_0$
such that for $n \geq n_0$, we have
$ B_{r_n}^\M(y_n) \subset g(U)$, and
by (\ref{int-form}),
%
%e4.7 #&#
\begin{equation}\label{0103a}
\int_{B^\M_{r_n}(y_n) } \,dy = D_g(\0) \int
_{U_n} \bigl(D_g(u)/D_g(\0)\bigr) \,du
\sim D_g(\0) {\mathcal L}(U_n) %    \mbox{as}   n \to\infty,
\end{equation}
and
%
%e4.8 #&#
\begin{eqnarray}\label{0729c}
&&\int_{B^\M_{r_n}(y_n) \setminus B^\M_{s_n}(y_n) } \,dy
\nonumber
\\[-8pt]
\\[-8pt]
\nonumber
&&\qquad = D_g(\0) \int
_{U_n \setminus V_n} \bigl(D_g(u)/D_g(\0)\bigr) \,du
\sim D_g(\0) {\mathcal L}(U_n \setminus V_n)
\end{eqnarray}
%
% and
% \bean
% \int_{B_{r_n}^{\K(\tka)}(y_n) }dy
%D_g(\0) {\mathcal L}(U_n \cap(\R^+ \times\R^{m-1})),
%    \mbox{as}   n \to\infty.
where the asymptotics are as $n \to\infty$. Given $n \geq
n_0$, set $u_n := g^{-1}(y_n) \in U$.
We claim that
%
%e4.9 #&#
\begin{equation}\label{0103b}
\limsup_{n \to\infty} \sup_{v \in U_n} r_n^{-1}
\bigl\|g'(\0) (v-u_n)\bigr\| \leq 1.
\end{equation}
%
% \|g'(\0)(v-u_n)\| \geq1.
To see (\ref{0103b}), take
$v_n \in U_n$
for $n \in\N$.
By continuity of $g^{-1}$, we
have $v_n \to\0$ as $n \to\infty$
so by applying the mean value theorem
and using continuity of $g'$,
as in~(\ref{0105a})
we have that
$
\| g( v_n) - y_n - g'(\0)(v_n -u_n) \| = o(\|v_n
-u_n \|)
$
as
%    \mbox{as}
$
n \to\infty.
$
Therefore
since $g'(\0)$ has full rank,
%
%e4.10 #&#
\begin{equation}\label{0103d}
\bigl\| g( v_n) - y_n \bigr\| \sim g'(\0)
(v_n-u_n)     \qquad \mbox{as }   n \to\infty.
\end{equation}
Then (\ref{0103b}) follows because
$\|g(v_n)- y_n \| \leq r_n$ and the choice of $v_n \in U_n$
was arbitrary.

%To see \eq{0103c}, let $v_n \in U \setminus U_n$
%with $v_n \to\0$ as $n \to\infty$. Then
%$\| g(v_n) - y_n \| \geq r_n$ and \eq{0103c} follows since
%the choice of $v_n \in U \setminus U_n$
%was arbitrary.

By (\ref{0103b}), given $\eps>0$ we have for large enough $n$ that
\[
{\mathcal L}(U_n) \leq{\mathcal L} \bigl( \bigl\{ v \in
\R^m\dvtx \bigl\| g'(\0) (v-u_n) \bigr\| \leq
r_n (1+\eps) \bigr\} \bigr) = (1+\eps)^m
\omega_m r_n^m / D_g(\0),
\]
and (\ref{0728a}) follows by (\ref{0103a}).
%By \eq{0103c}, given $\eps>0$
%we have for large $n$ that
% \bean
% {\mathcal L}(U_n \cap(\R_+ \times\R^{m-1}) ) \geq
% {\mathcal L}(\{ v \in\R^m: \| g'(\0) (v-u_n) \| \leq r_n
%(1-\eps)^{-1} \})
%= \frac{1}{2}(1-\eps)^{-m} \omega_m r_n^m / D_g(\0), \eean and

Finally we prove (\ref{0729a}). Since
$(1-a)^m \leq1 - a$
for all $a \in[0,1]$,
for large $n$ we have
%
%e4.11 #&#
\begin{equation}\label{0729b}
r_n^m - s_n^m =
r_n^m\bigl(1-(1-a_n)^m\bigr) \geq
r_n^m a_n.
\end{equation}
We claim that
if
$(v_n)_{n \in\N}$ is a sequence
in $U$ with
$\|g(v_n ) - y_n\|
= r_n$,
then setting
$w_n := u_n + (1-3 a_n)(v_n-u_n) $, we have
$\|g(w_n) - y_n\| < s_n$ for large enough $n$.
Indeed, for any such sequence
by the mean value theorem
we have
$r_n = \|g(v_n) -y_n\| \sim\| g'(\0)(v_n-u_n)\| $
as in (\ref{0103d}), and also
\[
% \|g(v_n + (t_n - C' a_n t_n)e_n) -
\bigl\|g(w_n ) - % g(v_n + t_n e_n)
g(v_n ) \bigr\| \sim3
a_n \bigl\| g'(\0) (v_n -u_n) \bigr\|
\sim3 a_n r_n.
\]
Moreover
$g(w_n ) - y_n $ and
$g(v_n ) - y_n $ are almost in the same
direction,
% of $g'(0) e_n$,
so that for large $n$,
$
\|g(w_n ) - y_n \| \leq r_n( 1- 2 a_n) ,
$
and the claim follows.

By the preceding claim, there is a constant $C$ such that
the thickness of the deformed annulus
$U_n \setminus V_n$ in all
directions is bounded by $C a_n r_n$, so by using polar
coordinates we have ${\mathcal L}(U_n \setminus V_n) = O( r_n^m a_n)$,
and then using (\ref{0729c}) and (\ref{0729b}) we have that
\[
\int_{B^\M_{r_n}(y_n) \setminus B^\M_{s_n}(y_n) } \,dy = O \bigl( r_n^m
a_n \bigr) = O \bigl(r_n^m -
s_n^m \bigr),
\]
demonstrating (\ref{0729a}).
\end{pf}

Recall that $\K(\tka)$ denotes the
support of $\tk$, and if $\tk\in\mathbb{P}_c(\M)$, then %\Comment{
%}
$\K(\tka)$ is
locally conic and satisfies \eqref{lconic}.
Given $\M\in\mathbb{M}$
and $\kappa\in\mathbb{P}_c(\M)$,
set $ \Delta(\tka) := \diam( \K(\tka))$.
% := \sup\{\|x - y\|: \ x,y \in\K(\tka) \}$.
%
%
%le4.3 #&#
\begin{lemm}
\label{claim0} Suppose $\mathbb{M}\in\M$ and
$\tk\in\mathbb{P}_c(\M)$. Then there is a constant
$C_0 \in(0, \infty)$ such that for all $r \in(0, \Delta(\tka)]$ and
$w \in\K(\tka)$, we have
%
%e4.12 #&#
\begin{equation}
\label{bds0} C_0^{-1} r^m \leq \int
_{B_r^{\K(\tka)}(w) } \,dy \leq\int_{B_r^\M(w) } \,dy \leq
C_0 r^m.
\end{equation}
There are also positive finite constants $C_1 $ and $\rho_1 $
such that if $0 < s <r < \rho_1$ and $w \in\K(\tka)$, then
%
%e4.13 #&#
\begin{equation}
\int_{B_r^\M(w) \setminus B_s^\M(w)} \,dy \leq C_1 \bigl(r^m -
s^m\bigr). \label{0729d}
\end{equation}
\end{lemm}

\begin{pf} In the proof, set $\K:= \K(\tka)$ and $\Delta:=
\Delta(\tka)$. The first inequality in~(\ref{bds0})
(for large enough $C_0$) follows
from the
% holds by definition of
assumption that $\K$ is
locally conic~(\ref{lconic}).
% provided $C_0^{-1} \geq L$, with $L$ as in
Suppose the last inequality of
(\ref{bds0}) fails; then there must be a $(\K
\times(0,\Delta])$-valued
sequence $\{(y_n,r_n),n \in\N\}$
such that %either
% \bea
% \lim_{n \to\infty}
%r_n^{-m}
% or
%
%e4.14 #&#
\begin{equation}\label{0728g}
\lim_{n \to\infty} r_n^{-m} \int_{B_{r_n}^\M(y_n) }
\,dz = \infty.
\end{equation}

Since $\K\times[0,\Delta]$ is compact, by taking a subsequence
we may assume without loss of generality that $y_n \to y$ and
$r_n \to r$ for some $y \in\K$ and $r \in[0,\Delta]$.
If $r =0$,
then %\eq{0728f} would contradict \eq{0728b} and
(\ref{0728g}) would contradict (\ref{0728a}).
%If $r >0$
%and \eq{0728f} holds, then since $B_{r/2}^{\K}(y) \subset B_{r_n}^{
%for large enough $n$, we would have
%$\int_{B_{r/2}^{\K}(y) } \,dz =0$, which is impossible.
If $ r >0$ and (\ref{0728g}) holds, then since $B_{r_n}(y_n) \subset
B_{2r}(y)$ for large $n$ we have $\int_{B_{2r}(y)}\,dz =\infty$,
which is impossible: indeed by compactness $B_{2r}(y) \cap\M$
is covered by finitely many of the regions $g_i(U_i)$,
and $\int_{g_i(U_i)} \kappa_i(x)\,dx$ is finite for all $i$ (we may
assume the charts were chosen so all the regions $U_i$ are bounded).
Therefore we have a contradiction so
(\ref{bds0}) must hold.

It remains to prove there exists positive $\rho_1$ such that
(\ref{0729d}) holds for all $w \in\K$ and $0< s < r < \rho_1$.
Suppose this is not the case.
Then there is a sequence $\{(y_n,r_n,a_n), n \in\N\}$
taking values in $\K\times(0,\Delta] \times(0,1)$ such that $r_n
\to0$, and setting $s_n = r_n(1-a_n)$ we have
%
%e4.15 #&#
\begin{equation}\label{0729e}
\lim_{n \to\infty} \bigl(r_n^m - s_n^m
\bigr)^{-1} \int_{B_{r_n}^\M(y_n) \setminus B_{s_n}^\M(y_n)} \,dy = \infty.
\end{equation}
By taking a subsequence, we may
assume that $y_n \to y$ for some $y \in\K$, and either~$a_n$ is
bounded away from zero or $a_n \to0$ as $n \to\infty$.

If $ \inf_{n \in\N}\{a_n\} >0$, then $(r_n^m -s_n^m)^{-1} =
O(r_n^{-m})$ so (\ref{0729e})
would give a contradiction of
(\ref{bds0}).
If $a_n \to0$
and $r_n \to0$, then
(\ref{0729e}) would give a contradiction of~(\ref{0729a}).
\end{pf}

%
%re4.1 #&#
\begin{rema}
\label{bdyrmk}
A sufficient condition for $\K\subset\M$ to be
locally conic is that $\K$ be a compact
$m$-dimensional $C^1$ submanifold-with-boundary of $\M$.
%This can be proved by similar arguments to
%the proof of Lemmas \ref{Lem0729} and \ref{claim0}; for details
%see the proof of these lemmas in the earlier version
%of this paper \cite{Version1} (where
%the definition of $\p_c$ is different from here).
\end{rema}
This can be proved by similar arguments to
the proof of Lemmas \ref{Lem0729} and~\ref{claim0}; for details
see the proof of these lemmas in the earlier version
of this paper \cite{Version1} (where
the definition of $\mathbb{P}_c$ is different from here).

%s5 #&#
\section{Weak convergence lemmas}
\label{prfsec}

For all $d \in\N$, we put a topology $\T:= \T_d$ on locally
finite point sets
in $\R^d$.
As in Aldous and Steele
(\cite{AS}, page 250), we adopt a topology
whereby a sequence of locally
finite point sets
$(\yy_n)_{n \geq1}$ converges to a locally finite $\yy$, if and only if
(i) it is possible to list the elements of $\yy$ as a possibly
terminating sequence $(y_i, i \geq1)$ and the elements of $\yy_n$
as a possibly terminating sequence $(y_{n,i}, i \geq1)$ in such a
way that
%
%e5.1 #&#
\begin{equation}\label{ppconv1}
\lim_{n \to\infty} y_{n,i} = y_i\qquad    \forall i
\end{equation}
and (ii) for any $L$ with no point of
$\yy$ on the boundary of $B_L(\0)$,
we have
%
%e5.2 #&#
\begin{equation}\label{ppconv2}
\lim_{n \to\infty} \card\bigl( \yy_n \cap B_L(\0)
\bigr) = \card\bigl( \yy\cap B_L(\0) \bigr).
\end{equation}

We would like to know that whenever point sets $\yy_n \subset
\R^d$ are close to the point set $\yy$ in the topology $\T_d$,
then $\xi(y, \yy_n)$ is close to $\xi(y, \yy), y \in\R^d.$ This
motivates Definition \ref{hstab} below.
Recall that $\H$ denotes a homogeneous
Poisson point process of unit intensity on $\R^m$.

%
%de5.1 #&#
\begin{defn} \label{hstab}
$\xi$ is \emph{continuous} if
for any linear $F\dvtx \R^m \to\R^d$ of full rank,
for almost all $z \in\R^m$
both $F(\H)$ and
$F(\H^z)$ lie a.s. at continuity points of
$\xi(\0,\cdot)$ with respect to $\T_d$.
\end{defn}

If $\U_n$ and $\U$ are simple point processes (i.e.,
random locally finite point sets in $\R^d$), then, following the
discussion in \cite{AS} page 251, we shall say that $\U_n$
converges in distribution to $\U$ if the law of $\U_n$ converges
weakly to that of $\U$ under the (metrisable) topology $\T_d$,
which is the same as the notion of weak convergence of point
processes discussed in Daley and Vere-Jones \cite{DVJ}, Section
11.1.

Given
the atlas $((U_i,g_i), i \in{\mathcal I})$,
for $i \in\I$ define the function $\ka_i\dvtx U_i \to[0,\infty)$ by
%
%e5.3 #&#
\begin{equation}
\label{kappi}\ka_i(x) = \tk\bigl(g_i(x)\bigr)
D_{g_i}(x),   \qquad   x \in U_i.
\end{equation}
By using (\ref{int-form}) with a partition of unity for which
$\psi_i \equiv1$ on $g_i(U_i)$, we see that
for Borel $B \subseteq U_i$,
%
%e5.4 #&#
\begin{equation}\label{0506a}
\int_{g_i(B)} \tk(y) \,dy = \int_{B} \tk
\bigl(g_i(x)\bigr) D_{g_i}(x) \,dx = \int
_{B} \ka_i(x) \,dx.
\end{equation}

Given $a >0$ and $x \in
U_i$, by $g'_i(x)(\H_{a})$ we mean the point process in $\R^d$
obtained by applying to $\H_{a}$ the linear map $g'_i(x)$.
Similarly, $g'_i(x)(z)$ is the image of $z$ under the map
$g'_i(x).$
If $\M\in\mathbb{M}(m,d)$ or $\M$ is an open subset of $\R^m$,
and $f\dvtx \M\to\R$ is measurable,
then we say $w \in\M$ is a \emph{Lebesgue point} of $f$ if
$\epsilon^{-m} \int_{B_{\epsilon}(w) \cap\M} |f(y) - f(w)| \,dy$
tends to zero as $\epsilon\downarrow0$.

%
%le5.1 #&#
\begin{lemm}
\label{Penlem1}
Suppose $i \in\I$, $U_i$ is bounded,
and $u \in U_i$ is a Lebesgue point of $\ka_i$.
Suppose $\ell(n), n \in\N$ is a sequence of integers
with $\ell(n) \sim n$
as $n \to\infty$.
Set $y_0 = g_i(u)$.
Then
as
$n \to\infty$
we have
(in the above sense
of convergence of point processes in $\R^d$)
%
%e5.5 #&#
\begin{equation}\label{0506c}
n^{1/m}(-y_0 + \Y_{\ell(n)}) \stackrel{{\mathcal
D}} {\longrightarrow }g'_i(u) (\H_{\ka_i(u)}).
\end{equation}
\end{lemm}

\begin{pf}
By taking $B =U_i$ in (\ref{0506a}), we see that $
\int_{U_i} \ka_i(x) \,dx \leq1$, so that $\ka_i$ is a (possibly)
defective density function on $U_i$. Extend $\ka_i$ in an arbitrary
manner to a probability density function on $\R^m$.\eject

Let $\X_{n}$ be a point process in $\R^m$ consisting of
$\ell(n)$ independent identically distributed random $m$-vectors
$X_{i,1},\ldots,X_{i,\ell(n)}$
with density $\ka_i$. Then
%
%e5.6 #&#
\begin{equation}\label{0506b}
g_i^{-1}\bigl( \Y_{\ell(n)} \cap
g_i(U_i) \bigr) \stackrel {{\mathcal D}} {=}
\X_{n} \cap U_i
\end{equation}
because for Borel $B \subseteq U_i$, (\ref{0506a})
shows that $ P [ Y_1 \in g_i(B) ] = P[X_{i,1}\in B]. $

By Lemma 3.2 of \cite{PeBer} (restated as Lemma 3.2 of
\cite{PeEJP}), the distribution of the point process $n^{1/m} ( - u
+ \X_{i,n})$ converges weakly to that of $\H_{\ka_i(u)}$,
in the metric
of \cite{PeBer}, which is \emph{not} the same as the one we are using
here (as discussed in \cite{PeBer}, the metric in \cite{PeBer} is
complete but not separable). We claim that if $f\dvtx \R^m \to\R$ is
measurable with bounded support, then
%
%e5.7 #&#
\begin{equation}\label{MP1006a}
\qquad\lim_{n \to\infty} \E\exp\biggl( \sqrt{-1} \sum_{x \in n^{1/m} (-u + \X_{i,n} )}
f(x) \biggr) = \E \exp\biggl( \sqrt{-1} \sum_{ x \in\H_{\ka_i(u)} }
f(x) \biggr).
\end{equation}
Indeed, by
the proof of Lemma 3.2 of \cite{PeBer} there is a coupling in which
the random variables under the expectations on the left and right-hand
sides of (\ref{MP1006a}) are equal with probability tending to
1.

It follows from (\ref{MP1006a}) that for any finite collection of bounded
Borel sets $A_j$ $(1 \leq j \leq k)$ in $\R^m$, the joint
distributions of the variables $\card(n^{1/m} (-u + \X_{n} )
\cap A_j)$, $1 \leq j \leq k$, converge to those of the variables
$\card( \H_{\ka_i(u)} \cap A_j)$, $1 \leq j \leq k$.
Hence, since $U_i$ is a neighborhood of $u$,
the joint
distributions of the variables $\card(n^{1/m} (-u + (\X_{n} \cap
U_i) ) \cap A_j)$, $1 \leq j \leq k$, converge to those of the
variables
$\card( \H_{\ka(u)} \cap A_j)$, $1 \leq j \leq k$.

Therefore the point processes $n^{1/m} (-u + ( \X_{n} \cap U_i))$
converge weakly to $\H_{\ka_i(u)}$, in the sense discussed at the
start of this section; see Theorem 9.1.VI of \cite{DVJ}.

Now we argue as in \cite{AS}, page 251.
By the Skorohod representation theorem, we can
choose coupled point processes $\tX_{n}$ and $\tH_{\ka_i(u)}$, all
on the same probability space, such that $\tX_{n}$ has the same
distribution as $\X_{n} \cap U_i$, and $\tH_{\ka_i(u)}$ has the same
distribution as $\H_{\ka_i(u)}$, and such that $n^{1/m} (-u +
\tX_{n} )$
converges almost surely to $\tH_{\ka_i(u)}$.
That is [see (\ref{ppconv1}) and (\ref{ppconv2})],
we can list the points of
$\tH_{\ka_i(u)}$ as $x_1,x_2,\ldots$ and the points of
$\tX_{n} $ as $x_{n,1}, x_{n,2},\ldots,x_{n,N_n}$, in such a way that
for each $j$ we have almost surely
%
%e5.8 #&#
\begin{equation}
n^{1/m}(x_{n,j} - u) \to x_j \label{MP1006b},
\end{equation}
and for any $L > 0$ with no point of
$\tilde{\mathcal H}_{\ka_i(u)}$ in $B_L(\0)$, we have
almost surely
%
%e5.9 #&#
\begin{equation}\label{ASpp2}
\card\bigl( n^{1/m} (-u + \tX_{n} ) \cap %C^{(m)}_L
B_L(\0)\bigr) \to \card\bigl( \tH_{\ka_i(u)} \cap
B_L(\0) % C^{(m)}_L
\bigr).
\end{equation}
By (\ref{MP1006b}), $x_{n,j} \to u$ as $n \to\infty$, so by
the differentiability of $g_i$, we can write
%
%e5.10 #&#
\begin{equation}
\label{JY2} g_i(x_{n,j}) -g_i(u) =
g'_i(u) (x_{n,j}-u) + w_{n,j},
\end{equation}
with $\|w_{n,j}\|= o(\|x_{n,j}-u\|)$
so $n^{1/m} w_{n,j} \to\0$
as $n \to\infty$.
By
(\ref{MP1006b}) and
(\ref{JY2}),
%
%e5.11 #&#
\begin{equation} \label{MP1006c}
n^{1/m}\bigl(g_i(x_{n,j}) - g_i(
u) \bigr) \to g'_i(u) ( x_j)   \qquad  \mbox{as }  n \to\infty.
\end{equation}\eject\noindent
We claim that
the point process $n^{1/m}(-g_i(u) + g_i(\tX_{n}) )$ converges a.s.
to $g'_i(u) ( \tH_{\ka_i(u)} )$.
The condition
corresponding to (\ref{ppconv1})
follows from (\ref{MP1006c}).
To demonstrate the condition corresponding to (\ref{ppconv2}),
set $y := g_i(u)$, and
let $L$ be such that
no point $g'_i(u)(x_j)$ lies on the boundary of
$B_L(\0)$. We need to show that
%
%e5.12 #&#
\begin{equation}\label{MP1006d}
\lim_{n \to\infty} \card\bigl[ n^{1/m} \bigl(-y + g_i(
\tX_n)\bigr) \cap B_L(\0) \bigr] = \card\bigl[
g'_i(u) (\tH_{\ka_i(u)}) \cap B_L(
\0) \bigr].\hspace*{-35pt}
\end{equation}
Choose $\delta>0$
such that $B_{2 \delta}(y) \subset g_i(U_i)$.
By Lemma \ref{DGlem} we may define finite $K$ by
\[
K := \sup_{z,z' \in B_\delta(y), z \neq z'} \bigl\|g_i^{-1}
\bigl(z' \bigr) - g_i^{-1}(z)\bigr\|/
\bigl\|z' -z\bigr\|.
\]

Let $K' >K$, and suppose $x \in U_i$ with $\|x - u\| > n^{-1/m} K' L
$. By definition of~$K$, if $g(x) \in B_\delta(y)$, then $
n^{1/m}\|g(x) -y\| > L, $ and this also holds if $g(x) \notin
B_\delta(y)$, provided $\delta n^{1/m} > L$.

Hence, for $n$ large the contribution to the left-hand side of
(\ref{MP1006d}) comes only from $x \in\tX_n \cap B_{n^{-1/m} K'
L}(u)$. For large enough $n$, the set of such $x$ consists precisely
of those $x_{n,j}$ such that $x_j \in B_{K'L}(\0)$, provided $K'$ is
chosen so that no point of $\tH_{\ka_i(u)}$ lies on the boundary of
$B_{K'L}(\0)$.

By (\ref{MP1006c}), for large enough $n$ the set of $j$ such that
$x_{n,j}$ contributes to the left-hand side of (\ref{MP1006a}) is precisely
those $j$ such that $g'(u) (x_j) \in B_L(\0)$. Thus we have
(\ref{MP1006d}), and therefore $n^{1/m}(-y + g_i(\tX_{n}) )$ converges
almost surely to $g'_i(u)(\tH_{\ka_i(u)})$ as claimed. Hence
$n^{1/m}(-y + g_i(\X_{n} \cap U_i))$ converges in distribution to
$g'_i(u)(\H_{\ka_i(u)})$.
Together with (\ref{0506b}) this yields
%
%e5.13 #&#
\begin{equation}\label{0506d}
n^{1/m} \bigl(-y + \bigl(\Y_{\ell(n)} \cap g_i(
U_i) \bigr) \bigr) \stackrel{{\mathcal D}} {\longrightarrow}g'_i(u)
\bigl( \H'_{\ka_i(u)}\bigr).
\end{equation}
%
%g_i^{-1}( \Y_n \cap g(V_i) ) \eqd\X_{i,n} \cap U_i
By again using Lemma 9.1.VI of \cite{DVJ} (equivalence
of weak convergence and convergence of fidi distributions)
and the fact that $g_i(U_i)$ is a neighborhood of $y$ in~$\M$,
we
can deduce that (\ref{0506d}) still holds with $\Y_{\ell(n)} \cap g_i(U_i)$,
replaced by $\Y_{\ell(n)}$ on the left-hand side, so that (\ref{0506c})
holds as asserted.
\end{pf}

The next lemma is a two-dimensional version of
Lemma \ref{Penlem1}.

%
%le5.2 #&#
\begin{lemm}
\label{Penlem1a} Suppose $1 \in\I$ and $ 2 \in\I$, and $g_1(U_1)
\cap g_2( U_2) = \varnothing$ and also $U_1 \cap U_2 = \varnothing$ and
$U_1 \cup U_2$ is bounded. Suppose that for $i=1,2$, $x_i \in U_i$
is a Lebesgue point of $\ka_i$,
% and of $\kappa_j$ respectively, and
and set $y_i = g_i(x_i)$. Suppose $(\ell(n), n \in\N)$ is a
sequence of positive integers such that $\ell(n) \sim n$ as $n \to
\infty$. Then as $n \to\infty$,
%
%e5.14 #&#
\begin{eqnarray}\label{0508c}
&&\bigl[n^{1/m}(-y_1 + \Y_{\ell(n)} ),
n^{1/m}(-y_2 + \Y_{\ell(n)})\bigr]
\nonumber
\\[-8pt]
\\[-8pt]
\nonumber
&&\qquad \stackrel{{\mathcal
D}} {\longrightarrow} \bigl[g_1'(x_1) (
\H_{\ka_1(x_1)}),g_2'(x_2) (
\tH_{\ka_2(x_2)})\bigr], %\eqd\H_{\tk(y_1)},
\end{eqnarray}
where $\tH_a$ here is an independent copy of
$\H_a$.
\end{lemm}
\begin{pf} By (\ref{0506a}) we have $\int_{U_1} \ka_1(x) \,dx +
\int_{U_2} \ka_2(x) \,dx \leq1$. Since we assume $U_1 \cap U_2 =
\varnothing$ and $U_1 \cup U_2$ is bounded, we can therefore find a
probability density function $\ka$ on $\R^m$ which is an extension
of both $\ka_1$ and $\ka_2$, that is, with $\ka(x) = \ka_i(x)$ for $x
\in U_i$, $i \in\{1,2\}$.

Let $\X_{n}$ be a point process in $\R^m$ consisting of $n$ independent
identically distributed random $m$-vectors $X_{1},\ldots,X_{n}$
with density $\kappa$. Then
%
%e5.15 #&#
\begin{equation}\label{0508a}
\qquad\bigl(g_1^{-1}\bigl( \Y_{\ell(n)} \cap
g_1(U_1)\bigr), g_2^{-1}\bigl(
\Y_n \cap g_2(U_2) \bigr) \bigr) \stackrel{{
\mathcal D}} {=}( \X_{{\ell(n)}} \cap U_1, \X_{\ell(n)} \cap
U_2)
\end{equation}
because for $i =1,2$ and Borel $B \subseteq U_i$,
(\ref{0506a}) shows that $P[Y_1 \in g_i(B)] =  P[X_1 \in B]$.

By Lemma 3.2 of \cite{PeBer} (restated as Lemma 3.2 of
\cite{PeEJP}), the joint distribution of the point processes
$n^{1/m} ( - x_1 + \X_{{\ell(n)}})$,
$n^{1/m} ( - x_2 + \X_{{\ell(n)}})$,
converges weakly to that of
$\H_{\ka_1(x_1)}$,
$\tH_{\ka_2(x_2)}$,
in the metric
of \cite{PeBer}. We can then follow the proof of Lemma \ref{Penlem1}
with straightforward modifications to deduce (\ref{Penlem1a}).
\end{pf}

Before proceeding we shall re-express the Poisson processes
appearing in the limits (\ref{0506c}) and (\ref{0508c}) in a manner that
is intrinsic to $\M$, that is, not dependent on the choice of atlas.
Recalling that $\operatorname{Gr}_m(d)$ is the
Grassmannian, given $\M\in\mathbb{M}$ and
$y \in\M$, let $T_y\M\in \operatorname{Gr}_m(d)$ be the
hyperplane tangent to $-y +\M$ at $\0$, that is, the image of $\R^m$ under
the linear map $g'_i(x)$ when $x \in\R^m$
and
$(U_i,g_i)$ is any chart such
that $x \in U_i$ and $y = g_i(x)$. We normalize the Lebesgue measure
on $T_y\M$ (with volume element denoted $du$) in such a way that
for any orthonormal basis $(f_i)_{i=1}^m$ of the subspace, the set
$\{\sum_{i=1}^m a_i f_i\dvtx 0 \leq a_i \leq1\}$ has unit Lebesgue
measure.
For $y,z \in\M$, let $\H'_{y, \kappa(y)}$ denote a
homogeneous Poisson point process on ${T}_y\M$
with intensity $\tk(y)$, and let
$\tilde{\H}'_{z, \kappa(z)}$ denote a
homogeneous Poisson point process on ${T}_z\M$
with intensity $\tk(z)$, independent of $\H'_{y,\kappa(y)}$.

%
%le5.3 #&#
\begin{lemm}
\label{intrinswklem}
Suppose $z_1 \in\M$ and $z_2 \in\M$ are distinct Lebesgue points
for $\tk$. Suppose $(\ell(n),n \geq1)$ is a sequence of
positive integers with $\ell(n) \sim n$ as $n \to\infty$.
Then
%
%e5.16 #&#
\begin{equation} \label{0508e}
\qquad\bigl[n^{1/m}(-z_1 + \Y_{\ell(n)} ),
n^{1/m}(-z_2 + \Y_{\ell(n)})\bigr] \stackrel{{\mathcal
D}} {\longrightarrow} \bigl(\H'_{z_1,\kappa(z_1)},
\tH'_{z_2, \kappa(z_2)}\bigr).
\end{equation}
\end{lemm}
\begin{pf}
It is easy to see that we can choose our
atlas $(U_i,g_i)_{i \in\I}$
such that $z_i \in g_i(U_i)$ for $i=1,2$, and such that moreover
$g_1(U_1) \cap g_2(U_2) = \varnothing$, and $U_1 \cap U_2 = \varnothing$,
and $U_1 \cup U_2$ is bounded. Let $x_i:= g_i^{-1}(z_i)$ for $i =
1, 2.$ Then by Lemma \ref{Penlem1a}, to prove (\ref{0508e}), it
suffices to
demonstrate for $i=1,2$
the distributional equality
%
%e5.17 #&#
\begin{equation}\label{0508b}
g'_i(x_i) (\H_{\ka_i(x_i)} )
\stackrel{{\mathcal D}} {=}\H'_{z_i,\kappa(z_i)}.
\end{equation}

Let $i=1$ or $i=2$. By the mapping theorem on page
18 of \cite{Ki}, $g_i'(x_i) (\H_{\ka_i(x_i)}) $ is a Poisson
process on the linear space $g'_i(x_i)(\R^m)$ with intensity
measure $\mu$ where $\mu(B) $ is $\ka_i(x_i)$ times
$|(g'_i(x_i))^{-1}(B)|$, where $|\cdot|$ denotes the $m$-dimensional
Lebesgue measure.

Recall from Section \ref{sectermin} the definition of $D_{g_i}(x)$.
For bounded measurable $A\subset\R^m$,
it is a fact from linear algebra
that
%
%e5.18 #&#
\begin{equation}\label{0508d}
\bigl|g'_i(x_i) (A)\bigr| = D_{g_i}(x_i)|A|.
\end{equation}
Indeed, the columns of
the Jacobian matrix $J_{g_i}(x_i)$ are the images under $g'(x_i)$ of
the standard basis vectors of $\R^m$, so (\ref{0508d})
clearly holds when the standard basis vectors map to an orthonormal
system, but then it can be deduced in the general case using
standard properties of determinants. Equation (\ref{0508d}) is the
basis of the formula (\ref{int-form}) given earlier.

By (\ref{0508d}), if $A = (g'_i(x_i))^{-1}(B)$, then $|B| =
D_{g_i}(x_i) |A|$ so that $\mu(B)=\ka_i(x_i) |B|/D_{g_i}(x_i)$ so by
(\ref{kappi}), $\mu(B) = \kappa(z_i)|B|$ and
(\ref{0508b}) follows.
\end{pf}

Next we give weak convergence results for $\xi$. % on the point
Recall that $\0$ is the origin of $\R^m$.
The next lemma is an analog of Lemma
3.6 of \cite{PeEJP} and Lemma 3.6 of \cite{PeBer}.

%
%le5.4 #&#
\begin{lemm} \label{Penlem2}
Suppose $\xi$ is continuous in the sense of Definition
\ref{hstab}, and rotation invariant. Let $y \in\M$ and $z \in\M$
be a
pair of
distinct Lebesgue points for $\tk$ with $\tk(y)>0$ and
$\tk(z) >0$.
Suppose $(\ell(n),n \geq1)$ is a sequence of positive integers
such that $\ell(n) \sim n $ as $n \to\infty$. Let $k \in\Z^+$ and
$\rho\in(0,\infty]$. Then as $n \to\infty$ we have
%
%e5.19 #&#
\begin{equation} \label{MPlem2eq1}
\bigl[\xi_{n,k,\rho}(y , \Y_{\ell(n)}), \xi_{n,k,\rho}(z ,
\Y_{\ell(n)}) \bigr] \stackrel{{\mathcal D}} {\longrightarrow} \bigl[\xi(\0,
\H_{\tk(y)}), \xi(\0, \tH_{\tk(z)})\bigr].
\end{equation}
Also, if we choose a chart
$(U,g)$ and $u \in U$ such
that $y = g(u)$,
then for almost all (fixed) $x\in\R^m$,
setting $v_n:= u + n^{-1/m}x$, we have
as $n \to\infty$ that
%
%e5.20 #&#
\begin{eqnarray}\label{MPlem2eq2}
&&\xi_{n,k,\rho}\bigl(y, \Y_{\ell(n)}^{g(v_n)} \bigr)
\xi_{n,k,\rho}\bigl(g (v_n ), \Y_{\ell(n)}^{y}
\bigr)
\nonumber
\\[-8pt]
\\[-8pt]
\nonumber
&&\qquad \stackrel{{\mathcal D}} {\longrightarrow} \xi\bigl(\0,
\H_{y,\tk(y)}^{' g'(u)(x) } \bigr) \xi\bigl(g'(u) (x),
\H_{y,\tk(y)}^{'\0} \bigr).
\end{eqnarray}
\end{lemm}
\begin{pf}
First suppose $\rho= \infty$.
Then the left-hand side of (\ref{MPlem2eq1})
is equal to
\[
\bigl[\xi \bigl(\0, n^{1/m}(-y + \Y_{\ell(n)}) \bigr), \xi \bigl(
\0, n^{1/m}(-z + \Y_{\ell(n)}) \bigr) \bigr].
\]

Also, by rotation invariance, the right-hand side of (\ref{MPlem2eq1})
has the same distribution as $[\xi(\0,\H'_{y,
\tk(y)}),\xi(\0,\tH'_{z, \tk(z)})]$. Therefore,\vspace*{1pt}
under the assumptions given, (\ref{MPlem2eq1}) is
immediate from Lemma \ref{intrinswklem} and the continuous mapping
theorem (\cite{Bill}, Chapter 1, Theorem 5.1).

Next we prove (\ref{MPlem2eq2}) in the case $\rho=\infty$.
By (\ref{scalexi}) and translation invariance of~$\xi$,
%
%e5.21 #&#
\begin{eqnarray} \label{1007a}
&&\xi_n\bigl(y, \Y_{\ell(n)}^ { g( v_n)} \bigr)
\xi_n\bigl(g(v_n), \Y_{\ell(n)}^{y}\bigr)\nonumber\\
&&\qquad=\xi\bigl( \0, n^{1/m} ( -y + \Y_{\ell(n)} ) \cup\bigl\{
n^{1/m} \bigl(-y + g(v_n) \bigr) \bigr\} \bigr)
\\
&&\qquad\quad{}\times \xi\bigl(n^{1/m}\bigl(-y + g(v_n) \bigr) ,
n^{1/m}(-y + \Y_{\ell(n)}) \cup\{\0\}\bigr)
\nonumber\\
&&\qquad=  F\bigl( n^{1/m} \bigl(-y + g(v_n) \bigr),
n^{1/m}(-y + \Y_{\ell(n)} ) \bigr),\nonumber
\end{eqnarray}
where for any $y \in\R^d$ and any locally finite
$\yy\subset\R^d$, we set
\[
F(y,\yy) : = \xi \bigl( \0, \yy\cup\{y\} \bigr) \xi \bigl( y, \yy \cup\{\0\}
\bigr) = \xi \bigl( \0, \yy\cup\{y\} \bigr) \xi \bigl( \0, -y + \bigl(\yy \cup\{\0
\} \bigr) \bigr) .
\]
By definition, $n^{1/m} (v_n -u) = x$, and by the same argument
as for (\ref{MP1006c}) earlier on,
$ n^{1/m} (g(v_n) -y)
$
converges to $g'(u) (x)$ as $n \to\infty$.
Combining this with Lemma~\ref{Penlem1}, we
have by (\ref{0508b}) the convergence in distribution
\[
\bigl( n^{1/m} \bigl(g(v_n) -y \bigr), n^{1/m}[-y
+ \Y_{\ell(n)}] \bigr) \stackrel{{\mathcal D}} {\longrightarrow}
\bigl(g'(u) ( x), \H'_{y,\tk(y)} \bigr).
\]
By the continuity assumption, for almost every $x$ the point set $
\H'_{y,\tk(y)} \cup\{ g'(u)(x)\} $
is a.s.
a continuity point
of $\xi(\0,\cdot)$, and so is the point set\break
$-g'(u)(x) + (\H'_{y,\tk(y)} \cup\{\0\}) $.
Thus
$(g'(u)( x), \H'_{y,\tk(y)} )$ is a.s. at
a continuity point of $F$, for almost all~$x$. Hence we have the
desired convergence in distribution (\ref{MPlem2eq2})
(when $\rho=\infty$) by~(\ref{1007a}) and
the continuous mapping theorem.

Finally, we consider the case with $0 < \rho< \infty$. It is
easy to see that
\[
\lim_{n \to\infty} P \bigl[ \xi_{n,k,\rho}(y , \Y_{\ell(n)})
\xi_{n,k,\rho}(z , \Y_{\ell(n)}) \neq\xi_{n}(y ,
\Y_{\ell(n)}) \xi_{n}(z , \Y_{\ell(n)}) \bigr] = 0,
\]
so the general case of (\ref{MPlem2eq1}) follows
from the special case with $\rho=\infty$ (already proved) along
with Slutsky's theorem. The proof of the general case of
(\ref{MPlem2eq2}) is similar.
\end{pf}

%s6 #&#
\section{\texorpdfstring{Proofs of Theorems \protect\ref{smainLLN}--\protect\ref{PoCLT}}
{Proofs of Theorems 3.1--3.3}}
\label{sec6}

We first give some definitions. Assume $\M\in\mathbb{M}$ and $\tka
\in
\mathbb{P}(\M)$ are given, and set $\K:= \K(\tka)$. We adapt to the
manifold setting the definition of exponentially stabilizing
functionals \cite{BY2,PeEJP}. Suppose $k \in
\Z^+ ,r \in[0,\infty)$ are given, along with the density $\tka$.
For $y \in\K$ and locally finite $\Y\subset\K$,
define
\[
R_\la(y,\Y) := \cases{ \max \bigl[r,\NND_k \bigl(
\lambda^{1/m}y,\lambda^{1/m}\Y \bigr) \bigr],      &\quad $\mbox{if }
\card \bigl(\Y\setminus\{y\} \bigr) \geq k$, \vspace*{2pt}
\cr
\la^{1/m}
\diam(\K), &\quad  $\mbox{otherwise}$. }
\]
Thus
if $k =0$, then
$R_\la(y,\Y) =r$.

It is easy to see that $R:=R_\la(y,\Y)$ serves as a \emph{radius of
stabilization} for any $\xi\in\Xi(k,r)$, in the following sense:
for all
finite $\A\subset(\K\setminus B_{\la^{-1/m}R} (y))$, we have
%
%e6.1 #&#
\begin{equation}
\label{0107a} \xi_\la\bigl(y, \bigl( \Y\cap B_{\la^{-1/m} R}(y)
\bigr) \cup\A\bigr) = \xi_\la\bigl(y, \Y\cap B_{\la^{-1/m} R}(y)
\bigr).
\end{equation}
For all $k \in\Z^+, \rho\in(0, \infty)$, note that $R$ also
serves as a radius of stabilization for $\xi_{\la,k,\rho}$ in the
sense that (\ref{0107a}) holds if $\xi_\la$ is replaced by
$\xi_{\la,k,\rho}$. Recall the definition of point processes $\Y_n$
and $\P_\la$ in Section \ref{sectermin},
and recall that $\S_2$ is the collection of all subsets of
$\K(\kappa)$ of cardinality at most $2$, including the empty set.
Given $\eps>0$ and $t >0$, we define the tail probabilities
for $R_\lambda$ denoted $\tau(t)$ and $\tau_\eps(t)$,
for Poisson and binomial input, respectively, as follows:
\begin{eqnarray*}
\tau(t)&:=& \sup_{\la\geq1 } \mathop{\operatorname{ess\,sup}}_{ y
\in\K} P
\bigl[R_\la (y,\Po_\la) > t\bigr] ;
\\
\tau_{\epsilon}(t)&:=& \sup_{\la\geq1, n \in\N
\cap((1- \epsilon) \la, (1 + \epsilon) \la),  {\mathcal A} \in\S
_2 }   \mathop{\operatorname{ess\,sup}}_{y \in\K} P\bigl[R_{\la}(y, \Y_n \cup{
\mathcal A}) > t\bigr],
\end{eqnarray*}
where the $\operatorname{ess\,sup}$ denotes essential
supremum with respect to the measure $\tk(y)\,dy$.

%
%de6.1 #&#
\begin{defn} \label{stab}
Given $k$ and $r$, we say that every
$\xi\in\Xi(k,r)$ is \emph{exponentially stabilizing} for $\tk$
if
$ \limsup_{ t \to\infty}
t^{-1} \log\tau(t) < 0. $
We say that every
$\xi\in\Xi(k,r)$ is \emph{binomially exponentially stabilizing} for
$\tk$
if there exists $\eps>0$ such that
$ \limsup_{ t \to\infty}
t^{-1} \log\tau_{\epsilon}(t) < 0. $
\end{defn}

We next show that functionals in $\Xi(k,r)$ have the continuity
property of Definition \ref{hstab} as well as the binomial and
exponential stabilization
properties. % (recall Definitions \ref{stab} and \ref{binom-stab},

%
%le6.1 #&#
\begin{lemm} \label{contlemm}
Let $k \in\Z^+$ and $r \geq0$. Then every
$\xi\in\Xi(k,r)$
is continuous. If either $k =0$ or
$\tk\in\mathbb{P}_c(\M) $, then
every $\xi\in\Xi(k,r)$ is
exponentially stabilizing and binomially exponentially stabilizing
for $\tk$.
\end{lemm}

\begin{pf} %We prove this for the functional $\xi$ of Theorem
To prove continuity, let $\xi\in\Xi(k,r)$, let $z \in\R^m$, and
let $F\dvtx \R^m \to\R^d$ be linear and of full rank. Assume that the
points of $F(\H^z)$ have distinct Euclidean norms, and that
for all $n \in\N$ there are no points of $F(\H)$ on
the boundary of the ball $B_n(\0)$.
List the elements of
$F(\H)$
in order of increasing Euclidean norm as
as $x_1,x_2,\ldots.$
Suppose $(\yy_n)_{n \in\N}$ is a sequence of locally finite point
sets in
$\R^d$ converging in $\T$ to $F(\H)$, and list
the elements of $\yy_n$ in order of increasing Euclidean
norm as $y_{n,1},y_{n,2}, y_{n,3}, \ldots$ (possibly a terminating
sequence).

Given the realization of $\H$, we pick the smallest $K \in\N$ such
that
\[
K > \max \bigl(r , N_k \bigl(\0,F(\H) \bigr),\bigl\|F(z)\bigr\| \bigr).
\]
Let $N$ denote the number of points of $F(\H)$ in $B_K$, and assume
$(y_1,\ldots,y_N)$ lies at a continuity point of
the mapping $(x_1,\ldots,x_N) \mapsto\xi(\0,\{x_1,\ldots,x_N\})$.
By the convergence of $\yy_n$ to $\H_a$,
for all large enough $n$ we have $y_{n,N} \in B_K$ and
$y_{n,N+1} \notin B_K$, and moreover\vadjust{\goodbreak} $y_{n,j} \to y_j$ as
$n \to\infty$ for all $j \leq N$.
Therefore by the continuity assumption we have
\begin{eqnarray*}
\xi \bigl(\0,F(\H_a) \bigr) &=& \xi \bigl(\0,\{y_1,
\ldots,y_N\} \bigr) = \lim_{n \to\infty} \xi \bigl(\0,
\{y_{n,1},\ldots,y_{n,N}\} \bigr)\\
& =& \lim_{n \to\infty} \xi
\bigl(\0, F(\yy_n) \bigr) .
\end{eqnarray*}
Similarly, if
$(F(z),y_1,\ldots,y_N)$ lies at a continuity point of
the mapping $(x_0,x_1, \ldots,x_N) \mapsto\xi(\0,\{x_0,x_1,\ldots
,x_N\})$
and if
$(\vv_n)_{n \in\N}$ is any sequence of locally finite point sets in
$\R^d$ converging in $\T$ to $F(\H^z)$, then
$
\xi(\0,F(\H_a^z)) =
\lim_{n \to\infty} \xi(\0, \vv_n).
$
Thus, $\xi$ is continuous.

We prove the exponential stabilization of $\xi\in\Xi(k,r)$,
that is, the uniform exponential tail bound for $R_\la(y,\Po_\la)$,
as follows.
Suppose $y \in\K$, $\la\in[1, \infty)$ and
$ r < t \leq\lambda^{1/d}\diam(\K)$. Then
\[
P \bigl[ R_\la(y,\Po_\la) > t \bigr] = P \bigl[
\NND_k \bigl(\la^{1/m}y, \la^{1/m}\P_\la
\bigr) > t \bigr],
\]
and the last event occurs
if and only if
the number of
points from $\P_\la$ in
$B_{t \la^{-1/m}}(y) \setminus\{y\}$
is less than $k$.
The number of such points
is Poisson distributed
with parameter $\alpha(t):=
\alpha(t,y,\la)$ equal to the $\lambda\tk$ measure of
$B_{t \la^{-1/m}} (y) \cap\M$.
By Lemma \ref{claim0} there is a constant $C_2 >0$ such that
we have uniformly in $\la\in[1, \infty)$,
$y \in\K$ and $t \in(0, \la^{1/m} \diam(\K))$ that
$ \alpha(t) \geq C_2^{-1} t^m $.
Thus by a Chernoff bound for the Poisson distribution
(see, e.g., Lemma 1.2 of \cite{Pe}),
there is a constant $C_3 $ such that
for $\max(r,(2kC_2)^{1/m}) < t < \la^{1/m} \diam(\K)$
we have
\[
P \bigl[R_\la(y,\Po_\la) >t \bigr] \leq k \exp
\bigl(-C_3^{-1} t^m \bigr),
\]
and moreover this also holds for $t \geq\la^{1/m} \diam(\K)$ since
$P[R_\la(y,\Po_\la) >t]=0$ in this case. This gives the desired
exponential stabilization of $\xi$ for Poisson input. Modifications
of this argument yields the exponential stabilization of $\xi$ with
respect to binomial input.
\end{pf}

For finite $\Y\subset\R^d$ and $y \in\Y$, and $k \in\Z^+,n \in
\N,
\rho\in
(0,\infty]$, define
%
%e6.2 #&#
\begin{eqnarray}\label{1101b}
\xi^*_{n,k,\rho}(y,\Y) &:=& \xi_{n,k,\rho}(y,\Y)\mathbf{1}\bigl\{ \bigl|
\xi_{n,k,\rho}(y,\Y)\bigr| \leq n^{5/12} \bigr\};
\nonumber
\\[-8pt]
\\[-8pt]
\nonumber
H^{*}_{n,k,\rho}(\Y) &:=& \sum_{y \in\Y}
\xi_{n,k,\rho}^*(y, \Y).
\end{eqnarray}
Recall the (similar) definition of $H_{n,k,\rho}^\xi(\Y)$ at (\ref{rf}).
Given $n,i,\nu\in\N$ with $i \leq\nu$, define
%
%e6.3 #&#
%e6.4 #&#
\begin{eqnarray}
G_{i,\nu,n} &:=& H^\xi_{n,k,\rho}(\Y_\nu) -
H^\xi_{n,k,\rho} \bigl(\Y_\nu\setminus
\{Y_i\}\bigr), \label{1222d}
\\
G^*_{i,\nu,n} &:=& H^{*}_{n,k,\rho}(\Y_\nu) -
H^{*}_{n,k,\rho} \bigl(\Y_\nu\setminus
\{Y_i\}\bigr). \label{101222a}
\end{eqnarray}

%
%le6.2 #&#
\begin{lemm} \label{Lemm52}
Suppose $\xi$ is binomially exponentially stabilizing,
and $\kappa\in\mathbb{P}_b(\M)$.
Suppose $h(n)/n \to0$ as $n \to\infty$ and suppose
for some $p \in\N$
that
(\ref{mom5}) holds. Then
%
%e6.5 #&#
\begin{equation}\label{0921f}
\limsup_{n \to\infty} \sup_{ n - h(n) \leq\nu\leq n + h(n)} \E|G_{\nu,\nu,n}|^p
< \infty
\end{equation}
and
%
%e6.6 #&#
\begin{equation}\label{1222c}
\limsup_{n \to\infty} \sup_{ n - h(n) \leq\nu\leq n + h(n)} \E|G^*_{\nu,\nu,n}|^p
< \infty.
\end{equation}
\end{lemm}
\begin{pf} We prove only (\ref{0921f}); the proof
of (\ref{1222c}) is virtually the same.
Write $\xi_n$ for $\xi_{n,k,\rho}$, and define
$\Delta^x \xi_n (y,\Y) := \xi_n(y,\Y\cup\{x\}) - \xi_n(y,\Y)$.
Putting $Y= Y_\nu$, observe that
%
%e6.7 #&#
\begin{equation}\label{0917b}
|G_{\nu,\nu,n} | \leq\bigl| \xi_n(Y, \Y_{\nu-1})\bigr | + \sum
_{i=1}^{\nu-1} \bigl| \Delta^Y
\xi_n(Y_i,\Y_{\nu-1})\bigr|.
\end{equation}
The $p$th moment of the first term in the right-hand side of (\ref
{0917b}) is uniformly bounded by (\ref{mom5}). The $p$th
moment of
the sum in the right-hand side of (\ref{0917b})
is given by
%
%e6.8 #&#
\begin{eqnarray}\label{0919b}
&&(\nu-1) \E\bigl|\Delta^Y \xi_n(Y_1,
\Y_\nu)\bigr|^p \nonumber\\
&&\qquad{}+ (\nu-1) (\nu-2) \E\bigl|\Delta^Y
\xi_n(Y_1,\Y_\nu)\bigr|^{p-1} \bigl|
\Delta^Y \xi_n(Y_2,\Y_\nu) \bigr|+
\cdots\\
&&\qquad{}
+ \biggl( \frac{(\nu-1)!}{ (\nu-1-p)!} \biggr) \E\prod_{i=1}^p
\bigl|\Delta^Y \xi_n(Y_i,\Y_\nu)\bigr|.\nonumber
\end{eqnarray}
Let $\I, y_i, \delta_i, U_i$
be as in Section \ref{sectermin}. Using compactness, let $\I_0
\subset\I$ be a finite set such that $\K\subset\bigcup_{i \in
\I_0} B_{\delta_i}(y_i)$, and set $\delta:= \min_{i \in\I_0}(
\delta_i)$.
Then
%
%e6.9 #&#
\begin{eqnarray}\label{0919a}
 &&\nu %\E[ |\Delta^y \xi_n(Y_1,\Y_{\nu-1})|^p | Y =y]
\E\bigl|\Delta^Y \xi_n(Y_1,
\Y_{\nu-1})\bigr|^p
\nonumber
\\[-8pt]
\\[-8pt]
\nonumber
&&\qquad = \int_\M \int
_{g_\ell(U_\ell)} \E \bigl|\Delta^y \xi_n(z,
\Y_{\nu-2} )\bigr|^p \nu\tk(z) \,dz \tka(y) \,dy.
\end{eqnarray}

For $z,y \in\K$ with $\|z-y\| \geq\delta$,
by binomial exponential stabilization we have $P [ \Delta^y
\xi_n(z,\Y_{\nu-2}^z ) \neq0]$ decaying exponentially in
$n^{1/m}$, uniformly over such $(z,y)$
and over $\nu\in[n - h(n), n + h(n)]$. By this, bound
(\ref{mom5}) and H\"older's inequality,
the contribution to
(\ref{0919a}) from such $(z,y)$ tends to zero as $n \to\infty$
and is uniformly bounded.

Now take $y \in\M$, and choose $i \in\I_0$ such that
$y \in B_{\delta_i}(y_i)$. Then $B_{2\delta}(y) \subset
g_i(U_i)$.
Assume without loss of generality that $g_\ell(\0)=y$,
and let $U'_i :=  g_i^{-1}(B_\delta(y))$. Then
the contribution to the inner integral
in the right-hand side of~(\ref{0919a})
from $z \in B_{\delta}(y)$ can be rewritten as the expression
\begin{eqnarray*}
&&\nu\int_{U'_i} \E\bigl| \Delta^y \xi_n
\bigl( g_\ell(u),\Y_{\nu-2} \bigr)\bigr|^p
\ka_\ell(u) \,du %\\
\\
&&\qquad= \int_{\nu^{1/m} U'_i} \E\bigl|
\Delta^y \xi_n\bigl( g_\ell\bigl(
\nu^{-1/m} v\bigr), \Y_{\nu-2} \bigr)\bigr|^p
\ka_i\bigl(\nu^{-1/m}v\bigr) \,dv,
\end{eqnarray*}
and by binomial exponential stabilization, and (\ref{mom5}),
and H\"older's inequality and comparability of norms in $U_i$ and
in $g_i(U_i)$ (see Lemma \ref{DGlem}), there is a constant~$C$,
independent of $y$, such that
the integrand is bounded by
$ C \exp( - C^{-1} \|v\|^{1/m})$, which is integrable in $v$. This
shows that (\ref{0919a}) is uniformly bounded.

Turning to the second term in (\ref{0919b}), we condition on
$Y=y$ with $i \in\I_0$ again chosen
so that
$y \in B_{\delta_i}(y_i) $. By
H\"older's inequality we have
\begin{eqnarray*}
&&\nu^2 \E\bigl|\Delta^y \xi_n(Y_1,
\Y_\nu)\bigr|^{p-1} \bigl|\Delta^y \xi_n(Y_2,
\Y_\nu ) \bigr|
\\
&&\qquad= \nu^2 \int_\K\int_\K
\E\bigl|\Delta^y \xi_n\bigl(w,\Y_{\nu-2}^z
\bigr)\bigr|^{p-1} \bigl|\Delta^y \xi_n\bigl(w,
\Y_{\nu-2}^z \bigr)\bigr|\tka(w) \tka(z) \,dz \,dw
\\
&&\qquad\leq\nu^2 \int_\K\int_\K
\bigl(\E\bigl|\Delta^y \xi_n\bigl(w,\Y_{\nu-2}^z
\bigr)\bigr|^p\bigr)^{(p-1)/p} \\
&&\hspace*{36pt}\qquad\quad{}\times\bigl( \E\bigl| \Delta^y
\xi_n\bigl(z,\Y_{\nu-2}^w \bigr) \bigr|^p
\bigr)^{1/p} \tka(w) \tka(z) \,dz \,dw.
\end{eqnarray*}
The contribution to the last expression from $(w,z) \notin
B_\delta(y) \times B_\delta(y)$
is uniformly bounded (and in
fact tending to zero as $n \to\infty$) by a similar argument to the
one given above for the contribution to (\ref{0919a}) from $w \notin
B_\delta(y)$. Assuming $g_\ell(\0)=y$, the contribution to the last
integral from
$(w,z) \in B_\delta(y) \times B_\delta(y)$
% g_\ell(U_\ell) \times g_\ell(U_\ell)$
is given by
\begin{eqnarray*}
&&\int_{\nu^{1/m} U'_i} \int_{\nu^{1/m} U'_i} \bigl(\E\bigl|
\Delta^y \xi_n\bigl( g_\ell\bigl(
\nu^{-1/m} v\bigr), \Y_{\nu-2} \cup\bigl\{ g_\ell\bigl(
\nu^{-1/m }v\bigr)\bigr\} \bigr)\bigr|^p \bigr)^{(p-1)/p}
\\
&&\hspace*{41pt}\qquad{}\times \bigl(\E\bigl|\Delta^y \xi_n\bigl( g_\ell
\bigl(\nu^{-1/m} v\bigr), \Y_{\nu-2} \cup \bigl\{
g_\ell\bigl(\nu^{-1/m }u\bigr)\bigr\} \bigr)\bigr|^p
\bigr)^{1/p} \\
&&\hspace*{65pt}{}\times\ka_\ell\bigl(\nu^{-1/m}u\bigr) \,du
 \ka_\ell\bigl(\nu^{-1/m}v\bigr) \,dv.
\end{eqnarray*}
By binomial exponential
stabilization, H\"older's inequality and comparability of norms in
$U_\ell$ and in $g_\ell(U_\ell)$ (see Lemma \ref{DGlem}), there is a
constant $C$ such that the integrand is bounded by $ C \exp( -
C^{-1} (\|u\|^{1/m} + \|v\|^{1/m}))$, which is integrable in
$(u,v)$. This  shows that the integral is uniformly bounded, and
hence the second term in (\ref{0919b}) is bounded. The remaining
terms in (\ref{0919b}) are handled similarly, showing that all terms in
(\ref{0919b}) are uniformly bounded, and hence by (\ref{0917b}), we have
(\ref{0921f}).
\end{pf}

\begin{pf*}{Proof of Theorem \ref{smainLLN}} We first sketch the proof
for $q = 2$. Using (\ref{MPlem2eq1}) we obtain, as in the proof of
equations (4.3) and (4.4) of \cite{PeBer}, the distributional
convergence
\[
\xi_{n,k,\rho}(Y_1, \Y_n) \stackrel{{\mathcal D}} {
\longrightarrow }\xi(\0, \H_{\tk(Y_1)}),
\]
where $ \H_{\tk(Y_1)}$ is a Cox process in $\R^d$
whose distribution, conditional on the value~$y$ of $Y_1$, is that
of $\H_{\tk(y)}$, and
also,
\[
\xi_{n,k,\rho}(Y_1, \Y_n) \xi_{n,k,\rho}(Y_2,
\Y_n) \stackrel{{\mathcal D}} {\longrightarrow}\xi(\0,
\H_{\tk
(Y_1)}) \xi(\0, \tH_{\tk(Y_2)}),
\]
where $ \tH_{\tk(Y_2)}$ is an independent copy of the Cox process $
\H_{\tk(Y_1)}$.

Put $\mu:= \E\xi(\0, \H_{\tk(Y_1)})$. Under the assumed moment
condition it follows that
%
%e6.10 #&#
\begin{equation}
\label{covlim} \qquad\lim_{n \to\infty} \E\xi_{n,k,\rho}(Y_1,
\Y_n) \xi_{n,k,\rho}(Y_2, \Y_n) = \E \xi(
\0, \H_{\tk(Y_1)}) \xi(\0, \tH_{\tk(Y_2)}) = \mu^2,
\end{equation}
where the last equality
follows by independence. Recalling definition (\ref{rf})
of $H_{n,k \rho}^\xi(\cdot)$,
we have that
\begin{eqnarray*}
n^{-2 } \E H_{n,k,\rho}^{\xi}(\Y_n)^2
= n^{-1} \E \xi_{n,k,\rho}(Y_1, \Y_n)^2
+ \bigl(1 - n^{-1}\bigr) \E\xi_{n,k,\rho}(Y_1,
\Y_n) \xi_n(Y_2, \Y_n),
\end{eqnarray*}
so we obtain from (\ref{covlim})
that $n^{-2}
\E H_{n,k,\rho}^{\xi}(\Y_n)^2 \to\mu^2$ as $n \to\infty$. Since\break
$n^{-1} \E H_{n,k,\rho}^{\xi}(\Y_n) \to\mu$ as $n \to\infty$ it follows
that $n^{-1} H_{n,k,\rho}^{\xi}(\Y_n)$ converges in $L^2$ to $\mu$.

Since $\mu$ equals the right-hand side
of (\ref{sWLLN}), we have
(\ref{sWLLN})
with $L^2$ convergence
when $q=2$. To obtain $L^1$ convergence when $q =
1$, we use a truncation argument and follow the proof of
Proposition 3.2 in \cite{PY4}. We leave the details to the
reader.

It remains to prove
that if we assume (\ref{mom5}) holds for some $p>5$
and $\tka\in\mathbb{P}_b(\M)$, and that either
$\tka\in\mathbb{P}_c(\M)$ or $k=0$, then
(\ref{sWLLN}) holds with a.s. convergence.
Under these extra assumptions,
Lemmas \ref{contlemm} and \ref{Lemm52} show that
$\E|H^*_{n,k,\rho} (\Y_n) - H^*_{n,k,\rho}(\Y_{n-1})|^5$
is bounded by a constant that is independent of~$n$. Then (2.11) of~\cite{PeBer}
holds
with $\beta=4/3$ and $p' =5$ (and $f\equiv1$ in the notation of~\cite{PeBer}).
By following the proof of Theorem 2.2 of~\cite
{PeBer} we obtain for all $\eps>0$ that
%
%e6.11 #&#
\begin{equation}\label{1004a}
\sum_{n=1}^\infty P\bigl[
\bigl|H^{*}_{n,k,\rho}(\Y_n) - \E H^{*}_{n,k,\rho}(
\Y_n)\bigr | > \eps n\bigr] < \infty.
\end{equation}

Also, by (\ref{smoment0}) [i.e., by taking ${\mathcal A}= \varnothing$ in
(\ref{mom5})] and (\ref{1101b}) and Markov's inequality,
%
%e6.12 #&#
\begin{eqnarray}\label{1004b}
&&\sum_{n=1}^\infty P\bigl[H^{\xi}_{n,k,\rho}(
\Y_n) \neq H^{*}_{n,k,\rho
}(\Y_n) \bigr]\nonumber\\
&&\qquad\leq\sum_{n=1}^\infty nP\bigl[
\xi_{n,k,\rho}(Y_1,\Y_n)\neq\xi^*_{n,k,\rho}(Y_1,
\Y_n)\bigr]
\\
&&\qquad\leq\sum_{n=1}^\infty n^{1-25/12} \E
\bigl[ \bigl|\xi_{n,k,\rho}(Y_1,\Y_n)\bigr|^5
\bigr] < \infty.\nonumber
\end{eqnarray}
By (\ref{1004a}), (\ref{1004b}) and the Borel--Cantelli lemma,
%
%e6.13 #&#
\begin{equation}\label{110101c}
\lim_{n \to\infty} n^{-1} \bigl( H^{\xi}_{n,k,\rho}(
\Y_n) - \E H^{*}_{n,k,\rho}(\Y_n) \bigr) =
0    \qquad \mbox{a.s.}
\end{equation}
Also, $\{ \xi_{n,k,\rho}(Y_1,\Y_n) , n \geq1\}$
are uniformly integrable by (\ref{smoment0}), so
%
%e6.14 #&#
\begin{eqnarray}\label{1004c}
\qquad &&n^{-1} \bigl|\E H^{\xi}_{n,k,\rho}(\Y_n) - \E
H^{*}_{n,k,\rho}(\Y_n)\bigr| \nonumber\\
&&\qquad\leq\E\bigl|\xi_{n,k,\rho}(Y_1,
\Y_n) - \xi^*_{n,k,\rho}(Y_1,\Y_n)\bigr|
\\
&&\qquad= \E\bigl[\bigl|\xi_{n,k,\rho}(Y_1,\Y_n)\bigr | \mathbf{1} \bigl
\{ \bigl|\xi_{n,k,\rho}(Y_1,\Y_n)\bigr | > n^{5/12}
\bigr\} \bigr] \to0   \qquad  \mbox{as }   n \to \infty.\nonumber
\end{eqnarray}
Finally,
by (\ref{sWLLN}) we have that
%the expected value of of
%
\[
\E \bigl[n^{-1} H^\xi_{n,k,\rho}(\Y_n)
\bigr] \to\int_{\M} \E \bigl[ \xi(\0, \H_{\kappa(y)})
\bigr] \kappa(y) \,dy,
\]
and by
(\ref{110101c}) and (\ref{1004c}) we have
(\ref{sWLLN}) with a.s. convergence.
\end{pf*}

The following lemma will be used in the proof of
Theorem \ref{smainCLT}. It can be proved by following
verbatim the proof of Lemma 4.2 of \cite{PeEJP}, so we omit details
here.
Again we write
$\Y^y$ for $\Y\cup\{y\}$.

%
%le6.3 #&#
\begin{lemm}
\label{lemcorr}
Suppose that $\M\in\mathbb{M}(m,d)$ and
$\tka\in\mathbb{P}_b(\M)$, and
$\xi$ is exponentially stabilizing.
Let $k \in\Z^+$, $\rho\in(0,\infty]$ and suppose
for some $p > 2$ that $\xi$ satisfies
(\ref{moment1}).
Then there is a constant $C>0$ such that for all $\la\geq1$ and
all $y,z \in\K$,
\begin{eqnarray*}
&&\bigl| \E\xi_{\la,k,\rho} \bigl(y, \P_\la^z \bigr)
\xi_{\la,k,\rho} \bigl(z, \P_\la^y \bigr) - \E
\xi_{\la,k,\rho} (y,\Po_\la) \E\xi_{\la,k,\rho}(z,
\Po_\la)\bigr| \\
&&\qquad\leq C \exp \bigl( - C^{-1} \la^{1/m}
\|z-y\| \bigr).
\end{eqnarray*}
\end{lemm}

\begin{pf*}{Proof of Theorem \ref{PoCLT}}
We show first the
asymptotic variance convergence~(\ref{sVAR}).
Recalling from (\ref{PPP}) that $\P_\la$ is the Poisson point
process on $\M$ having intensity measure $\la\tk(y)\,dy$,
we have (cf.
the proof of Lemma 4.1 of \cite{PeEJP})
%
%
%e6.15 #&#
\begin{eqnarray}\label{0508f}
&&\la^{-1} \Var \bigl[ H^\xi_{\la, k , \rho} \bigr]\nonumber\\
 &&\qquad= \int
_\M\E \bigl[ \xi_{\la,k,\rho}(y, \P_\la)^2
\bigr] \tk(y) \,dy
\nonumber
\\[-8pt]
\\[-8pt]
\nonumber
&&\qquad\quad{}+ \int_\M\int_{\M} \bigl\{ \E
\xi_{\la,k,\rho} \bigl(z, \P_\la^{y} \bigr)
\xi_\la \bigl(y, \P_\la^z \bigr)\\
&&\hspace*{48pt}\qquad\quad{} - \E
\xi_{\la,k,\rho}(z, \P_\la) \E\xi_{\la,k,\rho} (y,
\P_\la) \bigr\} \la\tk(z) \tk(y) \,dz \,dy.
\nonumber
\end{eqnarray}

Since we assume $\tk\in\mathbb{P}_b(\M)$, we can choose an index
set $\I
\subset\N$, and a set of quadruples
$(y_i,\delta_i,U_i,g_i)_{i \in\I}$ as in
Section \ref{sectermin}, such
that
the support $\tk$ of $\K$ is contained in
finitely many $B_{\delta_i}(y_i), i \in\I_0:= \{1,\ldots,m\} \subset\I$.
Also, we define our partition of unity here by
$\psi_i=
\mathbf{1}\{B_{\delta_i}(y_i) \setminus\bigcup_{j < i} B_{\delta_j}(y_j)
\}.$
Let $\delta:= \min_{i \in\I_0} \delta_i$.

Suppose $y,z \in\K$
with
$\|z - y\| > \delta$. Then by Lemma \ref{lemcorr},
the integrand \mbox{inside} the braces in
the double integral in (\ref{0508f})
is bounded by\break $C \la\exp(-C^{-1} (\delta\la)^{1/m})$, where the
constant $C$ does not depend on $\la,y$ or $z$.
Hence the contribution to the double integral from
such $y,z$ tends to zero.

To estimate the remaining contribution to the double
integral, given $y$ take $i=i(y)$ such that $\psi_i(y)=1$
[so in particular $B_\delta(y) \subset g_i(U_i)$], and
let $U'_i := g_i^{-1}(B_\delta(y))$ and $u := g_i^{-1}(y)$.
Then the contribution to inner integral in the double integral
from $z \in B_\delta(y)$ is given by
\begin{eqnarray*}
&&\int_{U'_i} \bigl\{ \E\xi_{\la,k,\rho}\bigl(y,
\P_\la^{g_i(x)} \bigr) \xi_{\la,k,\rho}\bigl(g_i(x),
\P_\la^{y} \bigr) \\
&&\quad{}- \E\xi_{\la,k,\rho}
\bigl(g_i(x), \P_\la\bigr) \E \xi_{\la,k,\rho}(y,
\P_\la) \bigr\} \la\ka_i(x) \,dx
\\
&&\qquad= \int_{\la^{1/m}(-u + U'_i)} F_\la(v,y) \ka_i
\bigl(u + \la^{-1/m}v\bigr) \,dv,
\end{eqnarray*}
where we set
\begin{eqnarray*}
F_\la(v,y) &:= &\E \xi_{\la,k,\rho}\bigl(y, \P_\la^{g_i(u + \la^{-1/m}v)}
\bigr) \xi_{\la,k,\rho}\bigl(g_i\bigl(u + \la^{-1/m}v
\bigr), \P_\la^{y} \bigr)
\\
&&{}- \E\xi_{\la,k,\rho}(y, \P_\la) \E\xi_\la
\bigl(g_i\bigl(u + \la^{-1/m}v\bigr), \P_\la
\bigr).
\end{eqnarray*}
By the Poisson analog of Lemma \ref{Penlem2}, together with the
moment condition
(\ref{moment1}), provided
$y \in g_i(U_i)$ is a Lebesgue point of $\kappa_i$,
for almost all $v \in\R^m$ we have
as $\la\to\infty$ that
$F_\la(v,y) \to F(g'_i(u)(v),y)$,
where for $x \in{T}_y\M$ we set
\[
F(x,y):= \E\xi\bigl(\0, \H_{\tk(y)}^{x} \bigr) \xi\bigl(x,
\H_{\tk(y)}^{\0} \bigr) - \bigl(\E\xi(\0, \H_{\tk(y)} )
\bigr)^2.
\]
Also,
by (\ref{kappi}) and the assumed a.e. continuity of $\tk$,
almost every $y \in g_i( U_i)$ is a
Lebesgue point of $\ka_i$.
Moreover, by
Lemmas \ref{contlemm}, \ref{lemcorr} and \ref{DGlem},
we have for some constant $C$,
independent of $(\la,y)$, that $|F_\la(v,y)| \leq C \exp(-
C^{-1} \|v\|)$. Hence, using the representation (\ref{int-form}) with
our chosen partition of unity,
by dominated convergence
the double integral in (\ref{0508f})
converges as $\la\to\infty$ to
%
%e6.16 #&#
\begin{equation}
\label{Fint} \sum_{i \in\I} \int_{g(U_i)}
\tk(y) \psi_i(y) \,dy \int_{\R^m} F
\bigl(g'_i(u) (v),y\bigr) \ka_i(u) \,dv.
\end{equation}

We can simplify this limit
by the change of variable $w = g'_i(u) (v)$. Then $dw =
D_g(u) \,dv$ and by (\ref{kappi}) and then (\ref{int-form}) the
expression (\ref{Fint}) equals
%
%e6.17 #&#
\begin{eqnarray} \label{0508g}
&&\sum_{i \in\I} \int_{g(U_i)}
\tk(y)^2 \psi_i (y) \,dy \int_{{\mathcal T}_{y_0}\M}
F(w,y) \,dw
\nonumber
\\[-8pt]
\\[-8pt]
\nonumber
&&\qquad= \int_\M\int_{T_{y}M} F(w,y)
\,dw \tk(y)^2 \,dy,
\end{eqnarray}
so that the double integral at (\ref{0508f}) tends to this.

On the other hand, by
the Poisson analog of Lemma \ref{Penlem2}, the moment bound~(\ref{moment1}), and dominated convergence, the single integral at
(\ref{0508f}) tends to
%
%e6.18 #&#
\begin{equation}\label{3-31a}
\lim_{\la\to\infty} \int_{\M} \E\bigl[
\xi_\la(y, \P_\la)^2 \bigr] \tk(y) \,dy = \int
_\M\E\bigl[ \xi(\0, \H_{\tk(y)})^2
\bigr] \tk(y) \,dy.
\end{equation}
Combining the right-hand sides of (\ref{0508g}) and (\ref{3-31a}),
recalling the definition of $V^\xi$ at (\ref{Vdef}), we get from
(\ref{0508f}) that $\lambda^{-1} \Var[H_{\lambda, k , \rho}^\xi]$
tends to $ \int_{\M} V^\xi(y,\tk(y))
\tk(y) \,dy $, which is (\ref{sVAR}) as desired.

Now assume (\ref{moment1}) holds for some $p >3$.
To prove the normal approximation result (\ref{sSteinbd}), we adapt the
proof of Corollary 2.4 of \cite{PY5} (Corollary 2.1 in the arXiv
version), putting $f \equiv1$. Recall that the partition of unity
has been chosen in such a way that $\psi_i(y) \in\{0,1\}$ for all
$y \in\M$ and all $i \in{\mathcal I}$. For $i \in{\mathcal I}$ let
\[
\K_i := \bigl\{y \in\K\dvtx \psi_i (y) =1 \bigr\}.
\]
As in Sections 4.2 and 4.3 of \cite{PY5}, let
$\rho_\lambda:= \alpha\log\lambda$ for $\lambda>0$,
where $\alpha$ is a suitably chosen large constant; see right after
(4.9) of \cite{PY5}.
For each $\ell\in{\mathcal I}_0$,
cover the bounded set $g_\ell^{-1}(\K_\ell)$ by a (minimal)
collection of
cubes of side $\lambda^{-1/m} \rho_\lambda$,
denoted $Q_{i,\ell}, 1 \leq i \leq V_\ell(\lambda)$,
where $V_\ell= O(\lambda\rho_\lambda^{-m})$.

Fix $\lambda$ for now. For $\ell\in{\mathcal I}_0$ and $1 \leq i
\leq
V_\ell(\lambda)$, let $N_{i,\ell}$ be the number of points of
$\Po_\lambda$ in $\K_\ell\cap g_\ell(Q_{i,\ell})$, a Poisson
variable with parameter $\nu_{i,\ell}$ given by
\[
\nu_{i,\ell} := \lambda\int_{\K_\ell\cap g_\ell(Q_{i,\ell})} \kappa(y)\,dy =
\lambda\int_{g_\ell^{-1}(\K_\ell)
\cap Q_{i,\ell}} \tilde\kappa_\ell(x)\,dx.
\]
Note that the densities $\tilde\kappa_\ell$ are uniformly
bounded because $D_{g_\ell}(\cdot)$ is uniformly bounded on
$g_\ell^{-1}(\K_\ell)$ by compactness.

Let $X_{i,\ell,j}$ denote the $j$th point of $\Po_\lambda\cap
\K_\ell\cap g_\ell(Q_{i,\ell})$, when these points are listed in a
randomized order; cf. Section 4.2 of \cite{PY5}. Then with obvious
modifications, Lemmas 4.2 and 4.3 of \cite{PY5} still hold in the
present setting.

Now follow Section 4.3 of \cite{PY5}, but now defining the graph
$G_\lambda:= ({\mathcal V}_\lambda,{\mathcal E}_\lambda)$ as
follows. The
set ${\mathcal V}_\lambda$ consists of pairs $(i,\ell), 1\leq i \leq
V_\ell(\lambda), \ell\in{\mathcal I}_0$, and the adjacency
${\mathcal
E}_\lambda$ is given by $\{(i,\ell),(j,\ell')\} \in{\mathcal
E}_\lambda
$ if and
only if the distance between $g_\ell(Q_{i,\ell}) \cap\K$ and
$g_{\ell'}(Q_{j,\ell'}) \cap\K$ is at most
$2 \alpha\lambda^{-1/m} \rho_\lambda$.
With $S_{i,\ell}$ defined similarly to $S_i$ in \cite{PY5}, the
variables $(S_{i,\ell}, (i,\ell) \in{\mathcal V}_\lambda)$ have
$G_\lambda$ as a dependency graph.

Next we show that the degrees of the graphs $G_\lambda$ to be bounded by
a constant, uniformly in $\lambda$.
By Lemma \ref{DGlem}, there exists a
constant $K$ such that for all (large enough) $\la$,
%
%e6.19 #&#
\begin{equation}\label{0824g}
\sup_{\ell\in\I_0} \sup_{y,z \in B_{\delta_\ell}(y_\ell) } \frac{
\|g_\ell^{-1}(y) - g_\ell^{-1}(z)\|
}{ \|y-z\|} \leq K
\end{equation}
and also
%
%e6.20 #&#
\begin{equation}\label{110101d}
\sup_{\ell\in\I_0} \sup_{1 \leq i \leq V_\ell(\la)} \frac{ \diam( g_\ell(Q_{i,\ell} ) ) }{\la^{-1/m} \rho_\la} \leq K.
\end{equation}
If $\{(i,\ell),(j,\ell)\} \in{\mathcal E}$, then
$\dist(Q_{i,\ell},Q_{j,\ell}) \leq2 \alpha K
\lambda^{-1/m} \rho_\lambda$, and for any $i$ the number of such
$j$ is bounded by a constant. Here, for subsets $E$ and $F$ of $\R^d$
we put
$\dist(E,F):= \inf\{ \|x - y\|\dvtx x \in E, y \in F \}$.

Now suppose $\ell\neq\ell'$, and fix
$i \leq V_\ell(\lambda)$.
Set $B':= B_{\delta_{\ell'}}(y_\ell')$.
Suppose $\dist(B',\break g_{\ell}(Q_{i,\ell})
\cap\K) > (3 \alpha+K) \lambda^{-1/m} \rho_\lambda$.
Then for each $j \leq V_{\ell'}(\la)$, since
$g_{\ell'}(Q_{j,\ell'}) \cap B' \neq\varnothing$, using
(\ref{110101d})
we have $\dist(g_\ell(Q_{i,\ell}),g_{\ell'}(Q_{j,\ell'}) \geq
3 \alpha\la^{-1/m} \rho_\la$, so
there are no $(j,\ell')$ adjacent to $(i,\ell)$.

Suppose instead that $\dist(B', g_{\ell}(Q_{i,\ell})
\cap\K) \leq(3 \alpha+K) \lambda^{-1/m} \rho_\lambda$.
Choose $w \in B'$
such that $\dist(w, g_\ell(Q_{i,\ell}) \cap\K) \leq
(3 \alpha+K) \lambda^{-1/m}
\rho_\lambda$.
Suppose
$j$ is such that $\{(i,\ell),(j,\ell')\} \in{\mathcal E}$.
Then by the triangle inequality and (\ref{110101d}),
\begin{eqnarray*}
\dist\bigl(g_{\ell'}(Q_{j,\ell'}),w\bigr) &\leq &\dist
\bigl(g_{\ell'}(Q_{j,\ell'}),g_{\ell}(Q_{i,\ell})
\bigr)\\
&&{} + \diam\bigl(g_\ell (Q_{i,\ell})\bigr) + (3 \alpha+K)
\lambda^{-1/m} \rho_\la
\\
&\leq&(5 \alpha+ 2K ) \lambda^{-1/m} \rho_\la
\end{eqnarray*}
so by (\ref{0824g}), $\dist(Q_{j,\ell'},g_\ell^{-1}(w)) \leq
K (5 \alpha+ 2K ) \lambda^{-1/m} \rho_\la$.
Hence
the number of such $j$
is bounded by a constant.

Thus the graphs $G_\lambda$ have degrees bounded uniformly by a
constant independent of $\lambda$, and we can follow the argument in
\cite{PY5} to complete the
proof of (\ref{sSteinbd}).

Finally, if the
condition that (\ref{moment1}) holds for some $p > 3$ is
weakened to (\ref{moment1}) holding for some
$p > 2$, then we may similarly adapt the proof of
Theorem 2.3 of~\cite{PY5} and show
that when $\sigma^2(\xi, \kappa)
> 0$, the left-hand side of (\ref{sSteinbd}) goes to zero, albeit
at a slower rate. That is, in this case (\ref{0324a}) holds.
\end{pf*}

\begin{pf*}{Proof of Theorem \ref{smainCLT}}%\Comment{checked the
We now prove the variance
asymptotics (\ref{sbinomCLTlim}). Recall the definition of $V^\xi(y,a)$
at (\ref{Vdef}),
and note that the definition (\ref{hatsig}) gives
%
%e6.21 #&#
\begin{equation}
\label{newsig} \qquad{\sigma}^2 (\xi, \tk) := \int_{\M}
V^{\xi}\bigl(y, \tk(y) \bigr) \tk(y) \,dy - \biggl( \int
_{\M} \de^{\xi} \bigl(y, \tk(y)\bigr) \tk(y) \,dy
\biggr)^2.
\end{equation}
The idea here is to follow the de-Poissonization argument
in Section 5 of \cite{PeEJP} (with $f \equiv1$).
To ease notation we write $\xi_{\la} $ for
$\xi_{\la,k,\rho}$ and also $\Delta^x \xi_\la(y,\Y)$
for $\xi_\la(y,\Y^x) - \xi_\la(y,\Y)$.
First we seek an analog of Lemma 5.1 of
\cite{PeEJP}. Set
\begin{eqnarray*}
\gamma_1 &:=& \int_\M\E\xi(\0,
\H_{\tk(y)} ) \tk(y) \,dy ;    \\
           \gamma_2&:=& \int
_\M\tk(y)^2 \,dy \int_{{T}_y\M}
\,dv\, \E \bigl[\Delta^v \xi \bigl(\0,\H'_{y,\tk(y)}
\bigr) \bigr]. % \cup\{v\}) - \xi(\0,\H_{\tk(y)}) ].
\end{eqnarray*}
We can show by a similar argument to that used already in the proof
of Theorem~\ref{smainLLN} that the analog of equation (5.6) of
\cite{PeEJP} holds; namely if $\ell\sim\lambda$ and $\tm\sim
\lambda$ with $\ell< \tm$, then
\[
\E\xi_\lambda(Y_{\ell+1}, \Y_{\ell+1} )
\xi_\lambda(Y_{\tm+1}, \Y_{\tm+1} ) \to
\gamma_1^2.
\]
Taking the same partition of unity $\{\psi_i\}$ as
in the proof of
(\ref{sVAR}) above,
analogously to (5.7) of
\cite{PeEJP} we have
%
%e6.22 #&#
\begin{eqnarray} \label{0511a}
&&\ell\E\bigl[ \xi_\la(Y_{\tm+1},\Y_\tm)
\Delta^{Y_{\ell+1}}\xi_\la (Y_1,\Y_\ell)
\bigr] \nonumber\\
&&\qquad= \ell\sum_{i \in\I} \int_\M
\tk(y) \,dy \int_{g_i(U_i)} \psi_i(x) \tk(x) \,dx
\\
&&\qquad\quad{}\times \int_\M\tk(z) \,dz\, \E\bigl[ \xi_\la
\bigl(y,\Y_{\tm-2}\cup\{x,z\}\bigr) \Delta^z
\xi_\la(x,\Y_{\ell-1}) \bigr] .\nonumber
\end{eqnarray}
Using the binomial exponential stabilization and moment conditions,
the contribution to (\ref{0511a}) from
$z \notin g_i(U_i)$
can be shown to vanish as $\la\to\infty$.
By (\ref{kappi}), the remaining contribution to (\ref{0511a})
can be written, using the change of variable
$u = g_i^{-1}(x) $ and $v= g_i^{-1}(z)$,
as
\begin{eqnarray*}
&&\ell\sum_{i \in\I} \int_\M\tk(y)
\,dy \int_{U_i} \psi_i\bigl(g_i(u)
\bigr) \ka_i(u) \,du \int_{U_i}
\ka_i(v) \,dv
\\
&&\qquad{}\times \E \xi_\la\bigl(y,\Y_{\tm-2}\cup\bigl
\{g_i(u),g_i(v)\bigr\}\bigr) \Delta^{g_i(v)}
\xi_\la\bigl(g_i(u),\Y_{\ell-1}\bigr).
\end{eqnarray*}
By the change of
variables $ w = \la^{1/m}(v-u)$, this
equals
%
%e6.23 #&#
\begin{eqnarray}\label{0511c}\qquad
&&\frac{\ell}{\la} \sum_{i \in\I} \int
_\M\tk(y) \,dy \int_{U_i}
\psi_i\bigl(g_i(u)\bigr) \ka_i(u) \,du
\int_{\la^{1/m}(U_i-u)} \ka_i\bigl(u+\la^{-1/m} w
\bigr) \,dw
\nonumber
\\
&&\qquad{}\times \E\bigl[ \xi_\la\bigl(y,\Y_{m-2}\cup\bigl
\{g_i(u),g_i\bigl(u+\la^{-1/m} w\bigr)\bigr\}
\bigr) \Delta^{g_i(u+\la^{-1/m} w)} \\
&&\hspace*{182pt}\qquad{}\times \xi_\la\bigl(g_i(u),
\Y_{\ell-1} \bigr) \bigr],\nonumber
\end{eqnarray}
and by the analog of Lemma 3.7 of \cite{PeEJP} [see also
(\ref{MPlem2eq1}) and (\ref{MPlem2eq2}) of the present paper], along with
the moments conditions and the binomial exponential stabilization to
provide a dominating function, as $\la\to\infty$ with $\ell\sim
\la$ and $\tm\sim\la$ and $\ell< \tm$,
expression (\ref{0511c}), and hence expression (\ref{0511a}),
tend to the expression
%
%e6.24 #&#
\begin{eqnarray}
&&\sum_{i \in\I} \int_\M\tk(y) \,dy
\int_{U_i} \psi_i\bigl(g_i(u)
\bigr) \ka_i(u)^2 \,du\, \E\xi(\0,\H_{\tk(y)})\nonumber\\
&&\quad{}\times  \int
_{\R^m} \,dw\, \E \Delta^{g'_i(u)(w)} \xi\bigl(\0,
\H'_{g_i(u),\tk(g_i(u))}\bigr)
\nonumber\\
&&\qquad= \gamma_1 \sum_{i \in\I} \int
_{U_i} \psi_i\bigl(g_i(u)\bigr)
\ka_i(u)^2 \,du \\
&&\hspace*{36pt}\qquad{}\times \int_{\T_{g_i(u)}\M} \E
\Delta^v \xi\bigl(\0,\H'_{g_i(u),\tk(g_i(u))} \bigr)
\bigl(D_{g_i}(u)\bigr)^{-1} \,dv
\nonumber\\
&&\qquad= \gamma_1 \gamma_2,\nonumber
\end{eqnarray}
which is analogous to (5.13) of \cite{PeEJP}.
Similarly,
$\ell\E\xi_\la(Y_{\ell+1}, \Y_\ell) \Delta^{Y_{\tm+1}} \xi_\la
(Y_1,\Y_\tm)$
converges to $\gamma_1 \gamma_2$,
analogously to equation (5.16) in \cite{PeEJP}.
Moreover, as in
(5.17) of~\cite{PeEJP},
$\E\Delta^{Y_{\ell+1}}\xi(Y_1,\Y_\ell)
\Delta^{Y_{\tm+1}}\xi(Y_2,\Y_\tm)$
is here equal to
\begin{eqnarray*}
&&\la^{-2} \sum_{i \in\I_0} \sum
_{j \in\I_0} \int_{U_i} \int
_{U_j} \psi_i\bigl(g_i(x)\bigr)
\ka_i(x)\,dx \psi_j\bigl(g_j(w)\bigr)
\ka_j(w) \,dw
\\
&&\hspace*{32pt}\qquad{}\times\int_{\la^{-1/d}(-x+ U_i)} \ka_i\bigl(x +
\la^{-1/d} u\bigr) \la \,du \\
&&\hspace*{32pt}\qquad{}\times\int_{\la^{-1/d}(-w+ U_i)}
\ka_j\bigl(w + \la^{-1/d} v\bigr) \la \,dv
\\
&&\hspace*{32pt}\qquad{}\times \E\Delta^{g_i(x+ \la^{-1/d} u)} \xi_\la\bigl(g_i(x),
\Y_{\ell-2}^{g_j(w)}\bigr) \Delta^{ g_j(w+ \la^{-1/d} v)} \xi_\la
\bigl(g_j(w),\\
&&\hspace*{200pt}\Y_{\tm-3}^{g_i(x)} \cup\bigl\{
g_i\bigl(x + \la^{-1/d} u\bigr) \bigr\} \bigr) + o(1)
\end{eqnarray*}
and by the analog of Lemma 3.7 of \cite{PeEJP} [see
also (\ref{MPlem2eq1}) of the present paper],
along with the moments
conditions and the binomial exponential stabilization to provide a dominating
function, as $\la\to\infty$ with $\ell\sim\la$ and $\tm\sim\la$,
this tends to the expression
\begin{eqnarray*}
&&\sum_{i \in\I_0} \sum_{j \in\I_0}
\int_{U_i} \int_{U_j} \psi_i
\bigl(g_i(x)\bigr) \ka_i(x)^2\,dx\,
\psi_j\bigl(g_j(w)\bigr) \ka_j(w)^2
\,dw
\\
&&\quad{}\times\int_{\R^m} \E \Delta^{g'_i(x)(u) } \xi(\0,
\H_{\tk(g_i(x))} ) \,du \times\int_{\R^m} \E
\Delta^{g'_j(w)(v)} \xi\bigl(\0,\H'_{g_i(w),\tk(g_i(w))} \bigr) \,dv
\\
&&\qquad= \biggl( \sum_{i \in\I_0} \int_{g_i(U_i)}
\psi_i(y) \tk(y)^2 \,dy \int_{{T}_y\M} \E
\Delta^z \xi\bigl(\0,\H'_{y,\tk(y)} \bigr) \,dz
\biggr)^2
\\
&&\qquad= \biggl( \int_{\M} \tk(y)^2 \,dy \int
_{{T}_y\M}\E\Delta^z \xi\bigl(\0,
\H'_{y,\tk(y)}\bigr) \,dz \biggr)^2 =
\gamma_2^2.
\end{eqnarray*}
Then by arguments similar to those in the proof
of Lemma 5.1 of \cite{PeEJP}, we have a similar result here with the
squared integral in equation (5.2) of \cite{PeEJP} replaced here by
\[
(\gamma_1 + \gamma_2)^2 = \biggl(\int
_\M\delta^\xi \bigl(y,\tk(y)\bigr) \,dy
\biggr)^2,
\]
and so we obtain
$\lim_{n \to\infty} n^{-1} \Var[H_n^\xi] = {\sigma}^2(\xi, \tk
). $
Given this, by using Theorem~\ref{PoCLT},
following verbatim the proof of Theorem~2.3 of
\cite{PeEJP}, using the case $p=2$ of (\ref{0921f})
(in place of Lemma 5.2 of \cite{PeEJP}),
we can obtain
the desired (\ref{sbinomialVAR})
and~(\ref{sbinomCLTlim}).
\end{pf*}

\begin{rema*}
The method of proof given in this section shows that our
general results can be extended
to a broader class of $\xi$, potentially providing the limit theory for
some of the statistics
mentioned in the penultimate paragraph of Section \ref{INTRO}. This goes as
follows. Consider the class $\F_c$ of
functionals $\xi$ which are continuous in the sense of Definition \ref
{hstab} and which stabilize over
homogeneous Poisson point processes in the sense that (\ref{0107a})
holds when $\Y$ is a homogeneous Poisson
point process and $R_\la(y, \Y) < \infty$ a.s.
Then the above proof of Theorem \ref{smainLLN} shows that the
conclusion of Theorem \ref{smainLLN} holds when
$\Xi(k,r)$ is replaced by $\F_c$. Likewise, if $\F_c(\kappa)
\subset\F_c$ consists of $\xi\in\F_c$
which are exponentially stabilizing and binomially exponentially
stabilizing for $\kappa$, then Theorems~\ref{smainCLT}
and~\ref{PoCLT} hold when $\Xi(k,r)$ is replaced by $\F_c(\kappa)$.
\end{rema*}

%s7 #&#
\section{\texorpdfstring{Proofs of Theorems \protect\ref{dimLLN}--\protect\ref{VRipsthm}}
{Proofs of Theorems 2.1--2.5}} %\label{PROOFS}

%s7.1 #&#
\subsection{\texorpdfstring{Proof of Theorem \protect\ref{dimLLN}}{Proof of Theorem 2.1}} %Recall from

We require some additional lemmas. Observe that
the functional $\zeta_k(y,\Y)$, defined
at (\ref{estim1}), is scale invariant, namely satisfies $\zeta_k(y,
\Y) =
\zeta_k(ay,a\Y)$
for all $a > 0$; cf. remark (ii) in Section \ref{gensec}. Recall that for all $a
\in(0,
\infty)$, $\H_a$ is a homogeneous Poisson point process of intensity
$a$ in $\R^m$.
%
%
%le7.1 #&#
\begin{lemm} \label{elemm}
Let $k >3$. For all $a > 0$ we have $\E\zeta_k(\0, \H_a) = m$ and
$\Var[\zeta_k(\0, \H_a)]= m^2/(k-3).$
\end{lemm}
\begin{pf} If $\NND_j:= \NND_j(\0, \H_a)$, then conditionally on
$\NND_k$, the random variables $(\NND_j/\NND_k)^m, 1 \leq j \leq k
-1$, are distributed as the order statistics of a sample of size $k
- 1$ from the uniform distribution on $[0,1]$, and therefore $-m \log
(\NND_j/\NND_k),  1 \leq j \leq k -1,$ are distributed as the order
statistics of a sample of size $k - 1$ from a standard exponential
distribution. Thus the  sum  $U := m \sum_{j=1}^{k-1} \log
({\NND_k/\NND_j})$ has a  Gamma$(k-1,1)$ distribution and
$\zeta_k(\0, \H_{a}) = (k - 2) m U^{-1}.$ Since
$\E[U^{-1}] = (k-2)^{-1}$ and since $\E[U^{-2}] = ((k-2)(k-3))^{-1}$,
we have $\E\zeta_k(\0,\H_a) = m$ and $\E\zeta^2_k(\0,\H_a) =
m^2(k-2)/(k-3)$, which gives the result.
\end{pf}

For all $k \in\Z^+, \la> 0, \rho> 0$ we put $
\zeta_{\la,k,\rho}(y, \Y):= \zeta_k(\la^{1/m}y, \la^{1/m} \Y
)\times{\mathbf{1}} \{\NND_k(y,\Y) \leq\rho\}$; by scale invariance
$\zeta_{\la,k,\rho}(y,
\Y) = \zeta_{\la,k,\rho}(y, \Y)$. The next lemma shows that
$\zeta_{\la,k,\rho}$ satisfies the moment bounds in the
hypotheses of Theorems \ref{smainLLN}--\ref{smainCLT}.
Recall that for $i \in\Z^+ $, $\S_i$ denotes
all subsets of $\K(\tka)$ of cardinality at most~$i$.

%
%le7.2 #&#
\begin{lemm}
\label{momlem} Let $\M\in\mathbb{M}$ and $\tk\in{\mathbb
{P}}_c(\M).$ With
$\rho_1$ as in Lemma \ref{claim0}, we have for $\rho\in(0,
\rho_1)$, $p > 1$, $k \in\Z^+$
and $i \in\Z^+$
that
%
%e7.1 #&#
\begin{eqnarray}\label{02}
\sup_{y \in\K, n \in\N, \A\in\S_i } \E\bigl[\zeta_{k,\rho}(y, \Y_n \cup
\A)^p\bigr] &<& \infty \qquad\mbox{if } k > p+1 +i;
\\
 \label{04a}
\sup_{y \in\K, \la\geq1, \A\in\S_i } \E\bigl[\zeta_{k,\rho} (y, \P_{\la} \cup
\A)^p \bigr] &<& \infty\qquad \mbox{if } k > p+1 + i.
\end{eqnarray}
\end{lemm}

\begin{pf}
We first show
(\ref{02}). Fix $p > 1$. Let $t_0 \in(2^p, \infty)$ be such that
$\exp(t^{-1/p}) < 1 + 2 t^{-1/p}$ for all $t \in(t_0, \infty)$. %,
Let $x \in\K, n \in\N, \A\in\S_i$.
For $j \in\{1,\ldots,k-1\}$, let
$\NND_j:=\NND_j(y;\Y_n \cup\A)$.
Then
$\zeta_{k,\rho}(y, \Y_n \cup\A) =
(k-2) (\sum_{j=1}^{k-1} \log(\NND_k/\NND_j))^{-1}\mathbf{1}\{\NND_k <
\rho\}$, and
\begin{eqnarray*}
(k-2)^{-p} \E\xi_{k,\rho}(y,\Y_n \cup
\A)^p &=& \E\Biggl[ \Biggl(\sum_{j=1}^{k-1}
\log\frac{\NND_k}{ \NND_j}\Biggr)^{-p}\mathbf{1}\{ \NND_k < \rho\}
\Biggr]
\\
&\leq& t_0 + \int_{t_0}^{\infty} P \Biggl[
\Biggl(\sum_{j=1}^{k-1} \log
\frac{\NND_k}{ \NND_j}\Biggr)^{-p}\mathbf{1}\{\NND_k < \rho\} > t
\Biggr] \,dt.
\end{eqnarray*}
Let $P_{N_k}$ denote the probability distribution of $N_k$.
Conditioning on $\NND_k$ we obtain
%
%e7.3 #&#
\begin{eqnarray}\label{0103e}
&&(k-2)^{-p} \E\xi_{k,\rho}(y,\Y_n \cup
\A)^p \nonumber\\
&&\qquad\leq t_0 + \int_{t_0}^{\infty}
\int_0^{\rho} P \Biggl[ \sum
_{j=1}^{k-1} \log \frac{\NND_k}{ \NND_j} <
t^{-1/p} | \NND_k = u \Biggr] \,dP_{\NND_k}(u) \,dt
\\
&& \qquad\leq t_0 + \int_{t_0}^\infty\int
_{0}^\rho P \biggl[ \log\frac{\NND_k}{ \NND_1} <
t^{-1/p} | \NND_k = u \biggr] \,d P_{\NND_k}(u) \,dt,\nonumber
\end{eqnarray}
and by choice of $t_0$,
for all $t \in(t_0, \infty)$,
we have
%
%e7.4 #&#
\begin{eqnarray}\label{JY1}
P \bigl[ %\log{\NND_k\over\NND_1}
\log(\NND_k / \NND_1) <
t^{-1/p} | \NND_k = u \bigr] &\leq& P \bigl[ (
\NND_k/ \NND_1) < 1 + 2t^{-1/p}|
\NND_k = u \bigr]
\nonumber
\\
&=& P \bigl[ \NND_1 > \NND_k/\bigl( 1 + 2t^{-1/p}
\bigr)| \NND_k = u \bigr]
\\
 &\leq& P\bigl[ \NND_1 > \NND_k\bigl( 1 -
2t^{-1/p}\bigr)| \NND_k = u\bigr].\nonumber
\end{eqnarray}

Given $N_k=u$, there are $k-1 - \card(\A\cap B_u)$ points of
$\Y_n $ in the interior of $B_u(y)$, and to have $N_1 >
u(1-2t^{-1/p})$ it is necessary for
these points to all lie in $A_{t,u}(y)$, where we set
$ A_{t,u}(y):= B_u(y) \setminus B_{u( 1 - 2t^{-1/p})}(y).
$
For $u \in(0, \rho]$ and $\rho< \rho_1,$ Lemma \ref{claim0} gives
\[
\int_{A_{t,u}(y) \cap\M}\,dy \leq C_1 \bigl(u^m -
\bigl(u \bigl( 1 - 2t^{-1/p} \bigr) \bigr)^m \bigr) \leq
2^m C_1 mt^{-1/p} u^m,
\]
and $\int_{B^\K_u(y_1)}\,dy \geq
C_0^{-1}u^m.$ %Thus for $u \in(0, \rho_1)$ we have
Thus by the boundedness assumptions on $\tk$,
for all $u \in(0, \rho]$,
%
%e7.5 #&#
\begin{eqnarray}
\label{bound} P\bigl[ \NND_1 > \NND_k\bigl( 1 -
2t^{-1/p}\bigr)| \NND_k = u\bigr] &\leq& \biggl(
\frac{\int_{A_{t,u}(y) \cap\M} \tka(z) \,dz }{
\int_{B_u(y) \cap\M} \tka(z) \,dz } \biggr)^{k-1-i}
\nonumber
\\[-8pt]
\\[-8pt]
\nonumber
 &\leq& C \bigl( t^{-1/p}
\bigr)^{k-1-i}.
\end{eqnarray}
Combining this with (\ref{0103e}) and
(\ref{JY1}) yields
\[
(k-2)^{-p} \E\zeta_{k,\rho}(y,\Y_n \cup
\A)^p \leq t_0 + \int_{t_0}^{\infty}
C \bigl( t^{-1/p} \bigr)^{k-1-i} \,dt,
\]
which is finite and
independent of $n$ if $k > p + 1 +i$. This gives us (\ref{02}),
and the proof of (\ref{04a}) is just the same.
\end{pf}

\begin{pf*}{Proof of Theorem \ref{dimLLN}} First we prove (\ref{close}).
Given $\rho>0$, we have
%
%e7.6 #&#
\begin{eqnarray}\label{1001a}
P\bigl[\hat{m}_{k,\rho} (\Y_n) \neq\hat{m}_{k} (
\Y_n)\bigr] &\leq& P\biggl[ \bigcup_{i \leq
n} \bigl\{
\NND_k(Y_i, \Y_n) \geq\rho\bigr\}\biggr]
\nonumber
\\[-8pt]
\\[-8pt]
\nonumber
 &\leq&
n P\bigl[\NND_k(Y_n, \Y_n) \geq \rho\bigr].
\end{eqnarray}
Now $\NND_k(Y_n, \Y_n) \geq\rho$ if and only if there are at
most $k - 1$ points of $\Y_{n-1} $ in $B_\rho(Y_n)$. By Lemma
\ref{claim0} there exists $p>0$ such that
$P[Y_1 \in B_\rho(y)] \geq p$ for all \mbox{$y \in\K$}.
Hence
by a Chernoff-type large deviations estimate for the binomial distribution
(see, e.g., Lemma 1.1 of \cite{Pe}), provided $(n-1)p \geq2(k-1)$ we
have for some constant $C$ independent of $n$ that
\[
n P \bigl[\NND_k(Y_n, \Y_n) \geq\rho \bigr]
\leq n \exp \bigl(- C^{-1} (n-1)p \bigr),
\]
which is summable in $n$, so that
assertion in (\ref{close}) follows by (\ref{1001a}) and the
Borel--Cantelli lemma.

We now prove the remaining assertions of Theorem \ref{dimLLN} with
$\rho_1$
as in Lemma~\ref{claim0}.
It suffices to show that\vadjust{\goodbreak}
$\zeta_{k,\rho}, k \in\N$, as defined at (\ref{estim1r}), satisfy
the conditions of Theorems \ref{smainLLN}--\ref{smainCLT}.
Since $\zeta_k$ is a continuous function of the $k$ nearest
neighbor
distances, it belongs to the class $\Xi(k,r)$. %Lemma \ref{contlemm}
Lemma
\ref{momlem} establishes that $\zeta_{k,\rho}$ satisfies the moment
condition  (\ref{smoment0}) when  $k > p + 1$. It follows by
Theorem~\ref{smainLLN} with  $q = 2$ there, and  Lemma
\ref{elemm} that
$\lim_{n \to\infty} \hat{m}_{k,\rho} (\Y_n) = m   \mbox{ in }  L^2. $ Combined with (\ref{close}), this implies the stated
convergence in probability (\ref{dimestLLN}) for $\hat{m}_{k}
(\Y_n)$.

If also $k \geq10$, then by taking $i=3$ in (\ref{02}), $\xi_{k,\rho}$
satisfies
(\ref{mom5}) for $p = 5.5$. Hence by Theorem \ref{smainLLN} we have
$\lim_{n \to\infty}
\hat{m}_{k,\rho} (\Y_n) = m $ almost surely. Combined with
(\ref{close}) this gives (\ref{dimestLLN}) with a.s. convergence.

To obtain the limits (\ref{devarlimit}) and (\ref{LBCLT}),
observe that for $k \geq7$,
the moment bound~(\ref{02}) (with $i=3$) is satisfied for $p = 2.5 $,
so taking $\xi\equiv\zeta_k$ we have conditions (\ref{mom5}) and
(\ref{moment1})
with $p=2.5$. Therefore
we can apply Theorem \ref{smainCLT} for this~$\xi$, and since
it is scale
invariant, by
(\ref{binomialedgevar}) and Lemma \ref{elemm}
the limiting variance $\sigma^2(\zeta_k,\kappa)$ appearing
in Theorem \ref{smainCLT} equals the right-hand side
of (\ref{devarlimit}), so (\ref{devarlimit}) and (\ref{LBCLT}) follow.
Finally, Lemma \ref{vpos} completes the proof
of Theorem \ref{dimLLN}.
\end{pf*}

%
%le7.3 #&#
\begin{lemm} \label{vpos}
With $\sigma_k^2$ given by
(\ref{devarlimit}),
it is the case
that
$\sigma_k^2 >0$.
\end{lemm}
\begin{pf}
Since they are given by (\ref{Videnti}) and (\ref{Didenti}),
the expressions $V^{\zeta_k}$ and $\delta^{\zeta_k}$
in~(\ref{devarlimit}) \,depend only on $m$ and not
on $d$ or $\tka$.
Hence by using
the part of Theorem~\ref{dimLLN} that we have already proved,
in the special case
where $d=m$ and $\tka$ is a uniform distribution
on the unit cube in $\R^m$, we have
$
\sigma_k^2 = \lim_{n \to\infty} n \Var[ \hat{m}_{k,\rho}(\X_n)
],
$
where $\X_n$ is a point process consisting of $n$
independent uniformly random vectors in a unit cube
in $\R^m$.
%Strictly speaking,
%the unit cube is not a submanifold-with-boundary in $\R^m$, but the
%geometric estimates we have used hold for the unit cube, so the
%proof of Theorem \ref{smainCLT} carries through to this case.
By the definition (\ref{estim}) and scale invariance of $\zeta_k$,
we have
\[
\sigma_k^2 = \lim_{n \to\infty} n^{-1} \Var
\sum_{X \in\X_n} \zeta_k \bigl(n^{1/m}X,n^{1/m}
\X_n \bigr),
\]
and so it is enough to show that
\begin{eqnarray*}
\liminf_{n \to\infty} n^{-1} \Var\sum_{X \in\X_n}
\zeta_k\bigl(n^{1/m}X,n^{1/m}\X_n\bigr)
>0.
\end{eqnarray*}
This
can be done, either by the method of Penrose and Yukich (Section 5
of \cite{PY1}), or by the method of Avram and Bertsimas (Proposition
5 of \cite{AB}, which can be adapted to binomial input). The
particular functional under consideration here is not considered in
those references, but the general approaches are well known so we
omit further details.
\end{pf}

%s7.2 #&#
\subsection{\texorpdfstring{Proof of Theorems \protect\ref{NNGCLT1}--\protect\ref{VRipsthm}}
{Proof of Theorems 2.2--2.5}}

Recall that $Y_i$ are i.i.d. with density $\tk$ and that $\NND_1(y,
\Y)$ is the Euclidean distance between $y$ and its nearest neighbor
in $\Y$, or $+\infty$ if $\Y\setminus y$\vadjust{\goodbreak} is empty.
To help deal with the possibility that a Poisson process $\Po_\la$
has no elements, define
\[
\tN_1(y,\Y):= \cases{ \NND_1(y,\Y), &\quad $ \Y\setminus y
\neq\varnothing$, \vspace*{2pt}
\cr
0, &\quad  $\Y\setminus y = \varnothing.$}
\]
The proofs of Theorems \ref{NNGCLT1}--\ref{ShCLT} depend
in part on the following lemmas.
%
%
%le7.4 #&#
\begin{lemm} \label{Dbound}
Suppose $\tk\in\mathbb{P}_c(\M)$. There is a constant $C$ such
that for all $n \geq3$, $\lambda\geq1$, $n \geq4$ and
$\ell\in[n/2,3n/2]$,
$y \in\K$, ${\mathcal A} \in\S_3$ and $t \in(0, \infty)$ we have
%
%e7.7 #&#
%e7.8 #&#
\begin{eqnarray} \label{0824a}
P\bigl[\NND_1\bigl(n^{1/m}y, n^{1/m}(
\Y_{\ell}\cup\A)\bigr) > t\bigr] &\leq&\exp\bigl(-C^{-1}t^m
\bigr);
\\
\label{0824b}P\bigl[\tN_1\bigl(\la^{1/m}y, \la^{1/m}
\P_\la\bigr) > t\bigr] &\leq&\exp\bigl(-C^{-1}t^m
\bigr).
\end{eqnarray}
\end{lemm}

\begin{pf} These bounds can be deduced from the proof of Lemma
\ref{contlemm}, but we prefer to argue directly, as follows.
Letting $\alpha(t,y,n)$ be the
$n \tk$ measure of
$B_{tn^{-1/m}}^{\M}(y) $, we have
%
%e7.9 #&#
\begin{eqnarray}
\label{D1bd} P\bigl[\NND_1\bigl(n^{1/m}y, n^{1/m}(
\Y_{\ell} \cup{\mathcal A})\bigr) > t\bigr] &=& \bigl(1 - \alpha(t,y,n)
\bigr)^{\ell}
\nonumber
\\[-8pt]
\\[-8pt]
\nonumber
 &\leq&\exp\bigl( -\ell\alpha(t,y,n)/n\bigr).
\end{eqnarray}
By Lemma
\ref{claim0}, there is a constant $C$ such that
uniformly in $n \geq3$, $y \in\K$, and $t \in
(0, \Delta)$, we have
$\alpha(t,y,n) \geq C^{-1} t^m,$
which gives (\ref{0824a}) for $t < \Delta$, and clearly~(\ref{0824a}) holds for $t \geq\Delta$.
The second bound (\ref{0824b}) is proved
similarly.
\end{pf}

%
%le7.5 #&#
\begin{lemm} \label{intparts1} If $\tk$ is bounded
and $\delta\in(0,m)$, then
$\sup_n \E[\NND_1(n^{1/m}Y_1,\break  n^{1/m}\Y_n)^{-\delta}] < \infty.
$
\end{lemm}

\begin{pf}
Set $F_{n,y}(t):= P[\NND_1(n^{1/m}y, n^{1/m}\Y_{n-1}) \leq t ]$. Then
\[
\E \bigl[\NND_1 \bigl(n^{1/m}Y_1,
n^{1/m}\Y_n \bigr)^{-\delta} \bigr] = \int
_{\M} \int_0^{\infty}
t^{-\delta} \,dF_{n,y}(t) \tk(y) \,dy.
\]
As in (\ref{D1bd}), we have for $n$ large and all $y \in\K$,
$t \in(0,1)$ that
\[
P \bigl[\NND_1 \bigl(n^{1/m}y, n^{1/m}
\Y_{n-1} \bigr) > t \bigr] = \bigl(1 - \alpha(t,y,n)
\bigr)^{n-1} \geq\exp \bigl( -2(n-1)\alpha(t,y,n) \bigr),
\]
where $\alpha(t,y,n) \leq C t^m/n$ by Lemma \ref{claim0}.
Thus for $t \in(0,1)$ we obtain that
%
%e7.10 #&#
\begin{eqnarray}
\label{D1BD} F_{n,y}(t)& =& 1 - P\bigl[\NND_1
\bigl(n^{1/m}y, n^{1/m}\Y_{n-1}\bigr) > t \bigr]
\nonumber
\\[-8pt]
\\[-8pt]
\nonumber
 &\leq&1
- \exp\bigl(-3 C t^m\bigr) \leq3 Ct^m.
\end{eqnarray}

Hence by Fubini's theorem we have for all $y \in\K$ that
\begin{eqnarray*}
\int_0^{\infty} t^{-\delta}
\,dF_{n,y}(t) = \delta\int_0^{\infty}
t^{-\delta-1} F_{n,y}(t)\,dt \leq\delta C \int_0^{1}
t^{m-\delta-1} \,dt + \delta \int_1^{\infty}
t^{-\delta-1} \,dt,
\end{eqnarray*}
which is finite if $\delta
\in(0,m)$. Integrating over $y \in\K$ gives the result.
\end{pf}\eject

\begin{pf*}{Proof of Theorem \ref{NNGCLT1}}
Let $q =1$ or $q = 2$.
If $\alpha>0$ and
$\kappa\in\mathbb{P}_c(\M)$, then
Lemma \ref{Dbound} shows
that $\sup_n \E[\NND_1(n^{1/m}Y_1, n^{1/m}\Y_n)^{\alpha p}] <
\infty$ for all $p
> 0$. If
$\alpha\in(-m/q,0)$ and $\kappa$ is bounded,
and if then $p > q$ is chosen so that $-m < \alpha p < 0$, then Lemma
\ref{intparts1} gives $\E[ \NND_1(n^{1/m}Y_1, n^{1/m}\Y_n)^{\alpha p}]
< \infty$.

Thus, in all these cases
the moment condition (\ref{smoment0}) holds for $\xi= N_1^\alpha$\break
and some $p >q$.
Since $\NND_1^\alpha$
belongs to the class $\Xi(k,r)$, limit (\ref{NNWLLN1})
(with $L^q$ convergence) thus
follows from Theorem \ref{smainLLN}, the fact that $\NND_1^{\alpha}$
is homo\-geneous of order $\alpha$ [see~(\ref{newedgeWLLN})], and the
identity
$\E[\NND_1^{\alpha}(\0, \H)] =\break m \omega_m \int_0^\infty
u^{\alpha+m-1} \exp(- \omega_m u^{m}) \,du
$ which yields
%
%e7.11 #&#
\begin{equation}\label{eqWade}
\E\bigl[\NND_1^{\alpha}(\0, \H)\bigr] = \pi^{-\alpha/2}
\biggl(\Gamma\biggl(1 +\frac{m }{ 2}\biggr) \biggr)^{\alpha/m} \Gamma
\biggl(1 + \frac
{\alpha
}{
m}\biggr) ;
\end{equation}
see also (15) of Wade \cite{Wa}.
%This concludes the proof of $L^q$ convergence.
Moreover,
in the first case
[$\alpha>0$ and $\kappa\in\mathbb{P}_c(\M)$],
Lemma \ref{Dbound} shows that
% then under C1,
the moment condition (\ref{mom5}) holds for $p > 6$, and so by Theorem
\ref{smainLLN}
we
obtain the a.s. convergence at (\ref{NNWLLN1}). This completes the
proof of Theorem \ref{NNGCLT1}.
\end{pf*}

\begin{pf*}{Proof of Theorem \ref{NNGCLT3}}
First assume $\alpha> 0$.
To prove variance
asymptotics and (\ref{binomNNCLT}) it suffices to show that the
functional $\tN_1^\alpha$ satisfies the conditions of Theorem~\ref{smainCLT}.
Since we assume $\tk\in\mathbb{P}_c(\M)$,
Lemma~\ref{Dbound} is applicable, showing that %Since also $\NND_1(y,
the moment conditions (\ref{mom5}) and
(\ref{moment1}) hold when $\xi\equiv\tN_1^\alpha$.
The conditions of Theorem~\ref{smainCLT} are all
satisfied, so that result gives variance asymptotics
(\ref{sbinomialVAR}) and the central limit theorem (\ref{sbinomCLTlim})
for $\xi\equiv\tN_1^\alpha$. Also, this choice of $\xi$ is
homogeneous of order $\alpha$, so that (\ref{binomialedgevar}) is
applicable with $\beta=\alpha$, and applying this identity to~(\ref{sbinomialVAR}) and (\ref{sbinomCLTlim}) gives the results
(\ref{binomNNvar}) and (\ref{binomNNCLT}) for $\alpha> 0$, subject to
showing that $\sigma^2 >0$ in (\ref{binomNNvar}).

Now assume $\alpha\in(-m/2, 0)$. We cannot directly apply Theorem
\ref{smainCLT} because the moment bound (\ref{moment1}) fails since
if $\A= \{z\}$, the distance between $y$ and $z$ can be made
arbitrarily small, and thus $\tN_{1, \lambda}^{\alpha}(y, {\mathcal
P}_{\la} \cup\{z\} )$ can be made arbitrarily large. Instead we
use a truncation argument and follow the approach of \cite{BPY}. Put
$\phi(x) = x^\alpha$ for $x >0$ with
$\phi(0)=0$. Given $\epsilon> 0$ define
the functions
\begin{eqnarray*}
\phi^{(\eps)} (x) := \cases{
\phi(x), &\quad
$\mbox{if } x \geq\eps$,
\vspace*{2pt}\cr
0, & \quad $\mbox{otherwise}$ }
\end{eqnarray*}
and
$\phi_{(\eps)}(x) := \phi(x) - \phi^{(\eps)}(x)$.
Let $\tN_1^{\alpha, \eps}(y, \Y):=
\phi^{(\eps)}(\tN_1(y,\Y))$.

Let $ \eps> 0$. Then the moment bounds (\ref{moment1})
and (\ref{mom5}) (with $\rho=\infty$) hold
for $\xi\equiv\tN_1^{\alpha, \eps}$ for
$p = 3$ (say),
because $\phi^{(\eps)}(y,\Y) \leq\eps^{\alpha}$ for all $y, \Y$.
Thus we may apply Theorem \ref{smainCLT} to deduce that as $n \to
\infty$, $n^{-1} \Var[H_{n}^{\tN_1^{\alpha, \eps}}(\Y_n)]
$ converges to $ {\sigma}^2(\tN_1^{\alpha, \eps}, \kappa)$ and
%
%e7.12 #&#
\begin{equation}\label{N0825d}
n^{-1/2} \bigl(H_{n}^{\tN_1^{\alpha, \eps}}(\Y_n) - \E
H_{n}^{\tN
_1^{\alpha
, \eps}}(\Y_n)\bigr) \stackrel{{\mathcal D}}
{\longrightarrow}\NN\bigl(0, {\sigma}^2\bigl(\tN_1^{\alpha, \eps},
\kappa\bigr)\bigr).
\end{equation}
To complete the proof, we adapt arguments in \cite{BPY1}, given in
detail in \cite{BPY}, as follows. To make the link with \cite{BPY},
%put $\phi(x):= \log(e^\gamma\omega_m x)$ in Section 5 of
for $1 \leq i \leq m$ we set $\NND_{i,m} := \NND_1(Y_i,\Y_m)$,
giving the identification $\tN_{1,n}^{\alpha}(Y_i, \Y_n) = \phi(n
\NND_{i,n}^m)$.

The equivalent of Lemma 5.1 of \cite{BPY} is given by Lemma
\ref{Nlemtrunc} below. Let us assume this for now. Lemma 5.2 of
\cite{BPY} remains valid here modulo some small notational
modifications; the ${\mathcal H}$ featuring in that result should be
considered now as a homogeneous Poisson process in $\R^m$. Likewise,
Lemma 5.3 of \cite{BPY} carries over with straightforward
modifications. By the proof of Lemma 5.5 of \cite{BPY}, we can show
that
%
%e7.13 #&#
\begin{equation}
\label{Nepsvar} \lim_{\eps\to0}\sigma^2\bigl(
\tN_1^{\alpha, \eps}, \tk\bigr) = \sigma^2\bigl(
\tN_1^{\alpha}, \tk\bigr).
\end{equation}
With (\ref{Nepsvar}) established, the variance asymptotics
(\ref{binomNNvar}) and central limit theorem~(\ref{binomNNCLT})
follow along the lines of the proof of
Theorem 2.1 of \cite{BPY}.

Now assume either $\alpha>0$ or $\alpha\in(-m/2,0)$.
To show positivity of the limiting variance in (\ref{binomNNvar}),
we would like to follow the approach used to show positivity of $\sigma_k^2$
in Lemma \ref{vpos}. In this case we do not have scale invariance. However,
since $\tka$ is a probability density function, for any
$\xi\in\Xi(k,r)$ we have from (\ref{hatsig}) and Jensen's
inequality that
%
%e7.14 #&#
\begin{equation}
\label{JensLB} {\sigma}^2(\xi, \tka) \geq\int_{\M}
\bigl\{ V^{\xi}\bigl(\tka(y) \bigr) - \de^{\xi} \bigl(\tka(y)
\bigr)^2 \bigr\} \tka(y) \,dy ,
\end{equation}
so it suffices to show that when we take $\xi\equiv N_1^\alpha$,
the integrand in the braces in the right-hand side of (\ref{JensLB})
is strictly positive. By what we have already proved,
\[
V^{N_1^\alpha}(a) - \bigl(\de^{N_1^\alpha}(a) \bigr)^2 =
\lim_{n \to\infty} n^{-1} \Var R^\alpha \bigl(n^{1/m}
\X_n \bigr),
\]
where now $\X_n$ is a point process consisting of $n$ independent,
uniformly distributed random vectors in a cube of volume $a^{-1}$ in
$\R^m$. This limit can be shown to be strictly positive by using
the methods of \cite{AB} or \cite{PY1}.

Finally we prove (\ref{NNde}). By (\ref{Didenti}) and (\ref{eqWade})
we have
%
%e7.15 #&#
\begin{equation} \label{0120a}
\delta^{N_1^\alpha} = \pi^{-\alpha/2} \biggl(\Gamma\biggl(1 +
\frac{m }{
2}\biggr) \biggr)^{1/m} \Gamma\biggl(1 +
\frac{\alpha}{ m}\biggr) + \int_{\R^m} f(x) \,dx,
\end{equation}
where for $x \in\R^m$ we here set
$f(x) := \E N_1^\alpha(\0,\H^x) - \E N_1^\alpha(\0,\H)$, so that
\begin{eqnarray*}
f(x) = \int_{\R^m \setminus B_{\|x\|}(\0)} \bigl(\|x\|^\alpha- \|y
\|^\alpha\bigr) \exp\bigl(- \omega_m \|y\|^m
\bigr) \,dy
\end{eqnarray*}
and hence by Fubini's theorem,
\begin{eqnarray*}
\int_{\R^m} f(x) \,dx& =& \int_{\R^m} \,dy \exp
\bigl(- \omega_m \|y\|^m\bigr) \int_{B_{\|y\|}(\0)}
\bigl(\|x\|^\alpha- \|y\|^\alpha\bigr) \,dx
\\
&=& \biggl( \frac{- \alpha\omega_m^2 m}{\alpha+m} \biggr) \int_0^\infty
r^{\alpha+ 2m -1} \exp\bigl(-\omega_m r^m\bigr) \,dr
\\
&=& - \biggl( \frac{ \alpha\omega_m^{-\alpha/m} }{\alpha+m} \biggr) \Gamma(2+ \alpha/m) = -(\alpha/m)
\omega_m^{-\alpha/m} \Gamma(1+\alpha/m).
\end{eqnarray*}
Substituting back into (\ref{0120a}) gives us (\ref{NNde}),
completing the
proof of Theorem \ref{NNGCLT3}.
\end{pf*}

We state the equivalent of Lemma 5.1 of \cite{BPY}. Recall
$\NND_{i,m} := \NND_1(Y_i,\Y_m)$.
%
%
%le7.6 #&#
\begin{lemm}
\label{Nlemtrunc}
Suppose $\K\in\mathbb{P}_c(\M)$.
Let
$\phi_{(\eps)} (\cdot)$ be as in
the proof of Theorem~\ref{NNGCLT3}.
Given $\delta>0$ there exists $\eps_0 >0$ and $n_0
>0$ such that for $\eps\in(0,\eps_0)$ and $n \geq n_0$ we have
$\Var\sum_{i=1}^n \phi_{(\eps)} (n\NND_{i,n}^m) \leq\delta$.
\end{lemm}
\begin{pf} We may proceed as in \cite{BPY} with only minor
modifications as far as~(5.4) of \cite{BPY}, but more effort is
required to adapt the proof of
(5.19) of \cite{BPY}. To do this,
set $Y = Y_{n+1}$ (corresponding to $X$ in \cite{BPY}), and for $1
\leq j \leq J$
define the open cones
$W_j$ with vertex $Y$, just as in \cite{BPY} (these cones cover
$\R^d$, not just the tangent hyperplane at $Y$).

Let $I_{j,n} := \mathbf{1}\{\Y_n \cap W_j(Y) \neq\varnothing\}$. If
$I_{j,n}=1$, then set $Z_{j,n}$ to be the nearest neighbor of $Y$ in
$\Y_n \cap W_j(Y)$; otherwise set $Z_{j,n} :=Y$. Set
$R_{j,n} := \|Z_{j,n}-Y\|$.
Noting that
$|\phi_{(\eps)}(\cdot)|$ is nonincreasing on
$(0,\eps)$,
we have as in (5.5) of \cite{BPY} that
%
%e7.16 #&#
\begin{equation}\label{N0820a}
\Biggl\llvert \sum_{i=1}^n \bigl(
\phi_{(\eps)}\bigl(n \NND_{i,n+1}^m\bigr) -
\phi_{(\eps)} \bigl(n \NND_{i,n}^m \bigr)\bigr)
\Biggr\rrvert \leq2 \sum_{j=1}^J
I_{j,n} \bigl |\phi_{(\eps)} \bigl(n R_{j,n}^m
\bigr)\bigr|.
\end{equation}

Using the last inequality in (\ref{bds0}) of the
current paper, we have, similarly to~(5.6) of \cite{BPY}, that
there is a constant $K_6$ such that $P[R_{j,n}^m > r ] \geq(1-K_6
r)^n$, and therefore there is a constant $K_7$ such that for $0 < t
\leq1$ and large enough $n$ we have
\[
P \bigl[ n R_{j,n}^m >t \bigr] \geq(1- K_6
t/n)^n \geq\exp(- K_7 t).
\]
We may follow the argument given after (5.7) of \cite{BPY} to obtain
the analog to (5.9) of \cite{BPY} (with $\NND_{i,n}^m$ instead of
the $D_{i,n,k}^g$ of
\cite{BPY}). We can then continue as in the proof of the case
${\mathcal
K}= \R^d$ of
Lemma 5.1 of \cite{BPY}, to complete the proof.
\end{pf}

Before proving Theorem \ref{ShCLT} we need some preliminary
results.
Recall that $\psi(y, \Y) :=
\log( e^{\gamma} \omega_m \NND_1^m(y, \Y))$.
If $\Y\setminus\{y\} = \varnothing$, let us
set
$\psi(y, \Y) := 0$.
For all $a \in
(0, \infty)$ we claim that
%
%e7.17 #&#
\begin{equation}\label{Shannonlemm}
\E\psi(\0, \H_a) = \gamma+ a \int_0^{\infty}
(\log u) e^{-a u} \,du = -\log a.
\end{equation}
The first equality of (\ref{Shannonlemm}) follows
since the
probability that the volume of the nearest neighbor ball around
$\0$ exceeds $u$ is $e^{-a u}$,\vadjust{\goodbreak}
and the second equality follows since $\int_0^\infty\log(w)
e^{-w} \,dw = - \gamma$; see, for example, \cite{Havil}, page 107.
Also, since $\int_0^\infty(\log w)^2 e^{-w} \,dw = \gamma^2 + \pi^2/6$,
we have
%
%e7.18 #&#
\begin{eqnarray}\label{Varlog}
\E\psi(\0,\H_a)^2 &=& \int_0^\infty(
\gamma+ \log u)^2 ae^{-au} \,du
\nonumber
\\[-2pt]
&=& \gamma^2 +2 (-\gamma- \log a) \gamma+ a \int_0^\infty(
\log u)^2 e^{-a u} \,du
\\[-2pt]
&=& \pi^2/6 + (\log a)^2.\nonumber
\end{eqnarray}

Let $V^\psi$ and $V^\psi(a)$ be given by setting
$\xi\equiv\psi$ in (\ref{Videnti}) and (\ref{Vdef}), respectively.
%
%
%le7.7 #&#
\begin{lemm}
\label{Vpsilem} It is the case that $\sigma^2(\psi,\tka) = V^\psi-
m^{-2} + \Var[\log(\tka(Y_1))]$.
\end{lemm}
\begin{pf}
Using
(\ref{Vdef}),
(\ref{Shannonlemm}) and (\ref{Varlog}), and
setting $v= a^{1/m}u$,
we have that
\begin{eqnarray*}
&&V^\psi(a) - ( \log a)^2 -\pi^2/6
\\[-2pt]
&&\qquad = \int
_{\R^m} \bigl\{ \E\psi\bigl(\0, \H_{a}^{a^{-1/m} v}
\bigr) \psi\bigl(a^{-1/m}v, \H_{a}^{\0}\bigr) -
\bigl(\E\psi(\0, \H_{a})\bigr)^2 \bigr\} \,dv
\end{eqnarray*}
and since $\H_a=a^{-1/m} \H$ in distribution, the last
displayed expression is
\begin{eqnarray*}
= \int_{\R^m} \bigl\{ \E\psi\bigl(\0, a^{-1/m}
\H^{v}\bigr) \psi\bigl(a^{-1/m}v, a^{-1/m}
\H^{\0}\bigr) - \bigl(\E\psi\bigl(\0, a^{-1/m} \H\bigr)
\bigr)^2 \bigr\} \,dv.
\end{eqnarray*}

By definition we always have $\psi(t y,t\Y)= \psi(y,\Y) + m \log
t$, so the preceding display is
\begin{eqnarray*}
&=& \int_{\R^m} \bigl\{ \E\bigl[ \bigl( \psi\bigl(\0,
\H^{v}\bigr) - \log a\bigr) \bigl( \psi\bigl(v, \H^{\0}
\bigr) - \log a\bigr) \bigr] - \bigl(\E\psi(\0, \H)- \log a \bigr)^2
\bigr\} \,dv
\\[-2pt]
&=& V^\psi- \pi^2/6 - 2 (\log a) \delta^\psi,
\end{eqnarray*}
so
%
%e7.19 #&#
\begin{equation} \label{0123a}
V^\psi(a) = (\log a)^2 + V^\psi- 2
\delta^\psi\log a.
\end{equation}
Using (\ref{ddef}) and (\ref{Shannonlemm}), and setting $v= a^{1/m}u$,
we have that
\begin{eqnarray*}
\delta^\psi(a) + \log a &=& \int_{\R^m} \E\bigl[\psi
\bigl( \0, \H_a^{a^{-1/m}v}\bigr) - \psi(\0, \H_a
)\bigr] \,dv
\\[-2pt]
&=& \int_{\R^m} \E\bigl[ \psi\bigl(\0,a^{-1/m}
\H^v\bigr) - \psi\bigl(\0,a^{-1/m} \H\bigr)\bigr] \,dv =
\delta^\psi.
\end{eqnarray*}
Setting
$
I_{1,j}(\tka) := \int_\M(\log\tka(y))^j \tka(y) \,dy
$
for $j=1,2$,
we may use (\ref{hatsig})\break and~(\ref{0123a}) to deduce that
%
%e7.20 #&#
\begin{eqnarray}\label{0121a}
\sigma^2(\psi,\tka) &=& V^\psi+ I_{1,2}(\tka) - 2
\delta^\psi I_{1,1} (\tka) - \bigl(\delta^\psi-
I_{1,1}(\tka)\bigr)^2
\nonumber
\\[-8pt]
\\[-8pt]
\nonumber
&=& V^\psi- \bigl( \delta^\psi\bigr)^2 +
I_{1,2}(\tka) - I_{1,1}^2(\tka) .
\end{eqnarray}
Moreover, by (\ref{Didenti}) and (\ref{Shannonlemm}) we have
$
\delta^\psi= \int_{\R^d} f(x) \,dx,
$
where here we set
\begin{eqnarray*}
f(x) &:=& - m \E\bigl[ \log\bigl(N_1(\0,\H) / N_1\bigl(\0,
\H^x\bigr)\bigr)\bigr]
\\[-2pt]
&=& - \int_0^\infty P\bigl[N_1(0,\H)
> \|x\| e^t \bigr] \,dt = - \int_0^\infty
\exp\bigl( - \omega_m \|x\|^m e^{tm} \bigr) \,dt
\end{eqnarray*}
so that setting $v= \omega_m \|x\|^m$ we have
\[
- \delta^\psi= \int_0^\infty\int
_0^\infty \exp\bigl(- v e^{tm} \bigr)
\,dt \,dv = \int_0^\infty e^{-tm} \,dt =
m^{-1}.
\]
Substituting this into
(\ref{0121a}) gives the result claimed.
\end{pf}

%
%le7.8 #&#
\begin{lemm}
\label{Shanlem} Suppose either \textup{(i)} that $\tka\in\mathbb{P}_b(\M)$,
%is bounded with bounded support,
or \textup{(ii)} that $m=d$ and $\M= \R^d$ and $r_c(\tka) >
0$. Then there exists $\delta>0$ such that $\sup_n
\E[\NND_1(n^{1/m} Y_1,\break  n^{1/m} \Y_n)^\delta] < \infty$.
\end{lemm}
\begin{pf}
First consider
case (ii).
% where $\M= \R^m = \R^d$, the
Choose $\delta$ small enough so that $r_c(\tk) > \delta d/(d
-\delta)$. Then by the proof of Theorem 13.3 of \cite{PY6}
(Theorem 2.3 in the arXiv version)
(with
$\delta$ corresponding to $\alpha p$ in that proof), $\E
[\NND_1(n^{1/m} Y_1,\break n^{1/m} \Y_n)^{\delta}]$ is bounded as asserted.

Now assume case (i) instead.
We adapt the proof of Theorem 13.3 of \cite{PY6}
to manifolds. Let $((y_i,\delta_i,U_i,g_i), i \in\I_0)$ be
as in Section \ref{sectermin}, and take
a finite $\I_0 \subset\I$
such that
$\K\subset\bigcup_{ i \in
\I_0} B_{\delta_i}(y_i)$.
Assume $\I_0 = \{1,\ldots,i_0\}$ for some $i_0$,
and write $A_i$ for $B_{\delta_i}(y_i) \setminus\bigcup_{j < i}
B_{\delta_j}
(y_j)$.
For any finite $\Y$ write $L^\delta(\Y)$ for $\sum_{y \in\Y}
\NND_1(y,\Y)^\delta$ with $L^\delta(\Y) = 0$ if $\card(\Y) \leq1$.
Similarly to (3.4) of \cite{PY6}, we have
%
%e7.21 #&#
\begin{equation}\label{0824c}
\E\bigl[ \NND_1\bigl(n^{1/m} Y_1,
n^{1/m} \Y_n\bigr)^\delta\bigr] = n^{\delta/m
-1}
\E\bigl[ L^\delta(\Y_n)\bigr],
\end{equation}
and using the
boundedness of $\K$ we have,
similarly to (3.6) of \cite{PY6},
%
%e7.22 #&#
\begin{equation}\label{0824d}
L^\delta( \Y_n) \leq \sum_{i=1}^{i_0}
L^\delta(\Y_n \cap A_i) + C.
\end{equation}
Similarly to (3.7) of \cite{PY6},
by combining (\ref{0824c}) and (\ref{0824d}) we have
\[
\E \bigl[ \NND_1 \bigl(n^{1/m} Y_1,
n^{1/m} \Y_n \bigr)^\delta \bigr] =
n^{\delta/m -1} \E \Biggl[ \sum_{i=1}^{i_0}
L^\delta(\Y_n \cap A_i) + C \Biggr],
\]
and by Jensen's inequality this remains bounded, provided we can
establish the deterministic bound, for all finite
$\Y\subset A_i $,
%
%e7.23 #&#
\begin{equation}
L^\delta(\Y) \leq C \bigl(\card(\Y)\bigr)^{1 - \delta/m}. \label{0824e}
\end{equation}

By Lemma \ref{DGlem}, there is a constant $C$ such that for all
finite $\Y\subset A_i$ and $x \in\Y$,
taking $y $ to be a nearest neighbor of $x$ in $\Y$,
and also taking $g_i^{-1}(z)$ to be a nearest neighbor of
$g_i^{-1}(x)$ in $g_i^{-1}(\Y)$,
we have
that
\[
\NND_1(x,\Y) = \|y-x\| \leq\|z-x\| \leq C\bigl \| g_i^{-1}(z)
- g_i^{-1}(x)\bigr\| = C \NND_1
\bigl(g_i^{-1}(x), g_i^{-1}(\Y)
\bigr).
\]
Hence, $ L^\delta(\Y) \leq C L^\delta(g_i^{-1}\Y) $, and
using Lemma 3.3 of \cite{Yubk}, we have
(\ref{0824e}) as asserted.
\end{pf}

\begin{pf*}{Proof of Theorem \ref{ShCLT}}
To prove the $L^2$ convergence at (\ref{ShLLN}) we shall apply
Theorem \ref{smainLLN}.
Note that $\psi$ belongs to $\Xi(k,r)$, and
thus it suffices to verify under either of the hypotheses of Theorem
\ref{ShCLT}
that $\xi\equiv\psi$ satisfies the moment condition
(\ref{smoment0}) for $p =3 $.
Set $\NND_1:= \NND_1(n^{1/m}Y_1, n^{1/m}\Y_n)$.
It will suffice to show
that $\sup_n \E| \log\NND_1^m|^p < \infty.$
Given $\delta\in(0,1)$, we can
choose a constant $C$ such that
$|\log t|^p \leq Ct^{-\delta}$ for $t \in[0,1]$ and
$|\log t|^p \leq Ct^{\delta}$ for $t \in[1, \infty).$
Then
%
%e7.24 #&#
\begin{equation}\label{0105b}
\E\bigl| \log\NND_1^m\bigr|^p \leq C \E\bigl[
\NND_1^{-\delta} + \NND_1^\delta\bigr].
\end{equation}

By Lemmas \ref{intparts1} and \ref{Shanlem}, the right-hand side of
(\ref{0105b}) is bounded, provided $\delta$ is chosen small enough.
Hence we can apply Theorem \ref{smainLLN}, yielding
\begin{eqnarray}
\lim_{n \to\infty} n^{-1} S \bigl(n^{1/m} \Y_n
\bigr) = \int_\M\E \bigl[\psi(\0, \H_{\tk(y)})
\bigr] \tk(y) \,dy = -\int_{\M} \log \bigl(\tk(y) \bigr) \tk(y)
\,dy  \nonumber\\[-2pt]
\eqntext{\mbox{in }   L^2,}
\end{eqnarray}
by (\ref{Shannonlemm}). Thus (\ref{ShLLN})
holds.

Now we suppose $\tka\in\mathbb{P}_c(\M)$ and prove the variance
asymptotics (\ref{ShVar}) and central limit theorem
(\ref{binomShCLT}). We cannot directly apply
Theorem \ref{smainCLT}
because the moment bound
(\ref{moment1}) fails since if $\A= \{z\}$, the distance between
$y$ and $z$ can be made arbitrarily small. % and thus $\psi_{\la}(y,
As in the proof of the case $\alpha<0$ of
Theorem \ref{NNGCLT3}
we use a truncation argument, defining for all
$\epsilon
> 0$ the function
\begin{eqnarray*}
\log^{(\eps)} (x) := \cases{ %
\log(x), &\quad
$\mbox{if } x \geq\eps$,
\vspace*{2pt}\cr
0, & \quad $\mbox{otherwise}.$}
\end{eqnarray*}
Let $\psi^{\eps}$ be defined as $\psi$, with $\log$
replaced by $\log^{(\eps)}$.

As in the proof of Theorem \ref{NNGCLT3} we may show that
the moment bounds (\ref{moment1})
and~(\ref{mom5}) hold (with $\rho=\infty$)
for $\psi^\eps$ for some $p > 2$. The case
${\mathcal A} = \varnothing$ of (\ref{moment1})
follows from the bound
%
%e7.25 #&#
\begin{equation}\label{0819b}
\sup_{\lambda\geq1,  y \in
\K} \E\bigl| \psi\bigl(\lambda^{1/d}y,\lambda^{1/d}
\P_\la\bigr) \bigr|^p < \infty,
\end{equation}
which is proved similar to
(\ref{smoment0}) which we have already established.
For the case $\A= \{z\}$, observe that
%
%e7.26 #&#
\begin{eqnarray}\label{0819a}
&&\psi^\eps\bigl(\lambda^{1/m}y, \lambda^{1/m}\bigl(
\Po_\lambda\cup\{z\}\bigr) \bigr) \nonumber\\
&&\qquad= \psi^\eps\bigl(
\lambda^{1/m}y,\lambda^{1/m} \Po_\lambda\bigr) \mathbf{1}
\bigl\{ \NND_1(y,\Po_\la) \leq\|y-z\|\bigr\}
\\
&&\qquad\quad{}+ \log^{(\eps)}\bigl(e^\gamma\omega_m \lambda\|y-z
\|^m\bigr) \mathbf{1}\bigl\{ \NND_1(y,\Po_\la) >
\|y-z\| \bigr\}.\nonumber
\end{eqnarray}
The first term on the right-hand side of (\ref{0819a}) has bounded
$p$th moment, which may be proved similarly to the proof of
(\ref{smoment0}). Since the function $|\log(\cdot)|$ is decreasing on
$(0,1)$ and increasing on $(1,\infty)$, provided $\eps<1$ the
absolute value of the last term in (\ref{0819a}) is bounded by $|\log
(\epsilon)|$ if $e^\gamma\omega_m \lambda\|y-z\|^m < 1$, and
otherwise is bounded
by $|\psi^\eps(\lambda^{1/m}y,\lambda^{1/m} \Po_\la)|$, which
has bounded $p$th moment, as noted already. This gives us
(\ref{moment1}), and the argument for (\ref{mom5}) is similar.

By Theorem \ref{smainCLT}, as $n \to\infty$, we deduce that
$n^{-1} \Var[H_{n}^{\psi^\eps}(\Y_n)] $ converges to $
{\sigma}^2(\psi^\eps, \kappa)$ and
%
%e7.27 #&#
\begin{equation} \label{0825d}
n^{-1/2} \bigl(H_{n}^{\psi^\eps}(\Y_n) - \E
H_{n}^{\psi^\eps}(\Y_n)\bigr) \stackrel{{\mathcal D}}
{\longrightarrow}\NN\bigl(0, {\sigma}^2\bigl(\psi^\eps,
\kappa\bigr)\bigr).
\end{equation}
To complete the proof,
we adapt arguments in \cite{BPY1},
given in detail in \cite{BPY}, as follows.
To make the link with \cite{BPY}, put
$\phi(x):= \log(e^\gamma\omega_m x)$
in Section 5 of \cite{BPY},
and also for $1 \leq i \leq m$, set $\NND_{i,m} := \NND_1(Y_i,\Y_m)$,
giving the identification
$\psi_n(Y_i, \Y_n) = \phi(n \NND_{i,n}^m)$.
Set $\log_{(\eps)}(x) = \log(x) - \log^{(\eps)}(x)$ and
$\phi_\eps(x):= \log_{(\eps)}(e^\gamma x)$. It is easily seen with
this choice of $\phi$ that Lemma \ref{Nlemtrunc} holds [follow that
proof verbatim, choosing $\epsilon$ small enough so that
$|\phi_\eps(\cdot)|$ is nonincreasing on $(0,\eps)$]. It now
suffices to follow the proof of Theorem \ref{NNGCLT3} for the case
$\alpha\in(-m/2, 0)$.

It remains only to show $\sigma^2(\psi, \tk) >
0$. We again follow the approach used in the
proof of positivity in Theorem \ref{NNGCLT3}.
By (\ref{ShVar}), it is enough to show
that the expression
$
V^{\psi}
- ( \delta^{\psi}
)^2,
$
is strictly positive.
This can be shown to be nonnegative as in
the proof of positivity in Theorem \ref{NNGCLT3},
and we
leave the details to the reader.
\end{pf*}

\begin{pf*}{Proof of Theorem \ref{VRipsthm}}
Let $\VRdelta\in(0,\infty)$.
Let $\varphi_k(y, \Y):=\varphi_k^{(\VRdelta)}(y, \Y)$
be $(k + 1)^{-1}$ times the number of $k$-simplices containing $y$
in ${\mathcal R}^\VRdelta(\Y)$,
that is, $(k+1)^{-1}$ times the number of
unordered $(k+1)$-tuples of points in $\Y$, all pairwise within
$\VRdelta$ of each other, and including $y$.
Then
%
%e7.28 #&#
\begin{equation}\label{0825a}
\Cl_k^{(\VRdelta)}(\Y) = \sum_{y \in\Y}
\varphi_k(y, \Y).
\end{equation}
We want to show that $\varphi_k$ satisfies the conditions of Theorems
\ref{smainLLN} and \ref{smainCLT}. Note that $\varphi_k
\in\Xi(0,r)$ if $r > \beta$. For any $y \in\K$, $\A\in\S_3$ and
$k,\ell\in\N$, $\varphi_k(n^{1/m}y,
n^{1/m}( \Y_\ell\cup\A) ) $ is bounded
by $(3 +\Cl_\VRdelta(n))^k$, where we set
$\Cl_\VRdelta(n):= \sum_{i=1}^\ell\mathbf{1}\{\|Y_i - Y_1\| < \VRdelta
n^{-1/m}\}$,
which is
binomially distributed with parameters $\ell$
and the $\tk$ measure of $B_{\VRdelta n^{-1/m}}(y)$.
Assuming $\tka$ is bounded,
there is a constant $C$ such that
the latter parameter is at most
$C \ell^{-1}$, uniformly in $y$.
Hence $\Cl_\VRdelta(n)$ is stochastically bounded by a binomial
random variable with parameters $\ell$ and $C \ell^{-1}$,
and thus for any $p \geq1$, $\varphi_k$ satisfies the moment condition
(\ref{mom5}), and hence also (\ref{smoment0}).
%Also, by Lemma \ref{contlemm}, $\varphi_k$ is continuous.
Therefore by using
Theorem \ref{smainLLN} and
(\ref{0825a}),
we have %\Comment{$H_n$ or $H_{n, k, \infty}$? JY.
%
%e7.29 #&#
\begin{eqnarray}\label{0825e}
\lim_{n \to\infty} n^{-1} \Cl_k^{\VRdelta}
\bigl(n^{1/m} \Y_n\bigr) &=& \lim_{n \to\infty}
n^{-1} H_n^{\varphi_k}(\Y_n)
\nonumber
\\[-8pt]
\\[-8pt]
\nonumber
&=& \int
_{\M} \E\varphi_k^{(\VRdelta)}(\0,
\H_{\kappa(y)}) \kappa(y) \,dy,
\end{eqnarray}
with both $L^2$ and a.s. convergence.

Define $h_k^{(\beta)}\dvtx (\R^m)^{k+1} \to\R$ by
$h_k^{(\beta)}(x_1,\ldots,x_{k+1}): = \prod_{1 \leq i < j \leq k+1}
\mathbf{1} \{ \|x_i-x_j\| \leq\beta\}$, that is,
the indicator
of the event that $x_1,\ldots,x_{k+1}$ are all within distance
$\beta$ of each other.
By the Palm theory of Poisson processes
(e.g., Theorem 1.6 of \cite{Pe})
we have
\begin{eqnarray*}
\E\varphi_k^{(\VRdelta)} \bigl(\0, \H'_\lambda
\bigr)& =& \frac{ \lambda^{k}}{ (k+1)!} \int_{\R^m} \cdots\int
_{\R^m} h_k^{(\beta)} (\0,x_1,
\ldots,x_k) \,dx_1 \cdots \,dx_k
\nonumber\\
&=& \frac{ (\lambda\VRdelta^m )^{k}}{ (k+1)!} \int_{\R^m} \cdots \int
_{\R^m} h_k^{(1)} (\0,x_1,
\ldots,x_k) \,dx_1 \cdots \,dx_k \\
&=&
\frac{ (\lambda\VRdelta^m )^{k}}{ (k+1)!} J_{k,k+1},\nonumber
\end{eqnarray*}
where the last equality comes from (\ref{Jdef}).
Combined with (\ref{0825e}), this gives us~(\ref{CliqLLN}).

A slight modification of
the above argument shows that $\varphi_k$ satisfies the Poisson moment
condition (\ref{moment1})
for any $p \geq1$.
By Theorem \ref{smainCLT},
the variance asymptotics~(\ref{binomcliqvar})
and central limit theorem (\ref{binomcliqCLT}) follow, with
$\sigma_k^2(\beta,\tka) := \sigma^2(\varphi_k^{(\beta)},\tka) $
given by (\ref{hatsig}). We need to show that this is consistent with
(\ref{1230a}).

The expression
$\delta^{\varphi_k}(a)$, given by (\ref{ddef}),
simplifies further as
%
%e7.30 #&#
\begin{equation}\label{0825b}
\delta^{\varphi_k}(a) = \bigl( a \beta^m\bigr)^k
J_{k,k+1} + a \int_{\R^m} \E\bigl[
\Delta^u \varphi_k(\0, \H_{a}) \bigr] \,du.
\end{equation}
Using the Palm theory
of the Poisson process again, we have for all $u \in\R^m$ that
\begin{eqnarray*}
&&\E\bigl[\Delta^u \varphi_k\bigl(\0,
\H'_{a} \bigr)\bigr] \\
&&\qquad= \frac{ a^{k-1}}{(k+1)(k-1)!} \int
_{\R^m} \cdots\int_{\R^m}
h_k^{(\beta)} (\0,u,x_1,\ldots,x_k)
\,dx_1 \cdots \,dx_{k-1}.
\end{eqnarray*}
Together with
(\ref{0825b}), this gives
%
%e7.31 #&#
\begin{eqnarray}\label{0121b}
\qquad\delta^{\varphi_k}(a) &=& \bigl(a\beta^m\bigr)^k
J_{k,k+1} + \frac{a^k k}{(k+1)!} \int\cdots\int h^{(\beta)}(
\0,x_1,\ldots,x_k) \,dx_1 \cdots
\,dx_k
\nonumber
\\[-8pt]
\\[-8pt]
\nonumber
&=& (k+1) \bigl(a \beta^m\bigr)^k J_{k,k+1}.
\end{eqnarray}
For the first
term in (\ref{hatsig}) we simplify the expression $V^{\varphi_k}(a)$ as
follows. Consider the special case where $m=d$ and $\M$ is a
smoothly bounded region of volume $a^{-1}$, and let $\Po'_\la$ be a
homogeneous Poisson point process in this region with expected total
number of points equal to $\lambda$. By applying (\ref{sVAR}) in this
case, recalling notation $\Cl_k^{(\beta)}$ from Section
\ref{subsecVR} we get that
\begin{eqnarray*}
V^{\varphi_k}(a)& =& \lim_{\la\to
\infty} \lambda^{-1} \Var\sum
_{y \in\Po'_\la} \varphi_k\bigl(
\la^{1/m}y, \la^{1/m} \Po'_\la\bigr) =
\lim_{\la\to\infty} \lambda^{-1} \Var\bigl[ \Cl_k^{(\beta)}
\bigl(\lambda^{1/m} \Po'_\la\bigr) \bigr]
\\
&=& \lim_{\la\to\infty} \lambda^{-1} \Var \bigl[ \Cl_k^{(\beta(a/\lambda)^{1/m})}
\bigl(a^{1/m} \Po'_\la\bigr) \bigr].
\end{eqnarray*}
Hence by Proposition 3.7 of \cite{Pe},
setting $k'=k+1$ we have
\[
V^{\varphi_k}(a) = \sum_{j=1}^{k'}
J_{k,j} \bigl(a \beta^m \bigr)^{2 k' -j -1 }.
\]
Using this
identity in the first term of (\ref{hatsig}),
and using
(\ref{0121b}) for the second term of (\ref{hatsig})
enables us to establish the identity (\ref{1230a}).

It remains to show that $\sigma^2(\varphi_k^{(\beta)},\tk) >
0.$ This can be done as in the proof of
Lemma \ref{NNGCLT3}, that is, using
(\ref{JensLB}) to reduce the problem
to showing positivity
in the case where $d=m$
and $\tka$ is a uniform distribution on a cube,
and using the methods of \cite{AB} or \cite{PY1}
to demonstrate positivity in this case.
This completes the proof of Theorem \ref{VRipsthm}.
\end{pf*}

\section*{Acknowledgments}
We thank Rob Neel for
discussions related to the volume estimates of Section
\ref{balls},
and Nikolai Leonenko for
discussions which stimulated our initial interest in entropy
estimators.

% imsref loaded by akundreckaite, 2013-01-22 10:58:49
% imsref loaded by akundreckaite, 2013-01-22 11:05:28

%suskaldyti doi

\printaddresses

\end{document}